\documentclass[review,hidelinks,onefignum,onetabnum]{siamart220329}

\usepackage{lipsum}
\usepackage{amsfonts}
\usepackage{graphicx}
\usepackage{epstopdf}
\usepackage{algorithmic}
\ifpdf
  \DeclareGraphicsExtensions{.eps,.pdf,.png,.jpg}
\else
  \DeclareGraphicsExtensions{.eps}
\fi

\newsiamremark{remark}{Remark}
\newsiamremark{hypothesis}{Hypothesis}
\crefname{hypothesis}{Hypothesis}{Hypotheses}
\newsiamthm{claim}{Claim}

\headers{Differentially Private Distributed Convex Optimization}{M. Ryu and K. Kim}

\title{Differentially Private \\ Distributed Convex Optimization\thanks{Submitted to the editors DATE.
}
}

\author{  
Minseok Ryu\thanks{Mathematics and Computer Science Division, Argonne National Laboratory
  (\email{mryu@anl.gov}, \email{kimk@anl.gov}).}
\and Kibaek Kim\footnotemark[2]}

\usepackage{amsopn}

\ifpdf
\hypersetup{
  pdftitle={An Example Article},
  pdfauthor={D. Doe, P. T. Frank, and J. E. Smith}
}
\fi

\begin{document}

\maketitle

\begin{abstract}
  This paper considers distributed optimization (DO) where multiple agents cooperate to minimize a global objective function, expressed as a sum of local objective functions, subject to some constraints. 
  In DO, each agent iteratively solves a local optimization model constructed by its own data and communicates some information (e.g., a local solution) with its neighbors in a communication network until a global solution is obtained.
  Even though locally stored data are not shared with other agents, it is still possible to reconstruct the data from the information communicated among agents, which could limit the practical usage of DO in applications with sensitive data. 
  To address this issue, we propose a privacy-preserving DO algorithm for constrained convex optimization models, which provides a statistical guarantee of data privacy, known as differential privacy, and a sequence of iterates that converges to an optimal solution in expectation.  
  The proposed algorithm generalizes a linearized alternating direction method of multipliers by introducing a multiple local updates technique to reduce communication costs and incorporating an objective perturbation method in the local optimization models to compute and communicate randomized feasible local solutions that cannot be utilized to reconstruct the local data, thus preserving data privacy.  
  Under the existence of convex constraints, we show that, while both algorithms provide the same level of data privacy, the objective perturbation used in the proposed algorithm can provide better solutions than does the widely adopted output perturbation method that randomizes the local solutions by adding some noise. 
  We present the details of privacy and convergence analyses and numerically demonstrate the effectiveness of the proposed algorithm by applying it in two different applications, namely, distributed control of power flow and federated learning, where data privacy is of concern.
\end{abstract}



\section{Introduction} \label{sec:introduction}
We consider distributed optimization (DO) where multiple agents cooperate to minimize a global objective function, expressed as a sum of local objective functions, subject to some constraints. 
In this paper we focus on the following DO model:
\begin{subequations}
\label{model:dist}
\begin{align}
\min_{w \in \mathbb{R}^n} \ &  \sum_{p=1}^P f_p(w) \\
\mbox{s.t.} \ 
& w \in \mathcal{W} := \cap_{p=1}^P \mathcal{W}_p, 
\end{align}
\end{subequations}
where $P$ is the number of agents, $w \in \mathbb{R}^n$ is a global variable vector, $\mathcal{W}_p$ is a local feasible region assumed to be a compact convex set for each agent $p$, and $f_p$ is a local convex objective function. 
The model \eqref{model:dist} appears in many applications, including distributed control in power systems \cite{molzahn2017survey, patari2021distributed} and distributed machine learning (ML) as well as federated learning (FL) \cite{li2020federated, 9599369}.
Numerous DO algorithms \cite{yang2019survey} have been developed, but most of the algorithms are susceptible to data leakage. 
Specifically, each agent iteratively solves a local optimization model constructed by its own data and communicates some information (e.g., local solutions or local gradients) with its neighbors in a communication network until a global solution $w^*$ is obtained.
Even though locally stored data are not shared with other agents, it is still possible to reconstruct the data from the information communicated among agents, as pointed out in \cite{shokri2017membership, zhu2019deep}.

To address the data privacy concern, we propose a privacy-preserving DO algorithm for solving \eqref{model:dist}.
The proposed algorithm provides a statistical guarantee of data privacy, known as differential privacy (DP), and a sequence of iterates that converges to an optimal solution of \eqref{model:dist} in expectation. 
Under the existence of convex constraints (i.e., $\mathcal{W}_p \subset \mathbb{R}^n$), we show that the proposed algorithm provides more accurate solutions compared with the existing DP algorithm, while both algorithms provide the same level of data privacy.
As a result, our algorithm can mitigate a trade-off between data privacy and solution accuracy (i.e., learning performance in the context of ML), which is one of the main challenges in developing DP algorithms, as described in \cite{dwork2014algorithmic}.

\subsection{Related literature}
DP is a privacy-preserving technique that randomizes an output of a query on data such that any single data sample cannot be inferred by the randomized output \cite{dwork2006calibrating,dwork2014algorithmic}.
Typically, the randomized output is constructed by adding either Laplacian noise or Gaussian noise to the true output of the query, known as Laplace and Gaussian mechanisms, respectively.
Because the true output is perturbed by adding the noise, this DP technique is often referred to as an output perturbation method.
We will discuss details of DP in Section \ref{sec:prelim_dp} since the algorithm proposed in this paper guarantees data privacy in an DP manner.

DP has been integrated in various learning algorithms that train ML models by using data centralized at a single server.
Depending on where to inject noise to randomize an output of the ML models, DP can be categorized by output and objective perturbation methods.
For ease of exposition, we denote the output and objective perturbation by \texttt{OutPert} (i.e., adds a noise vector to the model parameters) and \texttt{ObjPert} (i.e., adds a noisy linear function to the objective function before training), respectively.
These two methods have been applied to convex empirical risk minimization models in \cite{chaudhuri2011differentially, kifer2012private}, where the authors show that both theoretically and empirically, \texttt{ObjPert} is superior to \texttt{OutPert} in managing the inherent trade-off between privacy and learning performance. 

Within the context of distributed ML as well as FL, DP has been introduced in numerous algorithms \cite{agarwal2018cpsgd, naseri2020local, wei2020federated, zhang2016dynamic, huang2019dp, huang2020differentially} and frameworks \cite{ryu2022appfl,he2020fedml,ziller2021pysyft} mainly because of its computational efficiency \cite{wang2020hybrid}.
Specifically, local ML models solved for every iteration of the algorithms are randomized to communicate the randomized outputs (i.e., local model parameters or local gradients depending on the algorithms), which cannot be utilized for reconstructing any single data sample in the local dataset.
Most DP algorithms in the literature rely on \texttt{OutPert} to ensure DP \cite{naseri2020local,wei2020federated,huang2019dp,huang2020differentially}, which adds noise to the intermediate outputs of local ML models.
A few exceptions include \cite{zhang2016dynamic}, where \texttt{ObjPert} is applied to the alternating direction method of multipliers (ADMM) algorithm to randomize the outputs of local models to ensure DP.
For the unconstrained convex ML models, the authors show that ADMM with \texttt{ObjPert} outperforms that with \texttt{OutPert}.
However, the ADMM algorithm requires each agent to compute an exact solution of the convex subproblem for every iteration, which can be computationally expensive.
This computational challenge has been circumvented by employing a linearized ADMM \cite{huang2019dp}, where the convex loss function is linearly approximated, resulting in solving a convex quadratic subproblem that has a closed-form solution.
This work has been extended by introducing multiple local updates to reduce communication costs in \cite{huang2020differentially} and shown to outperform the other existing DP algorithms.

The aforementioned DP algorithms have been developed for unconstrained models; but, in general, constraints are necessary, for example, in most optimal control problems including distributed control of power flow. 
Also, imposing hard constraints on ML models (instead of penalizing constraints in the objective function as in \cite{shen2021agnostic,du2021fairness}) is increasingly considered for purposes such as improving accuracy, explaining decisions suggested by ML models, promoting fairness, and observing some physical laws (e.g., \cite{gallego2022controlled, goh2016satisfying, marquez2017imposing, zafar2019fairness, zhao2019physics}).
This calls for developing DP algorithms suitable for the general constrained optimization model, like the form \eqref{model:dist}.
Indeed the state-of-the-art DP algorithm \cite{huang2020differentially} developed for unconstrained models can be directly applied to solve \eqref{model:dist}, but \texttt{OutPert} used in the algorithm could make each agent  communicate \emph{infeasible} randomized outputs of the constrained local optimization models for every iteration of the algorithm, which in turn could negatively affect the overall convergence of the algorithm and the quality of the final solution.
To address this infeasibility issue arising from the existing DP algorithms applied to the general constrained optimization models, 
we propose a new DP algorithm based on \texttt{ObjPert} for solving \eqref{model:dist}, which allows each agent to communicate \emph{feasible} randomized outputs for every iteration of the algorithm while ensuring data privacy.
Specifically, the proposed algorithm generalizes the existing linearized ADMM by introducing a multiple local updates technique to reduce communication costs and incorporating \texttt{ObjPert} into local constrained optimization models to communicate feasible randomized ouputs that ensure DP. 
To the best of our knowledge, no previous work studies the privacy and convergence analyses of \texttt{ObjPert} in the linearized ADMM with multiple local updates for constrained convex distributed optimization models, which are included in this paper.



 
\subsection{Contributions} 
Our major contribution is to provide  theoretical and numerical analyses of \texttt{ObjPert} in the linearized ADMM with multiple local updates for the constrained model of the form \eqref{model:dist}.
First, under the existence of convex constraints $\mathcal{W}_p$ defined for every agent $p$, we propose how to perturb the local objective function $f_p$ such that the resulting outputs of the local model cannot be utilized for reconstructing data, thus preserving data privacy in an DP manner. 
Specifically, we propose Laplace and Gaussian mechanisms within the context of the linearized ADMM with \texttt{ObjPert}.
It turns out that the proposed mechanisms are different from the existing Laplace and Gaussian mechanisms used in \texttt{OutPert} with respect to the sensitivity computation.
Second, we show that the sequence of iterates produced by our algorithm converges to an optimal solution of \eqref{model:dist} in expectation and that its associate convergence rate is sublinear. 
Specifically, the rate is $O(\frac{1}{\bar{\epsilon} E \sqrt{T}})$ for a smooth convex function, $O(\frac{1}{\bar{\epsilon}^2 E  \sqrt{T}})$ for a nonsmooth convex function, and $O(\frac{1}{\bar{\epsilon}^2 E T})$ for a strongly convex function, where $T$ is the number of ADMM rounds, $E$ is the number of local updates, and $\bar{\epsilon}$ is a privacy budget parameter resulting in stronger data privacy as $\bar{\epsilon}$ decreases.
Third, we numerically demonstrate the effectiveness of our DP algorithm in two different application areas:  distributed control of power flow and federated learning.
 
\subsection{Notation and organization of the paper}
We use $\langle \cdot, \cdot \rangle$ and $\|\cdot \|$ to denote the scalar product and the Euclidean norm, respectively.
We define $[a]:= \{1,\ldots,a\}$.
We use $f'$ and $\nabla f$ to denote a subgradient and a gradient of $f$, respectively.
We use \textbf{pdf}, \textbf{det}, and \textbf{relint} to stand for probability density function, determinant, and relative interior, respectively.

We present some preliminaries in Section \ref{sec:preliminaries}, the proposed algorithm in Section \ref{sec:dpadmm}, privacy analysis in Section \ref{sec:privacy}, convergence analysis in Section \ref{sec:convergence}, numerical demonstration of the proposed algorithm in Section \ref{sec:experiments}, and conclusions in Section \ref{sec:conclusion}. 

\section{Preliminaries} \label{sec:preliminaries}
In this section we present some preliminaries on ADMM and DP for the completeness of this paper. 

\subsection{ADMM} \label{sec:ADMM}
Numerous distributed optimization algorithms \cite{liu2017convergence, nedic2009distributed, johansson2008subgradient, boyd2011distributed, shi2014linear} have been developed for solving \eqref{model:dist}. 
Among them, ADMM and its variants have been widely used for solving distributed optimization. 
Since the proposed algorithm in this paper is built on ADMM, we briefly describe a classical ADMM.
To this end, we first reformulate \eqref{model:dist} by introducing a local copy $z_p \in \mathbb{R}^{n}$ for every agent $p \in [P]$:
\begin{subequations}
\label{model:dist_1}  
\begin{align}
  \min_{w \in \mathbb{R}^n} \ &  \sum_{p=1}^P f_p(z_p)  \\
  \mbox{s.t.} \ 
  & z_p \in \mathcal{W}_p, \ \forall p \in [P], \\
  & w  = z_p, \ \forall p \in [P], \label{dist_1_consensus}
\end{align}
\end{subequations}
where $w$ represents global variables and $z_p$ represents local variables for agent $p$.

A typical iteration $t$ of ADMM for solving \eqref{model:dist_1} is given by
\begin{subequations}
\label{ADMM}
\begin{align}  
w^{t+1} & \leftarrow \arg \min_{w \in \mathbb{R}^{n}} \ \sum_{p=1}^P \{ \langle \lambda^t_p, w \rangle + \frac{\rho^t}{2} \|w-z^t_p\|^2 \}, \label{ADMM-1} \\
z^{t+1}_p & \leftarrow \arg \min_{z_p \in \mathcal{W}_p} \ f_p(z_p) - \langle \lambda^t_p, z_p \rangle + \frac{\rho^t}{2} \|w^{t+1}-z_p\|^2, \ \forall p \in [P], \label{ADMM-2} \\
\lambda^{t+1}_p &  \leftarrow  \lambda^{t}_p + \rho^t (w^{t+1}-z^{t+1}_p), \ \forall p \in [P], \label{ADMM-3}
\end{align}
\end{subequations}
where $\lambda_p$ is the Lagrange multiplier associated with the consensus constraints \eqref{dist_1_consensus} and $\rho^t > 0$ is a penalty parameter.

By noting that \eqref{ADMM-2} can be a computational bottleneck, solving \eqref{ADMM-2} inexactly has been proposed in the literature without affecting overall convergence. For example, linearized ADMM \cite{gao2018information, lin2017extragradient} has been proposed whose convergence rate is proven to be $\mathcal{O}(1/T)$ when the objective function to be linearized is smooth.
In this paper we approximate \eqref{ADMM-2} as follows:
\begin{align}
  z^{t+1}_p \gets \arg \min_{z_p  \in \mathcal{W}_p } \ \langle f'_p(z^t_p), z_p \rangle + \frac{1}{2\eta^{t}}\|z_p - z^t_p\|^2 +  \frac{\rho^t}{2} \|w^{t+1}-z_p + \frac{1}{\rho^t} \lambda^t_p  \|^2,  \label{ADMM-2-Prox} 
\end{align}
which is derived from \eqref{ADMM-2} by (i) replacing the convex function $f_p$ in \eqref{ADMM-2} with its lower approximation $\widehat{f}_p(z_p) := f_p(z^t_p) + \langle f'_p(z^t_p), \ z_p - z_p^t \rangle$, where $f'_p(z^t_p)$ is a subgradient of $f_p$ at $z^t_p$, and (ii) adding a proximal term $\frac{1}{2\eta^{t}}\|z_p - z^t_p\|^2$ with a parameter $\eta^{t} > 0$ that controls the proximity of a new solution $z^{t+1}_p$ from $z^t_p$ computed from the previous iteration. 
We note that when $\mathcal{W}_p := \mathbb{R}^n$ or $\mathcal{W}_p := \{ z_p \in \mathbb{R}^n: l \leq z_p \leq u \}$ for some $l$ and $u$, a closed-form solution expression can be derived from the optimality condition. 


\subsection{Differential privacy} \label{sec:prelim_dp}
In this section we describe DP, a privacy-preserving technique that provides statistical guarantee on data privacy.  
\begin{definition}[Differential privacy \cite{dwork2014algorithmic}] \label{def:differential_privacy_1}
A randomized function $\mathcal{A}$ satisfies $(\bar{\epsilon},\bar{\delta})$-DP if, for all neighboring datasets $\mathcal{D}$ and $\mathcal{D}'$ that differ in a single entry and for all possible events $\mathcal{S}$, the following inequality holds:
\begin{align}
 \mathbb{P}( \mathcal{A}(\mathcal{D}) \in \mathcal{S} ) \leq  e^{\bar{\epsilon}} \mathbb{P}( \mathcal{A}(\mathcal{D}') \in \mathcal{S} ) + \bar{\delta} ,  \label{def_differential_privacy}
\end{align}
where $\bar{\epsilon} \geq 0$, $\bar{\delta} \geq 0$, and $\mathcal{A}(\mathcal{D})$ (resp. $\mathcal{A}(\mathcal{D}')$) represent randomized outputs of $\mathcal{A}$ on input $\mathcal{D}$ (resp. $\mathcal{D}'$). 
Note that in this paper we refer to 
$\bar{\epsilon}$ as the  privacy budget,
$\ln ( \frac{ \mathbb{P}( \mathcal{A}(\mathcal{D}) \in \mathcal{S} ) }{ \mathbb{P}( \mathcal{A}(\mathcal{D}') \in \mathcal{S} )} )$ as the  privacy loss, and a function $\mathcal{A}$ satisfying DP (i.e., \eqref{def_differential_privacy}) as a  mechanism.
\end{definition}

If $\bar{\delta}=0$, we say that $\mathcal{A}$ is $(\bar{\epsilon},0)$-DP, namely, pure DP. 
One can easily see from \eqref{def_differential_privacy} that the absolute value of privacy loss, $|\ln ( \frac{ \mathbb{P}( \mathcal{A}(\mathcal{D}) \in \mathcal{S} ) }{ \mathbb{P}( \mathcal{A}(\mathcal{D}') \in \mathcal{S} )} )|$, is bounded by $\bar{\epsilon}$.
Therefore, as the privacy budget $\bar{\epsilon}$ decreases, it becomes harder to distinguish any neighboring datasets $\mathcal{D}$ and $\mathcal{D}'$ by analyzing the randomized outputs, thus providing stronger data privacy on any single entry in the dataset.
If $\bar{\delta}>0$, we say that $\mathcal{A}$ is $(\bar{\epsilon},\bar{\delta})$-DP, namely, an approximate DP that provides weaker privacy guarantee compared with the pure DP.
According to \cite[Lemma 3.17]{dwork2014algorithmic}, $(\bar{\epsilon},\bar{\delta})$-DP is equivalent to saying that the absolute value of privacy loss is bounded by $\bar{\epsilon}$ with probability $1-\bar{\delta}$. Note that $\bar{\delta} < 1$.

In what follows, we present the existing Laplace and Gaussian mechanisms.
\begin{lemma}[Laplace mechanism \cite{dwork2014algorithmic}]\label{lem:existing-laplace}
Given any function $\mathcal{R}$ that takes data $\mathcal{D}$ as input and outputs a vector in $\mathbb{R}^n$, the Laplace mechanism is defined as
\begin{subequations}
\label{LaplaceMechanism}
\begin{align}
\mathcal{A}(\mathcal{D}) := \mathcal{R}(\mathcal{D}) + \tilde{\xi}, \label{A}
\end{align}
where $\tilde{\xi} \in \mathbb{R}^n$ whose elements are independent and identically distributed (IID) random variables drawn from a Laplace distribution with zero mean and a variance $2 (\bar{\Delta}_1 / \bar{\epsilon})^2$. 
Here, $\bar{\epsilon}$ is a privacy budget from \eqref{def_differential_privacy}, and $\bar{\Delta}_1$ is an $L_1$-sensitivity defined as
\begin{align}
& \bar{\Delta}_1 := \max_{\forall (\mathcal{D},\mathcal{D}')} \| \mathcal{R}(\mathcal{D}) - \mathcal{R}(\mathcal{D}') \|_1,  \label{L1_sensitivity}
\end{align}
\end{subequations}
where $(\mathcal{D},\mathcal{D}')$ are two neighboring datasets.
The Laplace mechanism $\mathcal{A}$ in \eqref{LaplaceMechanism} satisfies $(\bar{\epsilon},0)$-DP.
\end{lemma}

\begin{lemma}[Gaussian mechanism \cite{dwork2014algorithmic}]\label{lem:existing-gaussian}
The Gaussian mechanism takes the form of \eqref{A} where $\tilde{\xi} \in \mathbb{R}^n$ is a noise vector whose elements are IID random variables drawn from a normal distribution with zero mean and a variance $2 \ln(1.25 / \bar{\delta}) (\bar{\Delta}_2 / \bar{\epsilon})^2$, where
\begin{align}
\bar{\Delta}_2:= \max_{\forall (\mathcal{D},\mathcal{D}')} \| \mathcal{R}(\mathcal{D}) - \mathcal{R}(\mathcal{D}') \|_2  \label{L2_sensitivity}
\end{align}
represents $L_2$-sensitivity. The Gaussian mechanism $\mathcal{A}$ satisfies $(\bar{\epsilon},\bar{\delta})$-DP, where $\bar{\delta} > 0$.
\end{lemma}
Even though the Gaussian mechanism that ensures  $(\bar{\epsilon},\bar{\delta})$-DP with $\bar{\delta}>0$ provides a weaker privacy guarantee compared with the Laplace mechanism that ensures $(\bar{\epsilon},0)$-DP, the Gaussian mechanism could lead to less noise since its variance relies on $L_2$ sensitivity, which is smaller than $L_1$ sensitivity used in the Laplace mechanism.
For this reason, the Gaussian mechanism has been widely adopted in numerous learning algorithms since the amount of noise introduced affects the learning performance.
Moreover, the Gaussian mechanism is favorable in FL algorithms where DP is applied for every iteration of the algorithms. 
The reason is that a mechanism that permits $T$ adaptive interactions with the Gaussian mechanisms ensures $( \bar{\epsilon} \sqrt{2T\ln(1/\bar{\delta}')} + 2k\epsilon^2, T \bar{\delta} + \bar{\delta}' ) $-DP, which is stronger than a mechanism that permits $T$ adaptive interactions with the Laplace mechanisms and ensures $(T\bar{\epsilon}, T\bar{\delta})$-DP, as shown in \cite{dwork2010boosting}.
Moreover, Abadi et al.~\cite{abadi2016deep} proposed a stronger composition theorem for the Gaussian mechanisms.

\begin{lemma}[Stronger composition theorem \cite{abadi2016deep}] \label{lem:composition}
An algorithm $\mathcal{A}$ that consists of a sequence of $T$ adaptive Gaussian mechanisms, namely, $\mathcal{A}^1, \ldots, \mathcal{A}^T$, each of which is proven to be $(\bar{\epsilon},\bar{\delta})$-DP with $\bar{\delta}>0$, is $(\epsilon', \bar{\delta})$-DP, where 
\begin{subequations}
\begin{align}
\ln(\bar{\delta}) = \ &  \min_{\tau \in \mathbb{Z}_+} \sum_{t=1}^T \alpha^t (\tau) - \tau \epsilon',  \label{lem:composition-1} \\
\alpha^{t} (\tau) := & \max_{\forall (\mathcal{D}, \mathcal{D}')} \ln \Bigg( \mathbb{E}_{\psi} \Bigg[ \Big( \frac{\textbf{pdf} \big( \mathcal{A}^{t} (\mathcal{D}) = \psi \big)}{\textbf{pdf} \big(  \mathcal{A}^{t} (\mathcal{D}') = \psi \big)} \Big)^{\tau} \Bigg] \Bigg), \ \forall t \in [T]. \label{lem:composition-2}
\end{align}
\end{subequations} 
Here $\alpha^{t} (\tau)$ represents all possible log of the moment-generating function at $\tau \in \mathbb{Z}_+$.
\end{lemma} 

In this paper we propose Laplace and Gaussian mechanisms in the context of the linearized ADMM with \texttt{ObjPert}, which are different from the Laplace and Gaussian mechanisms in Lemmas \ref{lem:existing-laplace} and \ref{lem:existing-gaussian}, respectively, with respect to the sensitivity computations as in \eqref{L1_sensitivity} and \eqref{L2_sensitivity}, respectively.
We also show that the  proposed Gaussian mechanism allows the stronger composition theorem to be applied.


\section{Differentially private linearized ADMM with multiple local updates} \label{sec:dpadmm}
In this section we propose a differentially private linearized ADMM algorithm for solving the distributed convex optimization model \eqref{model:dist}.
The proposed algorithm is equipped with (i) \textit{multiple local updates} that can reduce costs of communication between agents and (ii) \textit{objective perturbation} (i.e., \texttt{ObjPert}), which ensures data privacy in a DP manner without losing feasibility.
The proposed algorithm is especially effective when $\mathcal{W}_p \subset \mathbb{R}^n$, mainly because each agent communicates randomized (to ensure data privacy) but feasible local solutions that could affect the overall convergence of the algorithm.
In the rest of this section we describe the proposed algorithm (Section \ref{sec:proposed_algo}) and present a detailed comparison with the existing DP algorithm based on \texttt{OutPert} (Section \ref{sec:baseline}).
The privacy and convergence analyses are presented in Sections \ref{sec:privacy} and \ref{sec:convergence}, respectively.

\subsection{The proposed algorithm} \label{sec:proposed_algo}
The proposed algorithm is built on the linearized ADMM, that is, $\{ \eqref{ADMM-1} \rightarrow \eqref{ADMM-2-Prox} \rightarrow \eqref{ADMM-3} \}_{t=1}^T$.
Note that each agent solves \eqref{ADMM-2-Prox} using its local data and machines, which can be considered as a single local update.
In our algorithm, each agent conducts multiple local updates.
Specifically, for every $e \in [E]$, where $E$ is the number of local updates, each agent solves the following problem:
\begin{align}
z_p^{t,e+1} \leftarrow  \arg \min_{z_p  \in \mathcal{W}_p } \ &  \langle f'_p(z^{t,e}_p), z_p \rangle + \frac{1}{2\eta^{t}}\|z_p - z^{t,e}_p\|^2 + \frac{\rho^t}{2} \|w^{t+1}-z_p + \frac{1}{\rho^{t}} \lambda^t_p  \|^2.  \label{ADMM-2-MLU} 
\end{align}
The multiple local updates technique has been utilized in the federated averaging algorithm \cite{mcmahan2017communication} to reduce the communication costs. 
Inspired by that work, we introduce this technique in the linearized ADMM.

We note that the objective function of \eqref{ADMM-2-MLU} is strongly convex and $\mathcal{W}_p$ is a closed convex set; thus it has a unique optimal solution.
Therefore, the mapping from $f'_p(z_p^{t,e})$ to $z_p^{t,e+1}$ is injective (i.e., every $z_p^{t,e+1}$ has at most one $f'_p(z_p^{t,e})$). Therefore, the release of $z_p^{t,e+1}$ can lead to the leakage of $f'_p(z_p^{t,e})$, which contains the data information. 
It has been shown that the data can be reconstructed from the gradients (e.g., \cite{zhu2019deep, zhao2020idlg, jeon2021gradient}).

To protect data, one can simply add some calibrated noise $\tilde{\xi}^{t,e}_p$ to the output $z^{t,e+1}_p$, which is known as the output perturbation (i.e., \texttt{OutPert}) method. Note that the Laplace and Gaussian mechanisms described in Section \ref{sec:prelim_dp} belong to \texttt{OutPert} since they add Laplacian or Gaussian noise to the output directly to obtain randomized outputs to release. 
By noting that the resulting randomized outputs (that ensure data privacy) may be infeasible to $\mathcal{W}_p$, which may affect overall convergence of the algorithm, we alternatively add noise to the objective function of the constrained subproblem \eqref{ADMM-2-MLU} so that the resulting randomized outputs to release are always feasible to $\mathcal{W}_p$. 
Specifically, we add an affine function $\frac{1}{2\rho^t} \| \tilde{\xi}^{t,e}_p \|^2 - \langle w^{t+1} - z_p + \frac{1}{\rho^t} \lambda^t_p, \tilde{\xi}^{t,e}_p \rangle$ to \eqref{ADMM-2-MLU}, resulting in
\begin{align}
  z^{t,e+1}_p \leftarrow & \arg \min_{z_p \in \mathcal{W}_p} \ \langle f'_p(z^{t,e}_p), z_p \rangle + \frac{1}{2\eta^{t}} \| z_p - z_p^{t,e}\|^2 + \frac{\rho^t}{2} \|w^{t+1}-z_p + \frac{1}{\rho^t}(\lambda^t_p - \tilde{\xi}^{t,e}_p) \|^2,  \label{DPADMM-2-Prox}
\end{align}  
where $\tilde{\xi}^{t,e}_p \in \mathbb{R}^{n}$ is a noise vector that is guaranteed to preserve data privacy in a DP manner. In Section \ref{sec:privacy} we discuss how to generate such noise.
 

The proposed algorithm is presented in Algorithm \ref{algo:DP-IADMM-Prox}.
The computation at the central server is described in lines 1--9, while the local computation for each agent $p$ is described in lines 11--23.
In lines 2--3, the initial points are sent from the server to all agents.
In line 5, the central server computes a global solution $w^{t+1}$ via the equation derived from the optimality condition of \eqref{ADMM-1}.
In line 6, the server broadcasts $w^{t+1}$ to all local agents.
In lines 15--22, each agent receives $w^{t+1}$ from the server, conducts local updates for $E$ times, and sends the resulting local solution $z_p^{t+1}$ to the server.
Note that $z_p^{t+1}$ is an average of $E$ numbers of randomized outputs, each of which ensures data privacy.
In line 8 and line 23, the dual updates \eqref{ADMM-3} are performed at the server and the local agents individually. 
Note that these dual updates are identical since the initial points at the server and the local agents are the same.  

\begin{algorithm}[!ht]
\caption{Proposed DP linearized ADMM with multiple local updates}
\label{algo:DP-IADMM-Prox}
\begin{algorithmic}[1]
\STATE \textbf{(Server):}
\STATE Initialize $\lambda^1_1, \ldots, \lambda^1_P, z^1_1, \ldots, z^1_P$.
\STATE \texttt{Send} $\lambda^1_p, z^1_p$ to all agent $p \in [P]$ (\textbf{to line 12}).
\FOR{$t \in [T]$}
\STATE $w^{t+1} \leftarrow \frac{1}{P} \sum_{p=1}^P ( z_p^t - \frac{1}{\rho^t} \lambda^t_p) $.
\STATE \texttt{Send} $w^{t+1}$ to all agents (\textbf{to line 15}).
\STATE \texttt{Receive} $z^{t+1}_p$ from all agents (\textbf{from line 22}).
\STATE $\lambda^{t+1}_p \leftarrow \lambda^t_p + \rho^t(w^{t+1}-z^{t+1}_p)$ for all $p \in [P]$.
\ENDFOR
\STATE  
\STATE \textbf{(Agent $p \in [P]$):}  
\STATE \texttt{Receive} $\lambda^1_p, z^1_p$ from the server (\textbf{from line 3}).
\STATE Initialize $z^{0,E+1}_p = z^1_p$.
\FOR{$t \in [T]$}
\STATE \texttt{Receive} $w^{t+1}$ from the server (\textbf{from line 6}).    
\STATE Set $z_{p}^{t,1} \leftarrow z_p^{t-1,E+1}$
\FOR{$e \in [E]$}
\STATE Compute $z_p^{t,e+1}$ via \eqref{DPADMM-2-Prox}.
\ENDFOR        
\STATE $z^{t+1}_p \leftarrow \frac{1}{E} \sum_{e=1}^E z_p^{t,e+1}$.    
\STATE \texttt{Send} $z^{t+1}_p$ to the server (\textbf{to line 7}).
\STATE $\lambda^{t+1}_p \leftarrow \lambda^t_p + \rho^t(w^{t+1}-z^{t+1}_p)$.
\ENDFOR
\end{algorithmic}
\end{algorithm}
    
The benefits of Algorithm \ref{algo:DP-IADMM-Prox} are as follows.
First, 
solving \eqref{DPADMM-2-Prox} can be computationally cheaper than its exact version. Furthermore, for some special cases when $\mathcal{W}_p$ is a box constraint, \eqref{DPADMM-2-Prox} admits a closed-form solution expression.
Second, the quality of the solution can be improved via the multiple local updates that could result in saving communication rounds.
Third, the amount of communication is reduced by excluding the communication of the dual information.
This is different from the existing work \cite{zhou2021communication} that considers both multiple local primal and dual updates, namely, solving \eqref{ADMM-2-Prox} and \eqref{ADMM-3} multiple times, resulting in communicating not only local solution $z^{t+1}_p$ but also dual information $\lambda^{t+1}_p$.
Fourth and most important, data privacy is guaranteed in a DP manner for any communication rounds.

\subsection{Comparison with the existing DP algorithm} \label{sec:baseline}
We note that the authors in \cite{huang2019dp, huang2020differentially} proposed a differentially private linearized ADMM with multiple local updates when $\mathcal{W}_p=\mathbb{R}^n$.
Specifically, they employed a Gaussian mechanism (described in Section \ref{sec:prelim_dp}) that adds Gaussian noise $\tilde{\xi}^{t,e}_p$ to the output $z_p^{t,e+1}$ in \eqref{ADMM-2-MLU}.
We note that when $\mathcal{W}_p=\mathbb{R}^n$, adding Gaussian noise to the output (i.e., \texttt{OutPert}) is exactly the same as adding $- ( 1/\eta^t + \rho^t) \langle \tilde{\xi}^{t,e}_p, z_p \rangle$ to the objective function of \eqref{ADMM-2-MLU} (i.e., \texttt{ObjPert}).
We formalize the argument in the following remark.
\begin{remark}
Suppose that $\mathcal{W}_p=\mathbb{R}^n$ in \eqref{ADMM-2-MLU}.
The output perturbation is given by 
\begin{align}
z_p^{t,e+1} = \Big(\frac{1}{\eta^t} +\rho^t \Big)^{-1} \Big( - f'_p(z_p^{t,e}) + \frac{1}{\eta^t} z_p^{t,e} +\rho^t w^{t+1} + \lambda^t_p \Big) + \tilde{\xi}^{t,e}_p, \label{remark_opt_cond}
\end{align}
where $\tilde{\xi}^{t,e}_p$ is the Gaussian noise and the first term is derived by the optimality condition of \eqref{ADMM-2-MLU} with $\mathcal{W}_p=\mathbb{R}^n$.
Now consider the following problem constructed by the objective perturbation method:
\begin{align}
z_p^{t,e+1} \leftarrow  \arg \min_{z_p  \in \mathbb{R}^n } \  & \langle f'_p(z^{t,e}_p), z_p \rangle + \frac{1}{2\eta^{t}}\|z_p - z^{t,e}_p\|^2 + \frac{\rho^t}{2} \|w^{t+1}-z_p + \frac{1}{\rho^{t}} \lambda^t_p  \|^2 \nonumber \\
& - ( 1/\eta^t + \rho^t) \langle \tilde{\xi}^{t,e}_p, z_p \rangle. \nonumber
\end{align}
From its optimality condition, one can derive \eqref{remark_opt_cond}.
\end{remark}

We also note that when $\mathcal{W}_p = \mathbb{R}^n$, the multiple local update used in \cite{hager2020convergence} is exactly the same as solving \eqref{ADMM-2} by a subgradient descent algorithm up to $E$ iterations given by  
\begin{align}
z^{t, e+1}_p = z^{t,e}_p - s \Big( f_p'(z_p^{t,e})  - \rho^t (w^{t+1}-z_p^{t,e} + \frac{1}{\rho^t} \lambda^t )  \Big), \ \forall e \in [E], \nonumber
\end{align}
where $s = \eta^t/(\eta^t \rho^t +1)$ is a constant step size. 
To see this, one can derive the optimality condition of \eqref{ADMM-2-MLU} when $\mathcal{W}_p = \mathbb{R}^n$:
\begin{align*}
0 & = f'_p(z_p^{t,e}) + \frac{1}{\eta^t} (z_p^{t,e+1} - z_p^{t,e}) - \rho^t(w^{t+1} - z_p^{t,e+1}  + \frac{1}{\rho^t} \lambda^t_p ) \\  
& = f'_p(z_p^{t,e}) + \frac{1}{\eta^t} (z_p^{t,e+1} - z_p^{t,e}) - \rho^t(w^{t+1} - z_p^{t,e}  + \frac{1}{\rho^t} \lambda^t_p ) - \rho^t (z_p^{t,e} - z_p^{t,e+1}) \\
& = \eta^t \Big( f'_p(z_p^{t,e}) - \rho^t(w^{t+1} - z_p^{t,e}  + \frac{1}{\rho^t} \lambda^t_p )  \Big) - (\eta^t \rho^t+1) (z_p^{t,e} - z_p^{t,e+1}) \\
& = -\frac{\eta^t}{\eta^t \rho^t+1} \Big( f'_p(z_p^{t,e})  + \rho^t(w^{t+1} - z_p^{t,e}  + \frac{1}{\rho^t} \lambda^t_p )  \Big)  +  (z_p^{t,e} - z_p^{t,e+1}).
\end{align*}
To summarize, when $\mathcal{W}_p = \mathbb{R}^n$, there is no significant difference between \texttt{OutPert} and \texttt{ObjPert}.
Also, the multiple local updates simply involve solving the exact problem \eqref{ADMM-2} by the subgradient descent algorithm for a finite number of iterations.

\begin{figure}[!ht]  
\centering
\includegraphics[width=0.45\textwidth]{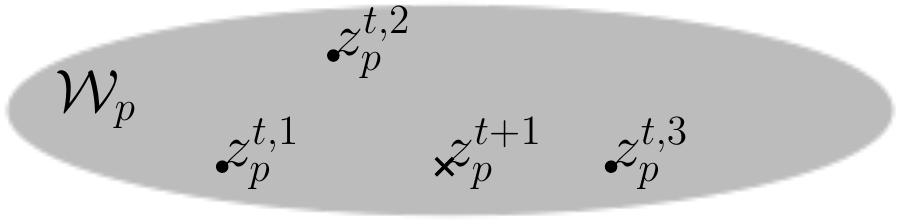}  
\includegraphics[width=0.45\textwidth]{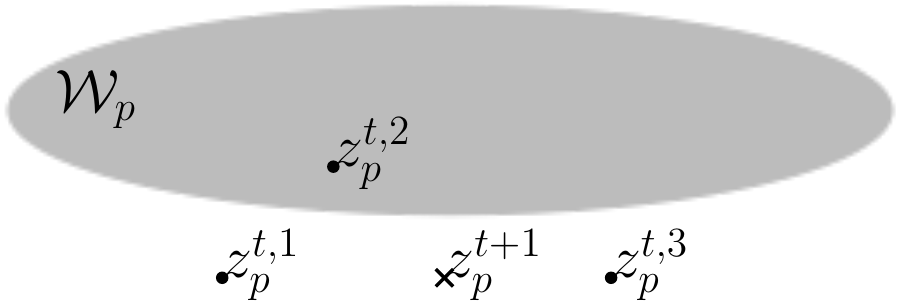}  
\caption{Comparison of \eqref{DPADMM-2-Prox} and \eqref{baseline} for fixed $t$, $p$, and $E=3$. 
Every randomized outputs produced by \eqref{DPADMM-2-Prox} is feasible, whereas those produced by \eqref{baseline} may not be.}
\label{fig:feasible}
\end{figure}

Different from the literature,  this paper focuses on a situation when $\mathcal{W}_p \subset \mathbb{R}^n$, which appears in many application areas. 
In this case, we aim to show that \texttt{ObjPert} outperforms \texttt{OutPert}.
To do so, we construct a baseline algorithm based on \texttt{OutPert}, which is a form of Algorithm \ref{algo:DP-IADMM-Prox} where \eqref{DPADMM-2-Prox} in line 18 is replaced with
\begin{align}
z_p^{t,e+1} \gets z_p^{t,e+1} \text{ from \eqref{ADMM-2-MLU} } + \tilde{\xi}^{t,e}_p, \label{baseline}
\end{align}
where $\tilde{\xi}^{t,e}_p$ should be sampled to ensure DP.
However, the randomized output $z_p^{t+1} = \frac{1}{E} \sum_{e=1}^E z_p^{t,e+1}$, where $z_p^{t,e+1}$ is computed by \eqref{baseline}, may not be feasible to $\mathcal{W}_p$, whereas the randomized output $z_p^{t+1}$ in line 20 of Algorithm \ref{algo:DP-IADMM-Prox} is always feasible (see Figure \ref{fig:feasible}).

While the privacy and convergence analyses of the baseline algorithm can be immediately borrowed from \cite{huang2020differentially}, it is not the case for the proposed algorithm (i.e., Algorithm \ref{algo:DP-IADMM-Prox}) because \texttt{ObjPert} is different from \texttt{OutPert} when $\mathcal{W}_p \subset \mathbb{R}^n$. 
Notably, noise generated for ensuring DP in Algorithm \ref{algo:DP-IADMM-Prox} is based on the sensitivity of the subgradients, which is different from the sensitivity of the true outputs (i.e., \eqref{L1_sensitivity} or \eqref{L2_sensitivity}) in the literature. 
Also, \texttt{ObjPert} under the existence of $\mathcal{W}_p \subset \mathbb{R}^n$ requires  construction of a set of mechanisms, each of which is $(\bar{\epsilon}, \bar{\delta})$-DP, that provides a sequence of the randomized outputs that converges to a feasible randomized output, which is not needed in the privacy analysis for the case when $\mathcal{W}_p = \mathbb{R}^n$ as in \cite{huang2020differentially}. 
More details on our privacy analysis will be discussed in Section \ref{sec:privacy}.
Furthermore, we show that the rate of convergence of Algorithm \ref{algo:DP-IADMM-Prox} is sublinear without any error bound, which is different from the convergence analysis in \cite{huang2020differentially}, where the rate of convergence is $\mathcal{O}(1/\sqrt{T}) + \mathcal{O}(1/\bar{\epsilon}^2)$, where the error bound $\mathcal{O}(1/\bar{\epsilon}^2)$ increases as $\bar{\epsilon}$ decreases for stronger data privacy. 
More details on our convergence analysis  will be discussed in Section \ref{sec:convergence}.
\section{Privacy analysis} \label{sec:privacy}
In this section we provide a privacy analysis of the proposed algorithm.
Specifically, we will use the following lemma to show that the constrained subproblem \eqref{DPADMM-2-Prox} for fixed $t$, $e$, and $p$ provides a randomized output that guarantees data privacy in a DP manner. 
 
\begin{lemma}[Theorem 1 of \cite{kifer2012private}] \label{lemma:kifer}
Let $\mathcal{A}$ be a randomized algorithm induced by the random noise $\tilde{\xi}$ that provides output $\phi(\mathcal{D},\tilde{\xi})$.
Let $\{\mathcal{A}_{\ell}\}_{\ell=1}^{\infty}$ be a sequence of randomized algorithms, each of which is induced by $\tilde{\xi}$ and provides output $\phi^{\ell}(\mathcal{D},\tilde{\xi})$.
If $\mathcal{A}_{\ell}$ is $(\bar{\epsilon},\bar{\delta})$-DP for all $\ell$ and satisfies a pointwise convergence condition, namely, $\lim_{{\ell} \rightarrow \infty} \phi^{\ell}(\mathcal{D},\tilde{\xi}) = \phi(\mathcal{D},\tilde{\xi})$, then $\mathcal{A}$ is also $(\bar{\epsilon},\bar{\delta})$-DP, where $\bar{\delta} \geq 0$.
\end{lemma} 

In what follows, we construct a sequence $\{\mathcal{A}_{\ell}\}_{\ell=1}^{\infty}$ of randomized algorithms that satisfy the pointwise convergence condition in Section \ref{sec:randomized_algo}, and we propose two different mechanisms $\mathcal{A}_{\ell}$ for any $\ell$ as in Section \ref{sec:mechanisms}.

\subsection{A sequence of randomized algorithms} \label{sec:randomized_algo}
For the rest of this section we fix $t \in [T]$, $e \in [E]$, and $p \in [P]$. 
Also, we denote the objective function of \eqref{DPADMM-2-Prox}, which is strongly convex, by  
\begin{align}
  & G^{t,e}_p(z_p; \tilde{\xi}^{t,e}_p) :=  \langle f'_p(z^{t,e}_p), z_p \rangle   +   \frac{1}{2\eta^{t}} \| z_p - z_p^{t,e}\|^2 + \frac{\rho^t}{2} \|w^{t+1}-z_p + \frac{1}{\rho^t}(\lambda^t_p - \tilde{\xi}^{t,e}_p) \|^2 \label{fn_G} 
\end{align}
and the feasible region of \eqref{DPADMM-2-Prox} by
\begin{align*}
     \mathcal{W}_p := \{ z_p \in \mathbb{R}^{n} : h_m (z_p) \leq 0, \ \forall m \in [M_p] \}, 
\end{align*}
where $h_m$ is convex and twice continuously differentiable and $M_p$ is the total number of inequalities.

\begin{figure}[!ht]  
\centering
\includegraphics[width=0.45\textwidth]{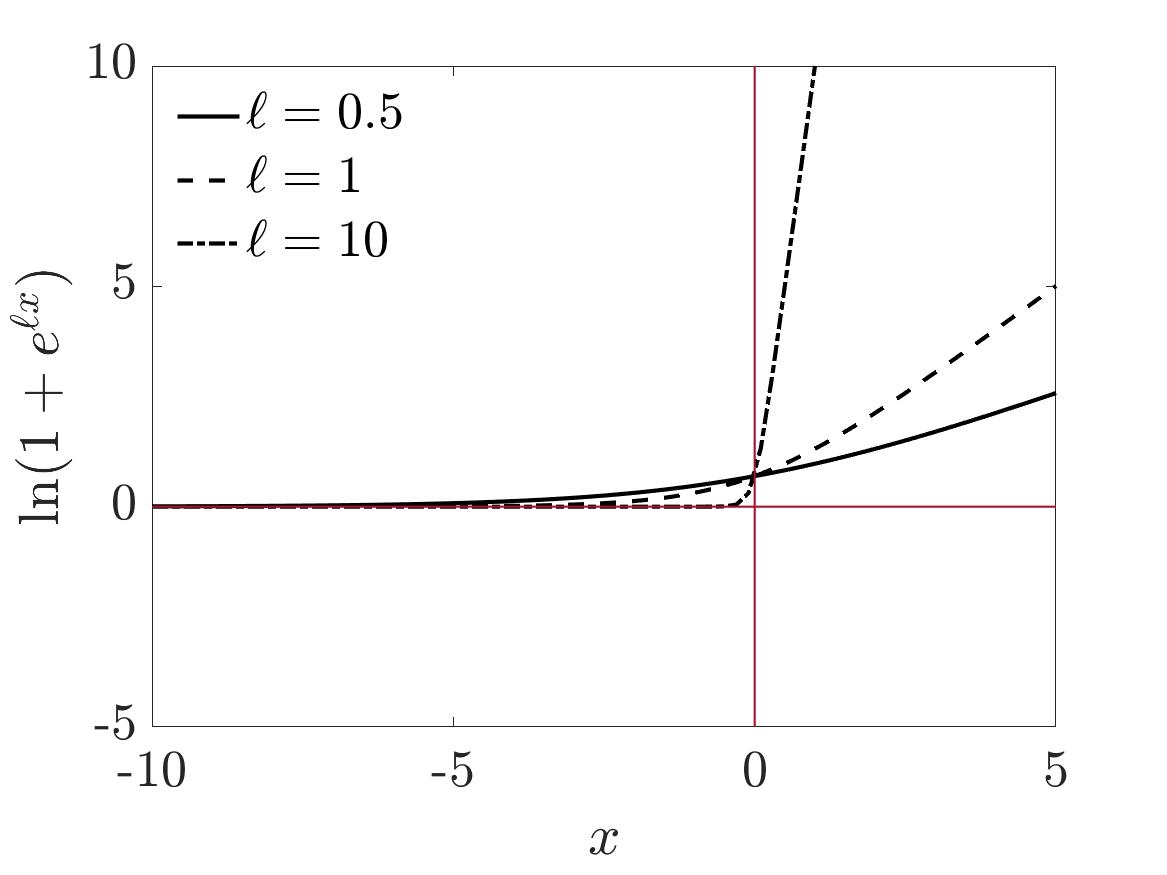}  
\includegraphics[width=0.45\textwidth]{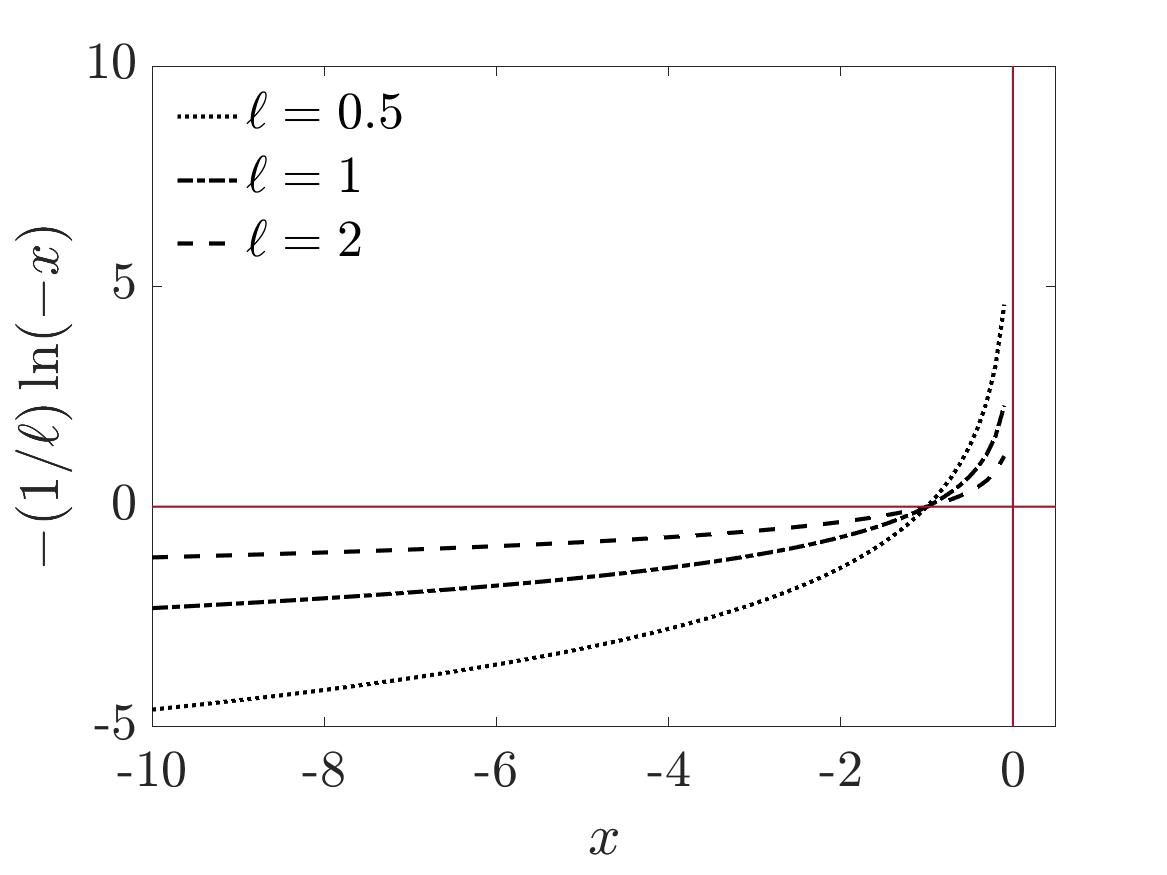}  
\caption{Example of the penalty function in \eqref{logfn-prox} with $M_p=1$ and $h_1(z_p) = z_p$ compared with the logarithmic barrier function.}
\label{fig:penalty}
\end{figure}

By utilizing an indicator function $\mathcal{I}_{\mathcal{W}_p}(z_p)$ that outputs zero if $z_p \in \mathcal{W}_p$ and $\infty$ otherwise, \eqref{DPADMM-2-Prox} can be expressed by the following problem:
\begin{align*}
\min_{z_p \in \mathbb{R}^{n}} \ G^{t,e}_p (z_p; \tilde{\xi}^{t,e}_p) + \mathcal{I}_{\mathcal{W}_p}(z_p).
\end{align*}
We observe that the indicator function can be approximated by the following penalty function:
\begin{align}
& g_{\ell} (z_p) :=  \sum_{m=1}^{M_p} \ln ( 1+ e^{\ell h_m(z_p)}), \label{logfn-prox}  
\end{align}
where $\ell > 0$. 
Increasing $\ell$ enforces the feasibility, namely, $h_m(z_p) \leq 0$, resulting in $g_{\ell}(z_p) \rightarrow 0$.
It is similar to the logarithmic barrier function (LBF) $-(1/\ell) \sum_{m=1}^M \ln (-h_m(z_p))$, in that the approximation becomes closer to the indicator function as $\ell \rightarrow \infty$.
However, the penalty function $g_{\ell}$ is different from LBF in that $z_p$ is not restricted.
See Figure \ref{fig:penalty} for an example of the penalty function compared with LBF.
By replacing the indicator function with the penalty function in \eqref{logfn-prox}, we construct the following \textit{unconstrained} problem (which is $\mathcal{A}_{\ell}$ in the context of Lemma \ref{lemma:kifer}):
\begin{align}
z^{t,e+1}_{p\ell} \leftarrow  \arg \min_{z_p \in \mathbb{R}^{n}} \ & G^{t,e}_p (z_p; \tilde{\xi}^{t,e}_p)  + g_{\ell}(z_p), \label{ADMM-2-Prox-log}
\end{align}
where the objective function is strongly convex because $g_{\ell}$ is convex over all domains and $G^{t,e}_p$ is strongly convex. Therefore, $z^{t,e+1}_{p\ell}$ is the unique optimal solution.

\begin{theorem}\label{thm:pointwise_convergence}
  For fixed $t$, $e$, and $p$,  
  \begin{align}
  \lim_{\ell \rightarrow \infty} z^{t,e+1}_{p \ell}  = z_p^{t,e+1},  
  \end{align}
  where $z_p^{t,e+1}$ and $z^{t,e+1}_{p \ell}$ are from \eqref{DPADMM-2-Prox} and \eqref{ADMM-2-Prox-log}, respectively.
\end{theorem}
\begin{proof}
  We develop the proof by contradiction.
  To this end, we suppose that $z^{t,e+1}_{p\ell}$ converges to $\widehat{z} \neq z^{t,e+1}_p$ as $\ell$ increases.
  Consider $\zeta := \| \widehat{z} - z^{t,e+1}_p \| / 2$.
  Since $z^{t,e+1}_{p\ell}$  converges to $\widehat{z}$, there exists $\ell' > 0$ such that $\| \widehat{z} - z^{t,e+1}_{p\ell} \| < \zeta$ for all $\ell \geq \ell'$.
  By the triangle inequality, we have
  \begin{subequations}
  \begin{align}
  \| z^{t,e+1}_{p\ell} - z^{t,e+1}_p \| \geq \| \widehat{z} - z^{t,e+1}_p \| - \| \widehat{z} - z^{t,e+1}_{p\ell} \| > 2\zeta - \zeta = \zeta, \ \forall \ell \geq \ell'. \label{triangle}
  \end{align}
  Since $G^{t,e}_p$ is strongly convex with a constant $\mu > 0$, we have
  \begin{align}
  G^{t,e}_p(z^{t,e+1}_{p\ell}) - G^{t,e}_p(z^{t,e+1}_p) \geq \textstyle \frac{\mu}{2} \| z^{t,e+1}_{p\ell} - z^{t,e+1}_p \|^2 >  \frac{\mu \zeta^2}{2}, \ \forall \ell \geq \ell', \label{strong_triangle}
  \end{align}
  where the last inequality holds by \eqref{triangle}.
  By adding $g_{\ell} (z^{t,e+1}_{p\ell}) \geq 0$ to the left-hand side of \eqref{strong_triangle}, we derive the following inequality:
  \begin{align}
  \big\{ G^{t,e}_p(z^{t,e+1}_{p\ell}) + g_{\ell} (z^{t,e+1}_{p\ell})  \big\} - G^{t,e}_p(z^{t,e+1}_p)  > \textstyle \frac{\mu \zeta^2}{2}, \ \forall \ell \geq \ell'. \label{diff_0}
  \end{align}  
  The continuity of $G^{t,e}_p: \mathcal{W}_p \mapsto \mathbb{R}$ at $z^{t,e+1}_p$ implies that for every $\epsilon > 0$, there exists a $\delta > 0$ such that for all $z \in \mathcal{W}_p$:
  \begin{align}
  z \in \mathcal{B}_{\delta}(z^{t,e+1}_p) := \{ z \in \mathbb{R}^{n} : \| z - z^{t,e+1}_p \| < \delta \} \ \Rightarrow \  G^{t,e}_p(z) - G^{t,e}_p(z^{t,e+1}_p)  < \epsilon. \label{continuity}
  \end{align}
  Consider $\tilde{z} \in \mathcal{B}_{\delta}(z^{t,e+1}_p) \cap \textbf{relint}(\mathcal{W}_p)$, where \textbf{relint} indicates the relative interior.
  Since $h_m(\tilde{z}) < 0$ for all $m \in [M_p]$, $g_{\ell}(\tilde{z})$ goes to zero as $\ell$ increases.
  Hence, there exists $\ell'' > 0$ such that
  \begin{align}
    g_{\ell}(\tilde{z}) =  \sum_{m=1}^{M_p} \ln (1 + e^{\ell h_m(\tilde{z})}) < \epsilon, \ \forall \ell \geq \ell''. \label{ineq_z_tilde}
  \end{align}
  For all $\ell \geq \ell''$, we derive the following inequalities:
  \begin{align}
    G^{t,e}_p(z^{t,e+1}_{p\ell}) + g_{\ell}(z^{t,e+1}_{p\ell}) \leq G^{t,e}_p(\tilde{z}) + g_{\ell}(\tilde{z}) < G^{t,e}_p(\tilde{z}) + \epsilon < G^{t,e}_p(z^{t,e+1}_p) + 2\epsilon, \label{sandwich_ineq}
  \end{align}
  where
  the first inequality holds because $z^{t,e+1}_{p\ell}$ is the optimal solution of \eqref{ADMM-2-Prox-log},
  the second inequality holds by \eqref{ineq_z_tilde}, and
  the last inequality holds by \eqref{continuity}.
  To see the contradiction, consider $\epsilon \in (0, \frac{\mu \zeta^2}{4})$.
  Then we have
  \begin{align}
  \big\{ G^{t,e}_p(z^{t,e+1}_{p\ell}) + g_{\ell} (z^{t,e+1}_{p\ell})  \big\} - G^{t,e}_p(z^{t,e+1}_p) \underbrace{<}_{\text{from } \eqref{sandwich_ineq}} 2 \epsilon < \frac{\mu \zeta^2}{2}, \ \forall \ell \geq \ell''. \label{diff_1}
  \end{align}  
  Therefore, for all $\ell \geq \max \{ \ell', \ell'' \}$, \eqref{diff_0} and \eqref{diff_1} contradict.
\end{subequations}
\end{proof}

\subsection{Two mechanisms} \label{sec:mechanisms}
In this section we propose two mechanisms, namely, variants of the existing Laplace and Gaussian mechanisms in Lemmas \ref{lem:existing-laplace} and \ref{lem:existing-gaussian}, respectively, based on the proposed \eqref{ADMM-2-Prox-log} for fixed $t$, $e$, $p$, and $\ell$.

First, we propose to sample $\tilde{\xi}^{t,e}_p$ in \eqref{ADMM-2-Prox-log} from a Laplace distribution with zero mean and a variance $2 (\bar{\Delta}^{t,e}_{p,1} / \bar{\epsilon})^2$, where $\bar{\Delta}^{t,e}_{p,1}$ is defined as follows:
\begin{align}
  \bar{\Delta}^{t,e}_{p,1} := \max_{ \forall (\mathcal{D}_p, \mathcal{D}_p')} \| f'_p(z^{t,e}_p;\mathcal{D}_p) - f'_p(z^{t,e}_p;\mathcal{D}'_p)\|_1, \label{redefined_L1sensitivity}
\end{align}
which is different from \eqref{L1_sensitivity} in that the proposed $L_1$ sensitivity \eqref{redefined_L1sensitivity} is based on the difference of subgradients (i.e., inputs of \eqref{ADMM-2-Prox-log}) while the existing $L_1$ sensitivity \eqref{L1_sensitivity} is based on the difference of true outputs.
In the following theorem, we show that the proposed Laplace mechanism, namely, \eqref{ADMM-2-Prox-log} with the Laplacian noise, can achieve DP.
 
\begin{theorem} \label{thm:laplace}  
  For fixed $t$, $e$, $p$, and $\ell$, the proposed \eqref{ADMM-2-Prox-log} with the aforementioned Laplacian noise is $(\bar{\epsilon},0)$-DP: 
  \begin{align}
    \label{DP_1}
    \mathbb{P} ( z^{t,e+1}_{p \ell} (\mathcal{D}_p )\in \mathcal{S} ) \leq  e^{\bar{\epsilon}} \mathbb{P} ( z^{t,e+1}_{p \ell} (\mathcal{D}'_p )\in \mathcal{S} ), 
  \end{align}  
  for all possible events $\mathcal{S} \subset \mathbb{R}^{n}$ and neighboring datasets $(\mathcal{D}_p, \mathcal{D}'_p)$.
\end{theorem}
\begin{proof}
  It suffices to show that the following is true:
  \begin{subequations}
  \begin{align}
    \textbf{pdf} \big( z^{t,e+1}_{p\ell} (\mathcal{D}_p) = \psi \big)
    \leq e^{\bar{\epsilon}} \ \textbf{pdf} \big( z^{t,e+1}_{p\ell} (\mathcal{D}'_p) = \psi \big), \ \forall \psi \in \mathbb{R}^{n}, \label{DP_PDF}    
  \end{align} 
  where \textbf{pdf} stands for probability density function.
  This is because integrating $\psi$ in \eqref{DP_PDF} over $\mathcal{S}$ yields \eqref{DP_1}.

  Consider $\psi \in \mathbb{R}^{n}$.
  If we have $z^{t,e+1}_{p \ell} (\mathcal{D}_p)=\psi$, then $\psi$ is the unique minimizer of \eqref{ADMM-2-Prox-log} because the objective function in \eqref{ADMM-2-Prox-log} is strongly convex.
  From the optimality condition of \eqref{ADMM-2-Prox-log}, we derive
  \begin{align}
    \tilde{\xi}^{t,e}_p (\psi; \mathcal{D}_p) = & - f'_p(z^{t,e}_p;\mathcal{D}_p) + \rho^t(w^{t+1}-\psi) + \lambda^t_p - \nabla g_{\ell}(\psi)  -  \frac{1}{\eta^{t}} \big( \psi - z_p^{t,e} \big), \label{correspondence}
  \end{align}
  where $\nabla g_{\ell}(\psi) = \sum_{m=1}^{M_p} \frac{\ell e ^{\ell h_m(\psi)}}{1+e^{\ell h_m(\psi)}} \nabla h_m(\psi)$, and we use $\tilde{\xi}^{t,e}_p (\psi; \mathcal{D}_p)$ and $\tilde{\xi}^{t,e}_p$ interchangeably.
  Note that the mapping from $\psi$ to $\tilde{\xi}^{t,e}_p$ via \eqref{correspondence} is injective. 
  Also the mapping is surjective because for all $\tilde{\xi}^{t,e}_p$, there exists $\psi$ (i.e., the unique minimizer of \eqref{ADMM-2-Prox-log}) such that \eqref{correspondence} holds.
  Therefore, the relation between $\psi$ and $\tilde{\xi}^{t,e}_p$ is bijective, which allows us to utilize the inverse function theorem \cite[Theorem 17.2]{billingsley}, namely,
    \begin{align}
      \textbf{pdf} \big( z^{t,e+1}_{p\ell} (\mathcal{D}_p) = \psi \big) \cdot \big|\textbf{det}[\nabla \tilde{\xi}^{t,e}_p (\psi;\mathcal{D}_p) ] \big| = \text{Lap} \big( \tilde{\xi}^{t,e}_p (\psi;\mathcal{D}_p); 0, \bar{\Delta}_{p,1}^{t,e} /\bar{\epsilon} \big), \label{inverse}
    \end{align}
    where 
    $\text{Lap}(\cdot; 0,\bar{\Delta}^{t,e}_{p,1} / \bar{\epsilon} ) $ represents a probability density function of Laplace distribution with zero mean and a variance $2(\bar{\Delta}^{t,e}_{p,1} / \bar{\epsilon})^2$, 
    \textbf{det} represents a determinant of a matrix, and
    $\nabla \tilde{\xi}^{t,e}_p (\psi;\mathcal{D}_p)$ represents a Jacobian matrix of the mapping from $\psi$ to $\tilde{\xi}^{t,e}_p$ in \eqref{correspondence}, namely,
    \begin{align}
    \nabla \tilde{\xi}^{t,e}_p (\psi;\mathcal{D}_p) =  (-\rho^t - 1/\eta^{t}) \mathbb{I}_{n} - \nabla \Big( \sum_{m=1}^{M_p} \frac{\ell e ^{\ell h_m(\psi)}}{1+e^{\ell h_m(\psi)}} \nabla h_m(\psi) \Big) , \label{Jacobian}
    \end{align}
    where $\mathbb{I}_{n}$ is an identity matrix of $n \times n$ dimensions.
    Since the Jacobian matrix is not affected by the dataset, we have
    \begin{align}
    \nabla \tilde{\xi}^{t,e}_p (\psi;\mathcal{D}_p) = \nabla \tilde{\xi}^{t,e}_p (\psi;\mathcal{D}'_p).  \label{jacobian}
    \end{align}
    Based on \eqref{inverse} and \eqref{jacobian}, we derive the following inequalities:
    \begin{align}
      & \ln \Bigg( \frac{\textbf{pdf} \big( z^{t,e+1}_{p \ell} (\mathcal{D}_p) = \psi \big)}{ \textbf{pdf}  \big( z^{t,e+1}_{p \ell} (\mathcal{D}'_p ) = \psi \big)}  \Bigg) = \ln \Bigg( \frac{\text{Lap} \big( \tilde{\xi}^{t,e}_p (\psi; \mathcal{D}_p); 0, \bar{\Delta}_{p,1}^{t,e}/\bar{\epsilon} \big) }{\text{Lap} \big( \tilde{\xi}^{t,e}_p (\psi; \mathcal{D}'_p); 0, \bar{\Delta}_{p,1}^{t,e} /\bar{\epsilon} \big) } \Bigg) \label{laplacian-1} \\
      = \ &    (\bar{\epsilon}/\bar{\Delta}_{p,1}^{t,e})(\| \tilde{\xi}^{t,e}_p (\psi;\mathcal{D}'_p) \|_1 - \| \tilde{\xi}^{t,e}_p (\psi;\mathcal{D}_p) \|_1)  \nonumber \\
      \leq \ &   (\bar{\epsilon}/\bar{\Delta}_{p,1}^{t,e})(\| \tilde{\xi}^{t,e}_p (\psi;\mathcal{D}'_p) - \tilde{\xi}^{t,e}_p (\psi;\mathcal{D}_p) \|_1)  \label{triangle_inequality} \\
      = \ &   (\bar{\epsilon}/\bar{\Delta}_{p,1}^{t,e})(\| f'_p(z^{t,e}_p;\mathcal{D}_p) - f'_p(z^{t,e}_p;\mathcal{D}'_p)\|_1) \leq \bar{\epsilon} , \label{4.8b}  
    \end{align}    
    \end{subequations}
    where
    the equality in \eqref{laplacian-1} holds because of \eqref{inverse} and \eqref{jacobian},
    the inequality \eqref{triangle_inequality} holds because of the triangle inequality,
    the equality in \eqref{4.8b} holds because of \eqref{correspondence}, and
    the inequality in \eqref{4.8b} holds by the definition of $\bar{\Delta}_{p,1}^{t,e}$ in \eqref{redefined_L1sensitivity}.    
\end{proof}

Next, we propose to sample $\tilde{\xi}^{t,e}_p$ in \eqref{ADMM-2-Prox-log} from a Gaussian distribution with zero mean and a variance $(\sigma^{t,e}_p)^2 := 2 \ln(1.25/\bar{\delta}) (\bar{\Delta}^{t,e}_{p,2} / \bar{\epsilon})^2$, where $\bar{\Delta}^{t,e}_{p,2}$ is defined as follows:
\begin{align}
  \bar{\Delta}^{t,e}_{p,2} := \max_{ \forall (\mathcal{D}_p, \mathcal{D}_p')} \| f'_p(z^{t,e}_p;\mathcal{D}_p) - f'_p(z^{t,e}_p;\mathcal{D}'_p)\|_2, \label{redefined_L2sensitivity}
\end{align}
which is different from \eqref{L2_sensitivity} in that the proposed $L_2$ sensitivity \eqref{redefined_L2sensitivity} is based on the difference of subgradients wheres the existing $L_2$ sensitivity is based on the difference of true outputs.
In the following theorem, we show that the proposed Gaussian mechanism, namely, \eqref{ADMM-2-Prox-log} with the Gaussian noise, can achieve DP.

\begin{theorem} \label{thm:gaussian}  
  For fixed $t$, $e$, $p$, and $\ell$, the proposed \eqref{ADMM-2-Prox-log} with the aforementioned Gaussian noise is $(\bar{\epsilon},\bar{\delta})$-DP: 
  \begin{align}
    \label{DP_2}
    \mathbb{P} ( z^{t,e+1}_{p \ell} (\mathcal{D}_p )\in \mathcal{S} ) \leq  e^{\bar{\epsilon}} \mathbb{P} ( z^{t,e+1}_{p \ell} (\mathcal{D}'_p )\in \mathcal{S} ) + \bar{\delta}, 
  \end{align}  
  for all possible events $\mathcal{S} \subset \mathbb{R}^{n}$ and neighboring datasets $(\mathcal{D}_p, \mathcal{D}'_p)$.
\end{theorem}

\begin{proof}  
  We develop the proof based on the proof in Theorem \ref{thm:laplace}. 
  Specifically, \eqref{laplacian-1} can be rewritten as 
  \begin{subequations}
  \begin{align}
    & \ln \Bigg( \frac{\textbf{pdf} \big( z^{t,e+1}_{p \ell} (\mathcal{D}_p) = \psi \big)}{\textbf{pdf} \big( z^{t,e+1}_{p \ell} (\mathcal{D}'_p ) = \psi \big)} \Bigg) = \ln \Bigg( \frac{\mathcal{N} \big( \tilde{\xi}^{t,e}_p (\psi; \mathcal{D}_p); 0, (\sigma^{t,e}_p)^2 \big) }{ \mathcal{N} \big( \tilde{\xi}^{t,e}_p (\psi; \mathcal{D}'_p); 0, (\sigma^{t,e}_p)^2 \big) } \Bigg) \label{gaussian-1}  \\
    = & \frac{1}{2 (\sigma^{t,e}_p)^2} \Big(\| \tilde{\xi}^{t,e}_p (\psi;\mathcal{D}'_p) \|^2 - \| \tilde{\xi}^{t,e}_p (\psi;\mathcal{D}_p) \|^2 \Big) \nonumber   \\
    = &   \frac{1}{2 (\sigma^{t,e}_p)^2} \Big(\| \tilde{\xi}^{t,e}_p (\psi;\mathcal{D}_p) + f'_p(z_p^{t,e};\mathcal{D}_p) - f'_p(z_p^{t,e};\mathcal{D}'_p)  \|^2 - \| \tilde{\xi}^{t,e}_p (\psi;\mathcal{D}_p) \|^2 \Big)  \label{gaussian-2} \\ 
    = &   \frac{1}{2 (\sigma^{t,e}_p)^2} \Big( 2 \langle \tilde{\xi}^{t,e}_p (\psi;\mathcal{D}_p), f'_p(z_p^{t,e};\mathcal{D}_p) - f'_p(z_p^{t,e};\mathcal{D}_p') \rangle \nonumber \\
    & + \|f'_p(z_p^{t,e};\mathcal{D}_p) - f'_p(z_p^{t,e};\mathcal{D}_p') \|^2 \Big), \nonumber 
  \end{align}
  where $\mathcal{N} \big( \cdot ; 0, (\sigma^{t,e}_p)^2 \big)$ in \eqref{gaussian-1} represents a probability density function of Gaussian distribution with zero mean and a variance $(\sigma^{t,e}_p)^2 := 2 \ln(1.25/\bar{\delta}) (\bar{\Delta}^{t,e}_{p,2} / \bar{\epsilon})^2$, where $\bar{\Delta}^{t,e}_{p,2}$ is from \eqref{redefined_L2sensitivity}, and 
  the equality in \eqref{gaussian-2} holds because of \eqref{correspondence}. 

  Let $x := \langle \tilde{\xi}^{t,e}_p (\psi;\mathcal{D}_p), f'_p(z_p^{t,e};\mathcal{D}_p) - f'_p(z_p^{t,e};\mathcal{D}_p') \rangle$. 
  Then $x \sim \mathcal{N}(0, (\sigma^{t,e}_p)^2 \|f'_p(z_p^{t,e};\mathcal{D}_p) - f'_p(z_p^{t,e};\mathcal{D}_p') \|^2 )$ since each element of $\tilde{\xi}^{t,e}_p$ is drawn from $\mathcal{N}(0, (\sigma^{t,e}_p)^2)$.    
  Let $y := x/(\sigma^{t,e}_p)^2 + \|f'_p(z_p^{t,e};\mathcal{D}_p) - f'_p(z_p^{t,e};\mathcal{D}_p') \|^2 / (2 (\sigma^{t,e}_p)^2)$ be the privacy loss random variable defined in \eqref{gaussian-1}. Then the mean and variance of $y$ are given by
  \begin{align}
    & \frac{\|f'_p(z_p^{t,e};\mathcal{D}_p) - f'_p(z_p^{t,e};\mathcal{D}_p') \|^2}{2 (\sigma^{t,e}_p)^2} \text { and } \frac{\|f'_p(z_p^{t,e};\mathcal{D}_p) - f'_p(z_p^{t,e};\mathcal{D}_p') \|^2}{(\sigma^{t,e}_p)^2}, \label{privacy_loss_rv_mean_var}
  \end{align}
  respectively.  
  Let $z \sim \mathcal{N}(0,1)$. Then the privacy loss random variable can be rewritten as 
  \begin{align}
    y := \frac{\|f'_p(z_p^{t,e};\mathcal{D}_p) - f'_p(z_p^{t,e};\mathcal{D}_p') \|}{\sigma^{t,e}_p} z + \frac{\|f'_p(z_p^{t,e};\mathcal{D}_p) - f'_p(z_p^{t,e};\mathcal{D}_p') \|^2}{2 (\sigma^{t,e}_p)^2}.
  \end{align}
  It suffices to prove the following:
  \begin{align}
    \mathbb{P}\big( |y| > \bar{\epsilon} \big) \leq \bar{\delta}, \label{eps_delta_dp}
  \end{align}
  which is equivalent to the definition of $(\bar{\epsilon},\bar{\delta})$-DP according to \cite[Lemma 3.17]{dwork2014algorithmic}.
  To this end, we derive an upper bound on the left-hand side of \eqref{eps_delta_dp} as 
  \begin{align}
    & \mathbb{P}\big( |y| > \bar{\epsilon} \big)  \nonumber \\
    = \ & \mathbb{P}\Big( |z| >\frac{\bar{\epsilon} \sigma^{t,e}_p}{\|f'_p(z_p^{t,e};\mathcal{D}_p) - f'_p(z_p^{t,e};\mathcal{D}_p') \|}   - \frac{\|f'_p(z_p^{t,e};\mathcal{D}_p) - f'_p(z_p^{t,e};\mathcal{D}_p') \|}{2 \sigma^{t,e}_p}\Big) \nonumber  \\
    \leq \ & \mathbb{P}\Big( |z| >  \sqrt{2 \ln(1.25/\bar{\delta})}   - \frac{\bar{\epsilon} }{2 \sqrt{2 \ln(1.25/\bar{\delta})}  }\Big)  \label{def_delta} \\
    \leq \ & \frac{2}{\sqrt{2 \pi}} \exp \Big(\frac{-1}{2} \Big(\sqrt{2 \ln(1.25/\bar{\delta})}   - \frac{\bar{\epsilon} }{2 \sqrt{2 \ln(1.25/\bar{\delta})}  } \Big)^2 \Big)  \label{tail_bound}  \\
    \leq \ & \bar{\delta},
  \end{align}
  where the inequality \eqref{def_delta} holds because of the definition of $\sigma^{t,e}_p$;
  the inequality \eqref{tail_bound} holds because of the tail bound of $\mathcal{N}(0,1)$ is given by $\mathbb{P}(|z| > r)  \leq \frac{2}{\sqrt{2\pi}} \exp(\frac{-r^2}{2})$; and
  the inequality holds because if $t=\sqrt{2 \ln(1.25/\bar{\delta})}$, then
  $(t - \bar{\epsilon}/t)^2 \geq 2 \ln (\frac{2}{\bar{\delta}} \frac{1}{\sqrt{2\pi}})$ for $\bar{\epsilon} \leq 1$ (see \cite[Appendix A]{dwork2014algorithmic} for more details). 
  \end{subequations}
\end{proof}

\subsection{Total privacy leakage}

The results from Sections \ref{sec:randomized_algo} and \ref{sec:mechanisms} allow us to use Lemma \ref{lemma:kifer} to prove that \eqref{DPADMM-2-Prox} preserves data privacy in a DP manner.

\begin{theorem} \label{thm:privacy}
  For fixed $t \in [T]$ and $e \in [E]$, each agent $p \in [P]$ can achieve (i) $(\bar{\epsilon}, 0)$-DP by solving \eqref{DPADMM-2-Prox}, where $\tilde{\xi}^{t,e}_p$ is drawn from a Laplace distribution with zero mean and a variance $2(\bar{\Delta}^{t,e}_{p,1}/\bar{\epsilon})^2$ and $\bar{\Delta}^{t,e}_{p,1}$ is from \eqref{redefined_L1sensitivity} and (ii) $(\bar{\epsilon}, \bar{\delta})$-DP with $\bar{\delta}>0$ by solving \eqref{DPADMM-2-Prox}, where $\tilde{\xi}^{t,e}_p$ is drawn from a Gaussian distribution with zero mean and a variance $2\ln(1.25/\bar{\delta})(\bar{\Delta}^{t,e}_{p,2}/\bar{\epsilon})^2$ and $\bar{\Delta}^{t,e}_{p,2}$ is from \eqref{redefined_L2sensitivity}.
\end{theorem} 

We note that the average randomized output $z^{t}_p$ in line 20 of Algorithm \ref{algo:DP-IADMM-Prox} is the communicated information among agents while the locally updated outputs $\{z_p^{t,1}, \ldots, z_p^{t,E} \}$ are not.
Nevertheless, we randomize all locally updated outputs to prevent data leakage from any one of them, which considers an extreme case when an adversary knows the average randomized output $z^{t}_p$ as well as any $E-1$ locally updated outputs (e.g.,  $\{z_p^{t,1}, \ldots, z_p^{t,E-1} \}$) that can be utilized to compute the  last unknown output (e.g., $z_p^{t,E}$).

Now we study a total privacy leakage of Algorithm \ref{algo:DP-IADMM-Prox}, which can be considered as a $TE$-fold adaptive algorithm.
This involves an extreme situation in which every communication round leaks the output.
According to the existing composition theorem \cite[Theorem 3.16]{dwork2014algorithmic}, Algorithm \ref{algo:DP-IADMM-Prox} guarantees $(TE \bar{\epsilon}, TE \bar{\delta})$-DP because the proposed mechanism \eqref{ADMM-2-Prox-log} guarantees $(\bar{\epsilon}, \bar{\delta})$-DP where $\bar{\delta} \geq 0$. 
This composition theorem implies that data privacy becomes weaker as $T$ and $E$ increase.
We note that this theorem holds for both Laplace and Gaussian mechanisms.

For the Gaussian mechanism, on the other hand, Abadi et al.~\cite{abadi2016deep} proposed a stronger composition theorem based on the moments accountant method, which states that a $T$-adaptive algorithm that guarantees $(\bar{\epsilon}, \bar{\delta})$-DP with $\bar{\delta} > 0$ for every iteration by the existing Gaussian mechanism in Lemma \ref{lem:existing-gaussian} is $(\mathcal{O}(\sqrt{T}) \bar{\epsilon}, \bar{\delta})$-DP, which is a tighter bound than $(T \bar{\epsilon}, T\bar{\delta})$-DP from the composition theorem \cite[Theorem 3.16]{dwork2014algorithmic}.
In the following theorem, we utilize the moments accountant method in \cite{abadi2016deep} to show that our algorithm with the Gaussian mechanism proposed in Theorem \ref{thm:gaussian} provides a tighter bound. 
\begin{theorem}
Consider Algorithm \ref{algo:DP-IADMM-Prox} where \eqref{DPADMM-2-Prox} in line 18 is replaced with \eqref{ADMM-2-Prox-log} where $\ell$ is sufficiently large such that $z^{t,e+1}_{p \ell} \approx z^{t,e+1}_p $ as shown in Theorem \ref{thm:pointwise_convergence}.
This algorithm guarantees $(\sqrt{\frac{TE \ln (1/\bar{\delta})}{\ln(1.25/\bar{\delta})}} \bar{\epsilon}, \bar{\delta} )$-DP if \eqref{ADMM-2-Prox-log} is $(\bar{\epsilon}, \bar{\delta})$-DP, as shown in Theorem \ref{thm:gaussian}.
\end{theorem}
\begin{proof}
\begin{subequations}
We develop the proof based on \cite[Theorem 2]{abadi2016deep}, which is restated in Lemma \ref{lem:composition} in this paper. 
Since we consider an algorithm composed of a sequence of $TE$ adaptive Gaussian mechanisms, we rewrite \eqref{lem:composition-2} within our context, for fixed $t$, $e$, $p$, and sufficiently large $\ell$, as
\begin{align}
  \alpha^{t,e}_p (\tau) := & \max_{\forall (\mathcal{D}_p, \mathcal{D}_p')} \ln \big( \mathbb{E} [ x^{\tau} ] \big),  \ \ \ \
  x := \frac{\textbf{pdf} \big( z^{t,e+1}_{p \ell} (\mathcal{D}_p) = \psi \big)}{\textbf{pdf} \big( z^{t,e+1}_{p \ell} (\mathcal{D}'_p ) = \psi \big)},
\end{align} 
where 
$(\mathcal{D}_p, \mathcal{D}_p')$ are two neighboring datasets, 
$z^{t,e+1}_{p \ell}$ is the optimal solution of \eqref{ADMM-2-Prox-log},
and
$x$ is a random variable that follows a log-normal distribution because in Theorem \ref{thm:gaussian} we show that $\ln(x)$ is a random variable that follows a Gaussian distribution whose mean $\mu$ and variance $\sigma^2$ are presented in \eqref{privacy_loss_rv_mean_var}.
Noting that the $\tau$th moment of the random variable $x$ is given by $\mathbb{E}[x^{\tau}] = \exp( \tau \mu + (\tau^2/2) \sigma^2 )$, we derive
\begin{align}
  \alpha^{t,e}_p(\tau) = & \max_{\forall (\mathcal{D}_p, \mathcal{D}_p')} \ln  ( \mathbb{E}[ x^{\tau}] )  =  \max_{\forall (\mathcal{D}_p, \mathcal{D}_p')}    (\tau+\tau^2) \Big( \frac{\|f'_p(z_p^{t,e};\mathcal{D}_p) - f'_p(z_p^{t,e};\mathcal{D}_p') \|^2}{2 (\sigma^{t,e}_p)^2}  \Big) \nonumber \\
  = &\frac{\tau (\tau+1) (\bar{\Delta}^{t,e}_{p,2})^2}{2 (\sigma^{t,e}_{p})^2} = \frac{\tau(\tau+1) \bar{\epsilon}^2}{4 \ln (1.25/\bar{\delta})}, \label{def_moment}
\end{align}
where the derivation in \eqref{def_moment} holds because of \eqref{redefined_L2sensitivity} and a variance $(\sigma^{t,e}_p)^2 := 2 \ln(1.25/\bar{\delta}) (\bar{\Delta}^{t,e}_{p,2} / \bar{\epsilon})^2$.

Next we restate \eqref{lem:composition-1} within our context as follows:
\begin{align}
  \ln(\bar{\delta}) = \ &  \min_{\tau \in \mathbb{Z}_+} \sum_{t=1}^T \sum_{e=1}^E \alpha^{t,e}_p(\tau)  - \tau \epsilon' = \min_{\tau \in \mathbb{Z}_+} \tau ( a \tau + (a-\epsilon') ), \label{tau_problem} \\
  \text{ where } a := \ &  \frac{TE \bar{\epsilon}^2}{4 \ln (1.25/\bar{\delta})}.
\end{align}
We note that the objective function $\tau ( a \tau + (a-\epsilon') )$ in \eqref{tau_problem} is strictly convex because $a > 0$, whose value is zero at $\tau = 0$ and $\tau = -(a-\epsilon')/a$. 
Since $\ln (\bar{\delta}) < 0$ as $\bar{\delta} \in (0,1)$, the optimal value should be negative; and this is possible when $\tau = 1$ is a feasible solution to the optimization problem in \eqref{tau_problem}.
To this end, we should have
\begin{align}
\frac{\epsilon'-a}{a} > 1  \Leftrightarrow  \epsilon' > \frac{TE \bar{\epsilon}^2}{2 \ln (1.25/\bar{\delta})}. \label{thm:composition-0}
\end{align} 
Now consider
\begin{align}
  \ln(\bar{\delta}) = \ &  \min_{\tau \in \mathbb{Z}_+} \sum_{t=1}^T \sum_{e=1}^E \alpha^{t,e}_p(\tau)  - \tau \epsilon'  \nonumber \\
  \geq \ & \min_{\tau \in \mathbb{R}} \sum_{t=1}^T \sum_{e=1}^E \alpha^{t,e}_p(\tau) - \tau \epsilon' = - \frac{(\epsilon')^2 \ln(1.25/\bar{\delta})}{TE \bar{\epsilon}^2} - \frac{TE \bar{\epsilon}^2}{16 \ln(1.25/\bar{\delta})} + \frac{\epsilon'}{2} \label{thm:composition-1} \\
  \geq \ & - \frac{(\epsilon')^2 \ln(1.25/\bar{\delta})}{TE \bar{\epsilon}^2} - \frac{\epsilon'}{8}  + \frac{\epsilon'}{2}  \geq  - \frac{(\epsilon')^2 \ln(1.25/\bar{\delta})}{TE \bar{\epsilon}^2},\label{thm:composition-2}
\end{align}
where 
the inequality in \eqref{thm:composition-1} holds as the integrality restriction is relaxed;
the equality in \eqref{thm:composition-1} holds because $\tau^* = (\epsilon'-b) / (2b)$, where $b = \bar{\epsilon}^2 / (4 \ln (1.25/\bar{\delta}))$, by the optimality condition; and
the first inequality in \eqref{thm:composition-2} holds because of \eqref{thm:composition-0}.
Therefore we have
\begin{align}
  \bar{\epsilon} \sqrt{ \frac{TE \ln (1/ \bar{\delta})}{\ln(1.25/\bar{\delta})} } \leq \epsilon'.
\end{align}
\end{subequations}
\end{proof}

\section{Convergence analysis} \label{sec:convergence}
In this section we show that a sequence of iterates generated by Algorithm \ref{algo:DP-IADMM-Prox} converges to an optimal solution of \eqref{model:dist} in expectation and that its associated convergence rate is sublinear without any error bound, unlike \cite{huang2020differentially} that results in a nontrivial positive error bound as stated in Section \ref{sec:baseline}.
Moreover, our results characterize the impact of privacy budget parameter $\bar{\epsilon}$ and the number of local updates $E$ as well as the number of ADMM rounds $T$ on the convergence rates.

\subsection{Assumptions}
We make the following assumptions:
\begin{assumption}\label{assump:convergence}  
  \noindent
  \begin{enumerate}
    \item[(i)] The ADMM penalty parameter $\rho^t$ is nondecreasing and bounded; that is, $\rho^1 \leq \rho^2 \leq \ldots \leq \rho^T \leq \rho^{\text{max}}$.
    \item[(ii)] There exists $\gamma > 0$ such that $\gamma \geq 2\| \lambda^*\|$, where $\lambda^*$ is a dual optimal. 
    \item[(iii)] The local objective function $f_p$ is $H$-Lipschitz over a set $\mathcal{W}_p$ with respect to the Euclidean norm.    
  \end{enumerate}  
\end{assumption}
Assumption \ref{assump:convergence} is typically used for the convergence analysis of ADMM (see Chapter 15 of \cite{beck2017first}). 
We adopt the assumptions because the proposed algorithm is a variant of ADMM.
Specifically, Algorithm \ref{algo:DP-IADMM-Prox} with $E=1$ and $\tilde{\xi}^{t,1}_p=0$ for all $t$ is equivalent to the linearized ADMM.
Under Assumption \ref{assump:convergence}, we show that the rate of convergence in expectation is $\mathcal{O}( 1 / (\bar{\epsilon}^2 E \sqrt{T}))$ (resp., $\mathcal{O}( 1 / (\bar{\epsilon} E \sqrt{T}))$) when the local objective function $f_p$ is smooth (resp. nonsmooth).

In addition, we show that the rate becomes $\mathcal{O}( 1 / (\bar{\epsilon}^2 E T))$ when the local objective function $f_p$ is strongly convex. 
This result requires additional assumptions:
\begin{assumption}\label{assump:convergence_stronglyconvex}
\noindent
\begin{enumerate}
\item[(i)] There exists $\gamma > 0$ such that $\| \lambda^t \| \leq \gamma, \ \forall t$.
\item[(ii)] $\rho^t \leq \frac{t}{t-1} \rho^{t-1}, \ \forall t$.
\end{enumerate}
\end{assumption}
Assumption \ref{assump:convergence_stronglyconvex} (i) can be strict in practice. As indicated in \cite{azadi2014towards}, however, it can be considered as a price that we have to pay for faster convergence (see \cite[Assumption 3]{azadi2014towards}).
Assumption \ref{assump:convergence_stronglyconvex} (ii) is not strict since it is satisfied by a constant penalty (i.e., $\rho^t = \rho, \ \forall t$) that is commonly considered in the literature (e.g.,~\cite{beck2017first}).

\subsection{Main results}
In this section we present the convergence analysis results in Theorems \ref{thm:smooth}, \ref{thm:nonsmooth}, and \ref{thm:strong} when the local objective function $f_p$ is smooth, nonsmooth, and strongly convex, respectively; we relegate all the detailed proofs to the Appendix.

\begin{theorem} \label{thm:smooth}
  Suppose that Assumption \ref{assump:convergence} holds, the local objective function $f_p$ is $L$-smooth convex, and the parameter $\eta^t$ of \eqref{DPADMM-2-Prox} is updated as
  \begin{subequations}    
  \begin{align}
    \eta^t  = \frac{1}{L + \sqrt{t}/\bar{\epsilon}}, \ \forall t, \label{smooth_proximity}
  \end{align}
  where $\bar{\epsilon} > 0$ is a privacy budget parameter.
  Then we have  
  \begin{align}
    & \mathbb{E} \Big[ F(z^{(T)}) - F(z^*) + \gamma \| Aw^{(T)} - z^{(T)} \| \Big]   \label{smooth_inequality} \\
    & \leq \frac{nPU_3 + U_2^2/(2E) }{ \bar{\epsilon} \sqrt{T} } + \frac{U_2^2  (\rho^{\text{max}} +L/E) +  (\gamma + \|\lambda^1\|)^2 / \rho^1 }{2T}, \nonumber
  \end{align}    
  where  
  \begin{align}    
    & z := [z_1^{\top}, z_2^{\top}, \ldots, z_P^{\top}]^{\top}, \ \ z^{(T)} :=  \frac{1}{TE} \sum_{t=1}^T \sum_{e=1}^E z^{t,e+1}, \ \ F(z) :=  \sum_{p=1}^P f_p(z_p), \label{Common_Notation} \\ 
    & A^{\top} := \begin{bmatrix}
        \mathbb{I}_J \ \cdots \  \mathbb{I}_J        
      \end{bmatrix}_{J \times PJ}, \ \ w^{(T)} :=  \frac{1}{T} \sum_{t=1}^T w^{t+1}, \nonumber \\
    & U_1 := \max_{u \in \mathcal{W}_p, p \in [P]} \| f'_p(u) \|, \ \ U_2 := \max_{u, v \in \mathcal{W}_p, p \in [P] } \|u-v\|, \nonumber \\
    & U_3 := 
    \begin{cases}
    \max_{p \in [P]}  2 ( \bar{\Delta}^{t,e}_{p,1} )^2  & \text{ for the Laplace mechanism} \\
    \max_{p \in [P]}  2 \ln(1.25/\bar{\delta}) ( \bar{\Delta}^{t,e}_{p,2}  )^2  &   \text{ for the Gaussian mechanism} 
    \end{cases}    
    \nonumber
  \end{align}  
  and $z^*$ is an optimal solution of \eqref{model:dist_1}.
  \end{subequations}  
\end{theorem}
\begin{proof}
See Appendix \ref{apx-thm:smooth}.
\end{proof}
The inequality \eqref{smooth_inequality} implies that the rate of convergence in expectation is 
$\mathcal{O}(1/ (\bar{\epsilon} \sqrt{T}))$ for $\bar{\epsilon} \in (0, \infty)$, 
while in a nonprivate setting it is $\mathcal{O}(1/T)$ because the first term in \eqref{smooth_inequality} is zero when $\bar{\epsilon}=\infty$.

\begin{theorem} \label{thm:nonsmooth}
Suppose that Assumption \ref{assump:convergence} holds, the local objective function $f_p$ is a nonsmooth convex function, and the parameter $\eta^t$ of \eqref{DPADMM-2-Prox} is updated as
\begin{subequations}
\begin{align}
\eta^{t} = \frac{1}{\sqrt{t}}, \ \forall t.  \label{nonsmooth_proximity}
\end{align}
Then we have
\begin{align}
& \mathbb{E} \Big[ F(z^{(T)}) - F(z^*) + \gamma \| Aw^{(T)} - z^{(T)} \| \Big] \leq  \label{nonsmooth_inequality} \\
&  \frac{ nPU_3 /\bar{\epsilon}^2 + PU_1^2 + U_2^2 / (2E) }{ \sqrt{T} } + \frac{U_2^2 \rho^{\text{max}} +  (\gamma + \|\lambda^1\|)^2 / \rho^1 + 2\gamma U_2}{2T}, \nonumber
\end{align}  
where $z$, $F(z)$, $A$, $w^{(T)}$, $U_1$, $U_2$, $U_3$ are from \eqref{Common_Notation} while we define
\begin{align}
z^{(T)} :=  \frac{1}{TE} \sum_{t=1}^T \sum_{e=1}^E z^{t,e}. \label{nonsmooth_z}
\end{align}
\end{subequations}
\end{theorem}
\begin{proof}
  See Appendix \ref{apx-thm:nonsmooth}.
\end{proof}
The inequality \eqref{nonsmooth_inequality} implies that the rate of convergence in expectation is 
$\mathcal{O}(1/ (\bar{\epsilon}^2 \sqrt{T}))$ for $\bar{\epsilon} \in (0, \infty)$, 
while in a nonprivate setting it is $\mathcal{O}(1/\sqrt{T})$ because the first term in \eqref{nonsmooth_inequality} is $(PU_1^2 + U_2^2/(2E))/\sqrt{T}$ when $\bar{\epsilon}=\infty$. 
 
\begin{theorem} \label{thm:strong}
Suppose that Assumptions \ref{assump:convergence} and \ref{assump:convergence_stronglyconvex} hold, the local objective function $f_p$ is $\alpha$-strongly convex, and the parameter $\eta^t$ of \eqref{DPADMM-2-Prox} is updated as
\begin{subequations}
\begin{align}
\eta^{t} = \frac{2}{\alpha(t+2)}, \ \forall t. 
\end{align}
Then we have
\begin{align}
  & \mathbb{E} \Big[ F(z^{(T)}) - F(z^*) + \gamma \| Aw^{(T)} - z^{(T)} \| \Big] \leq  \label{strong_inequality} \\
  &\frac{1}{T+1} \Big\{2 U_2 \gamma + U_2^2 \rho^{\text{max}} + 4 \gamma^2 / \rho^1 + \alpha U_2^2/(2E)  + 2P(U_1^2 + nU_3/\bar{\epsilon}^2)/\alpha \Big\}, \nonumber
\end{align}
where $z$, $F(z)$, $A$, $w^{(T)}$, $U_1$, $U_2$, $U_3$ are from \eqref{Common_Notation} while we define
\begin{align}
  & z^{(T)} := \frac{2}{T(T+1)} \sum_{t=1}^T t  (\frac{1}{E} \sum_{e=1}^E z^{t,e}), \ \ w^{(T)} := \frac{2}{T(T+1)} \sum_{t=1}^T t w^{t+1}. \label{strong_z}   
\end{align}  
\end{subequations}
\end{theorem}
\begin{proof}
  See Appendix \ref{apx-thm:strong}. 
\end{proof}
The inequality \eqref{strong_inequality} derived under Assumptions \ref{assump:convergence} and \ref{assump:convergence_stronglyconvex} implies that the rate of convergence in expectation is 
$\mathcal{O}(1/ (\bar{\epsilon}^2 T))$ for $\bar{\epsilon} \in (0, \infty)$, 
while in a nonprivate setting it is $\mathcal{O}(1/T)$ when $\bar{\epsilon}=\infty$. 
  
From the convergence analysis we note that increasing the number $E$ of local updates decreases the values on the right-hand side of \eqref{smooth_inequality}, \eqref{nonsmooth_inequality}, and \eqref{strong_inequality}.
This implies that the gap between $F(z^{(T)})$ and $F(z^*)$ can become smaller by increasing $E$ for fixed $T$.
This may result in better performance by introducing the multiple local updates, thus reducing communication costs, as will be numerically demonstrated in Section \ref{sec:experiments}.

\section{Numerical demonstration} \label{sec:experiments}
In this section we  numerically demonstrate the effectiveness of the proposed algorithm based on the objective perturbation and multiple local updates for distributed optimization models with constraints.
Specifically, we compare our approach with the popular output perturbation method, which has been widely used in the literature \cite{dvorkin2020differentially, huang2019dp, huang2020differentially}.
For ease of exposition, we denote by $\texttt{ObjL}$ and $\texttt{ObjG}$ the proposed algorithm with Laplace and Gaussian mechanisms, respectively, which will be compared with $\texttt{OutL}$ and $\texttt{OutG}$, the output perturbation methods.
For the comparison of the algorithms, we consider two applications: distributed control of power flow in Section \ref{sec:numerical_power_flow} and federated learning in Section \ref{sec:numerical_fl}.

\subsection{Privacy-preserving distributed control of power flow} \label{sec:numerical_power_flow}
In power systems, an increasing penetration of distributed energy resources in the network has motivated the distributed control of power flow.
Such a problem is a form of \eqref{model:dist}.
We specifically consider that the power network can be decomposed into several zones, each of which is controlled by the agent $p \in [P]$. 
These agents cooperate to find a global power flow solution $w^*$ by iteratively providing a local power flow solution resulting from the model composed of local objective function $f_p$ and local constraints $\mathcal{W}_p$.

In particular we consider an optimization problem that determines a distributed control of power flow such that the deviation of power balance is minimized as considered in \cite{holzer2021grid, dandurand2021bilevel, kim2016data}.
Such a problem can be formulated with the following local objective function:
\begin{align}
f_p (z_p) := \sum_{i \in \mathcal{N}_p} (a_{pi}^{\top} z_{pi} + d_i)^2, \ \forall p \in [P] , \label{localobj}
\end{align}
where $\mathcal{N}_p$ is a set of buses controlled by agent $p$, $d_i$ is a given demand data, $z_{pi}$ is a vector composed of power flow and generation variables connected to a node $i \in \mathcal{N}_p$, and $a_{pi}$ is a given coefficient vector whose element has a value in $[-1,1]$.
Note that the  objective function is derived from the power balance equations $a_{pi}^{\top} z_{pi} = d_i$ for all $i \in \mathcal{N}_p$. 
The local constraint $\mathcal{W}_p$ is composed of convex relaxation constraints for observing some physical laws and operational constraints.
Detailed formulation is presented in Appendix \ref{apx:distpf}.

In this setting, agents are not required to share electric power loads consumed at their zones, but it is still possible to reverse-engineer local solutions communicated between agents to estimate the sensitive power load data \cite{ryu2021privacy}.
To protect data, one can utilize DP algorithms that admit an inevitable trade-off between data privacy and solution quality.
We aim to show that our algorithms \texttt{ObjG} and \texttt{ObjL} provide better solution quality than do the existing algorithms \texttt{OutG} and \texttt{OutL} while providing the same level of data privacy.

\subsubsection{Experimental settings} \label{sec:1_exp_settings}
For the power network instances, we consider case 14 and case 118, each of which is decomposed into three zones as in \cite{ryu2021privacy}.
For the ADMM parameters, we set $\rho^t = 100$ for all $t$ and $\eta^t = 1/\sqrt{t}$.
We consider various $\bar{\epsilon} \in \{0.05, 0.1, 0.5, 1\}$ where smaller $\bar{\epsilon}$ ensures stronger data privacy, as described in Definition \ref{def:differential_privacy_1}. 
For the Gaussian mechanism, we fix $\bar{\delta} = 10^{-2}$ because $\bar{\delta} < n^{-2}$ is typically chosen where $n$ is the number of data points.
For the algorithms based on output perturbation (i.e., \texttt{OutG} and \texttt{OutL}), we compute the sensitivity $\bar{\Delta}_1$ from \eqref{L1_sensitivity} and $\bar{\Delta}_2$ from \eqref{L2_sensitivity} by considering the worst-case scenario as described in \cite[Proposition 1]{dvorkin2020differentially}, resulting in $\bar{\Delta}_1 = \bar{\Delta}_2 = \beta$, where $\beta$ measures the adjacency of two neighboring datasets, namely, $\mathcal{D}_p := \{d_i\}_{i \in \mathcal{N}_z}$ and $\mathcal{D}_p':= \{d_i'\}_{i \in \mathcal{N}_z}$, which are $\beta$-adjacent datasets if  $\exists i \text{ s.t. } |d_i - d_i'| \leq \beta, \ d_j = d_j', \forall j \neq i$.
For a fair comparison, we also compute the worst-case sensitivity for our algorithms \texttt{ObjG} and \texttt{ObjL}, resulting in $\bar{\Delta}_1 = \bar{\Delta}_2 = 2\beta$, which is twice  the sensitivity used in \texttt{OutG} and \texttt{OutL} based on \eqref{redefined_L1sensitivity}, \eqref{redefined_L2sensitivity}, and \eqref{localobj}.
In this experiment we set $\beta = 1 \%$ as in \cite{dvorkin2020differentially}.
For all experiments, we solve optimization models by Ipopt \cite{wachter2006implementation} via Julia 1.8.0 on Bebop, a 1024-node computing cluster at Argonne National Laboratory.
Each computing node has 36 cores with Intel Xeon E5-2695v4 processors and 128 GB DDR4 of memory.


\begin{figure*}[!ht]  
\centering
\begin{subfigure}[b]{0.23\textwidth}
\centering        
\includegraphics[width=\textwidth]{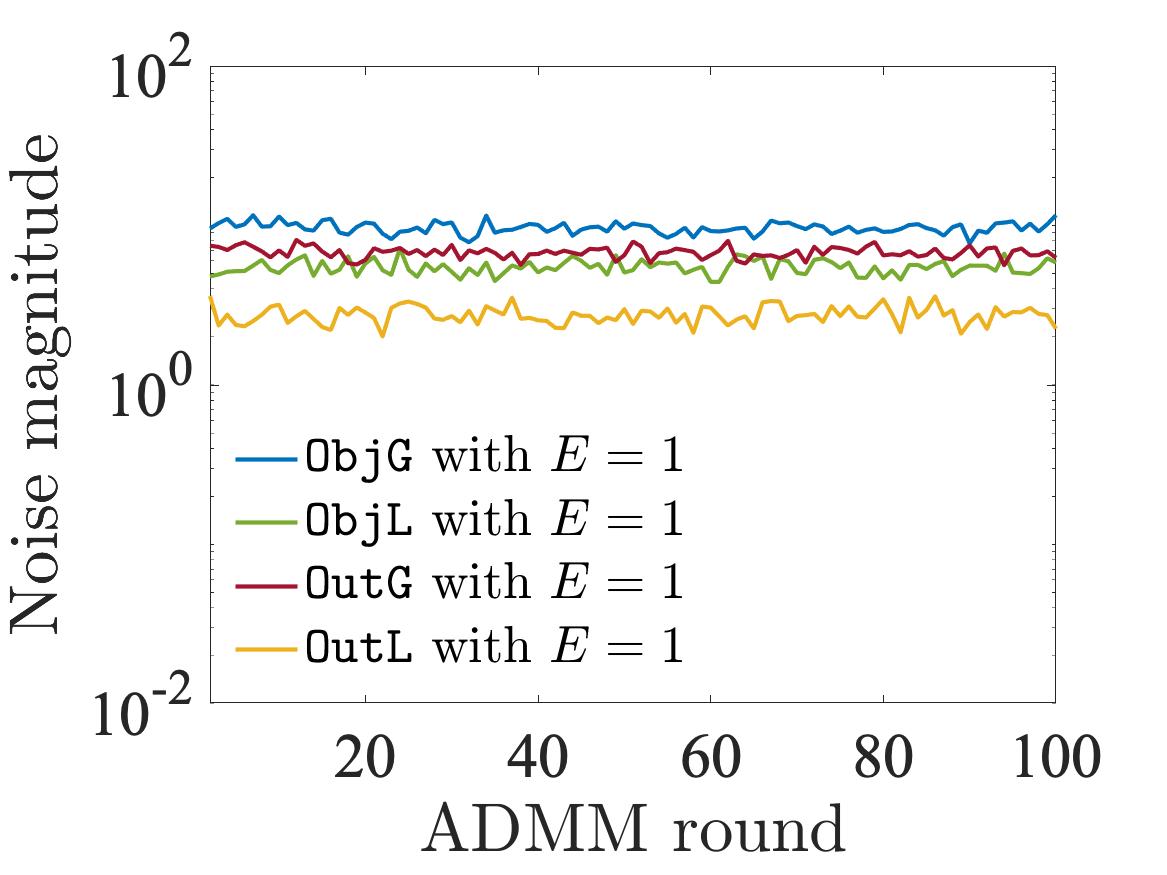}  
\includegraphics[width=\textwidth]{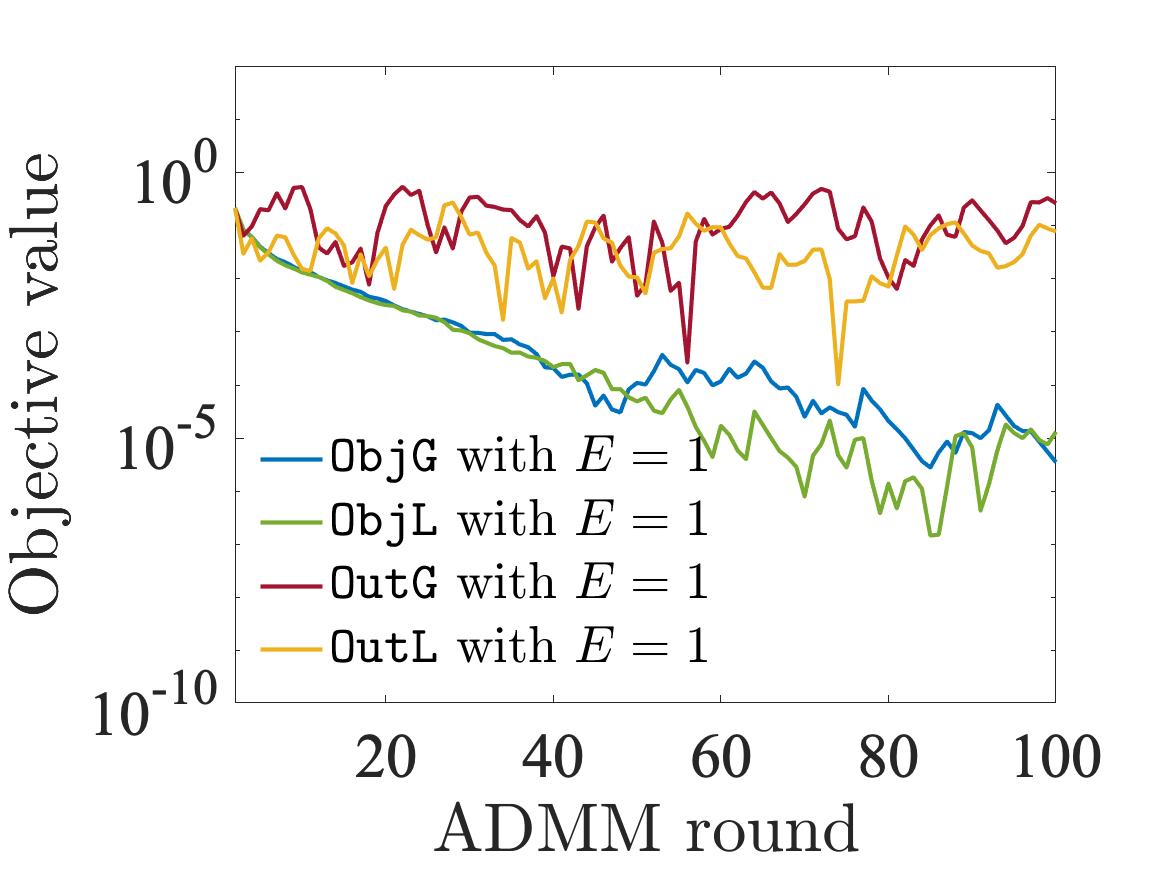}   
\includegraphics[width=\textwidth]{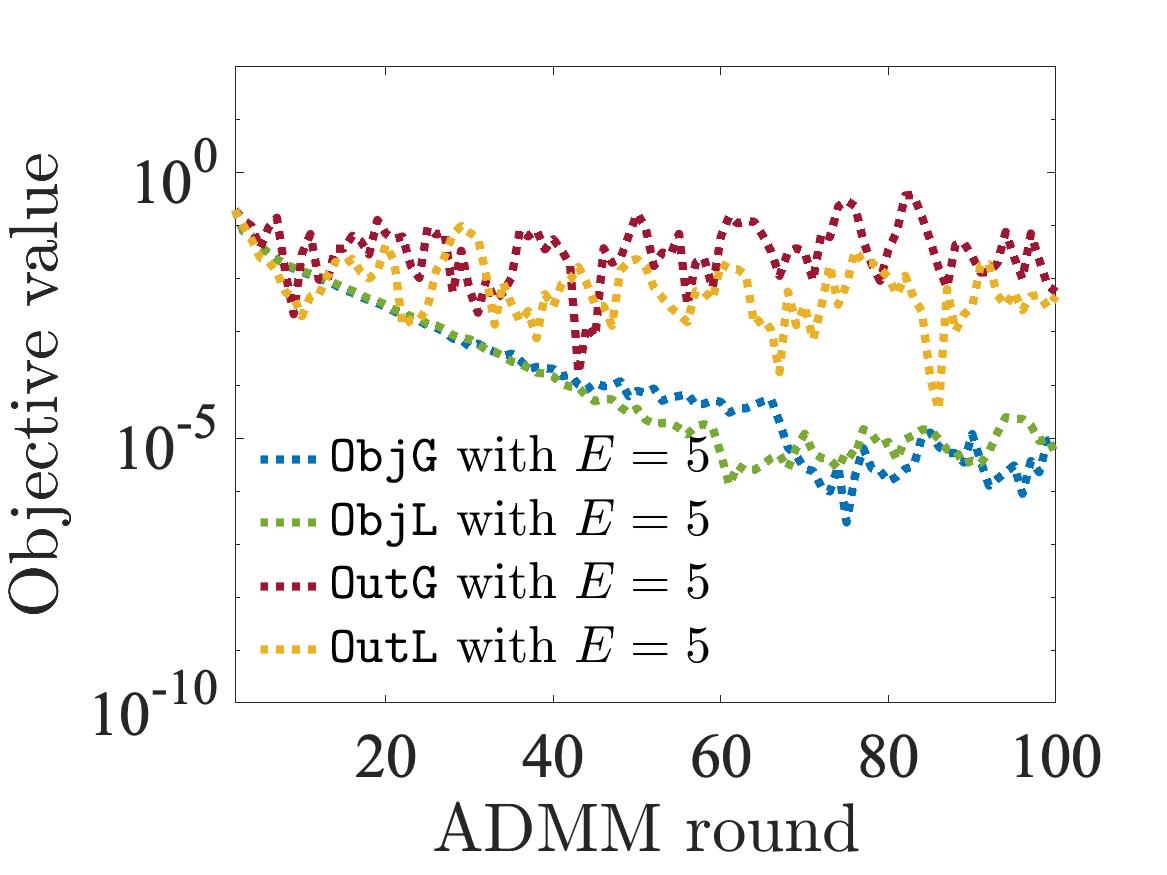}   
\includegraphics[width=\textwidth]{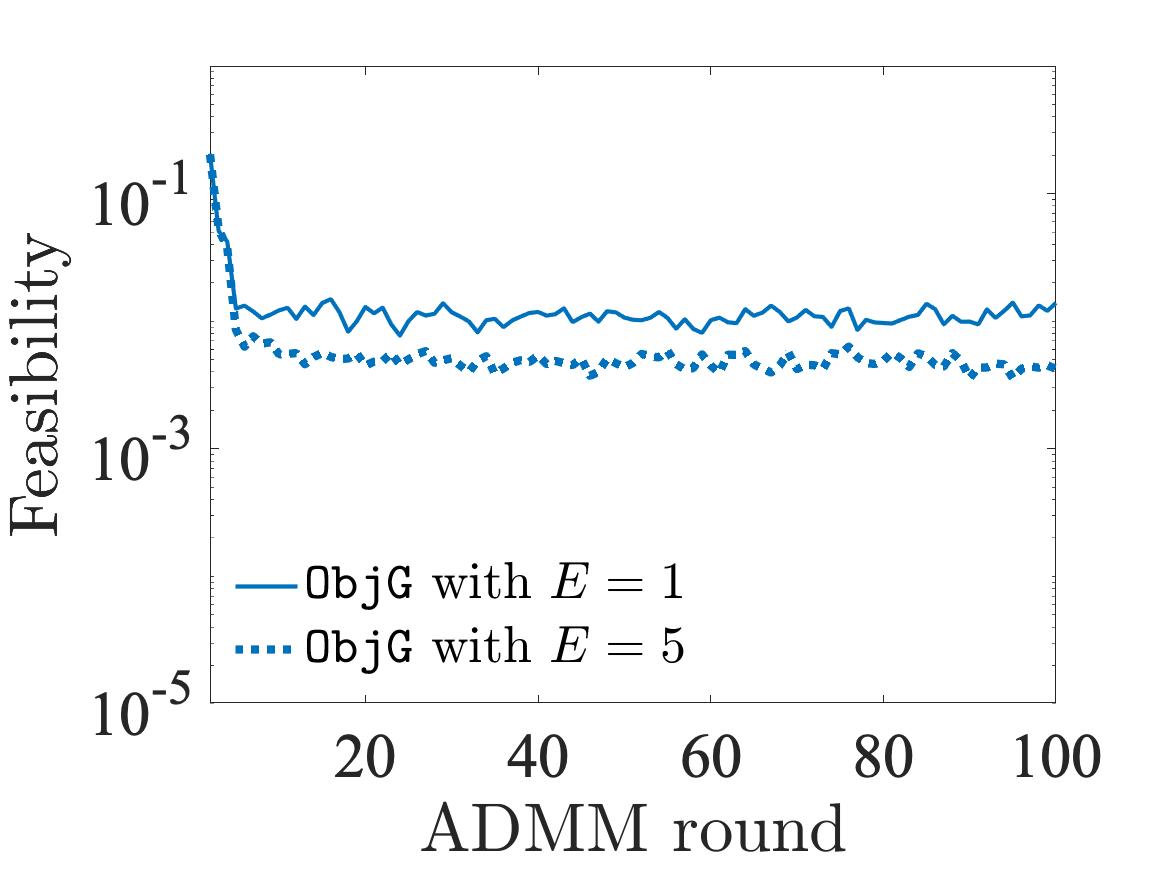}  
\caption{$\bar{\epsilon}=0.05$}
\end{subfigure}
\begin{subfigure}[b]{0.23\textwidth}
\centering        
\includegraphics[width=\textwidth]{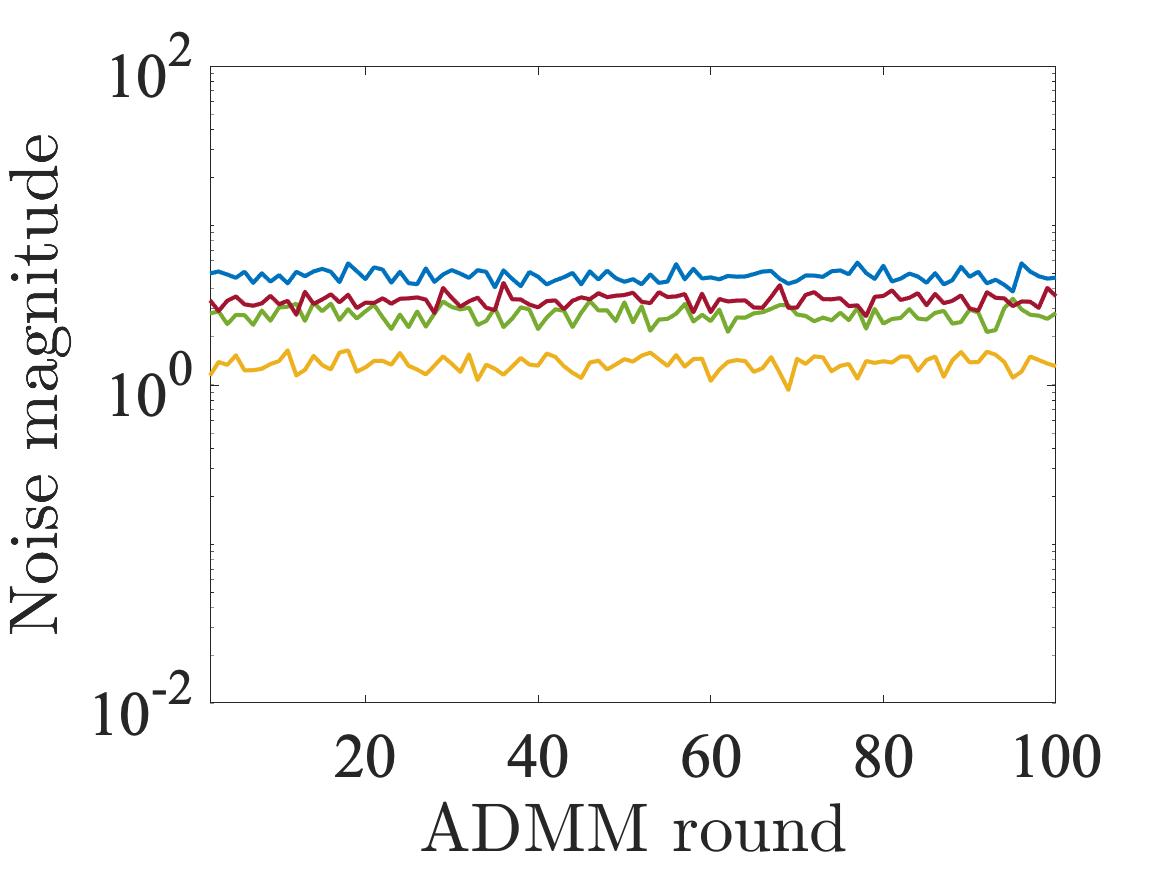}  
\includegraphics[width=\textwidth]{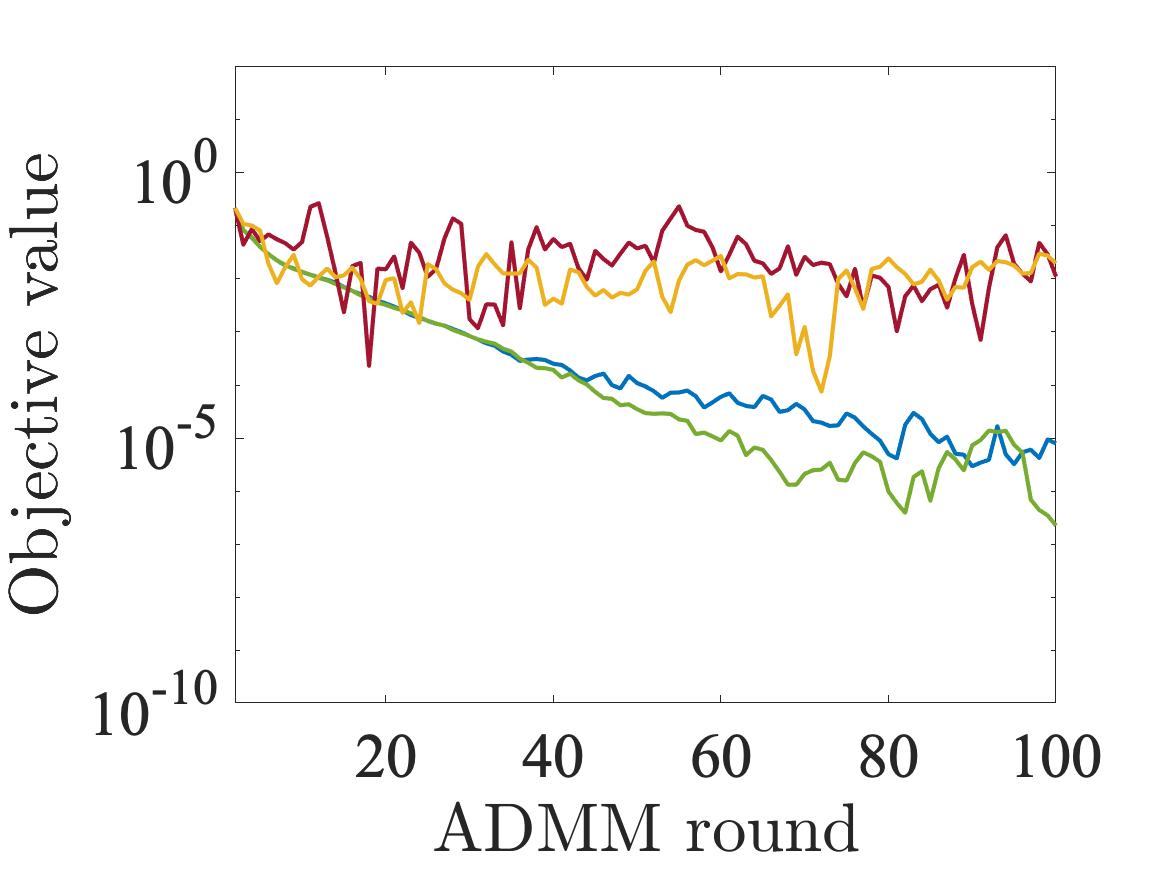}   
\includegraphics[width=\textwidth]{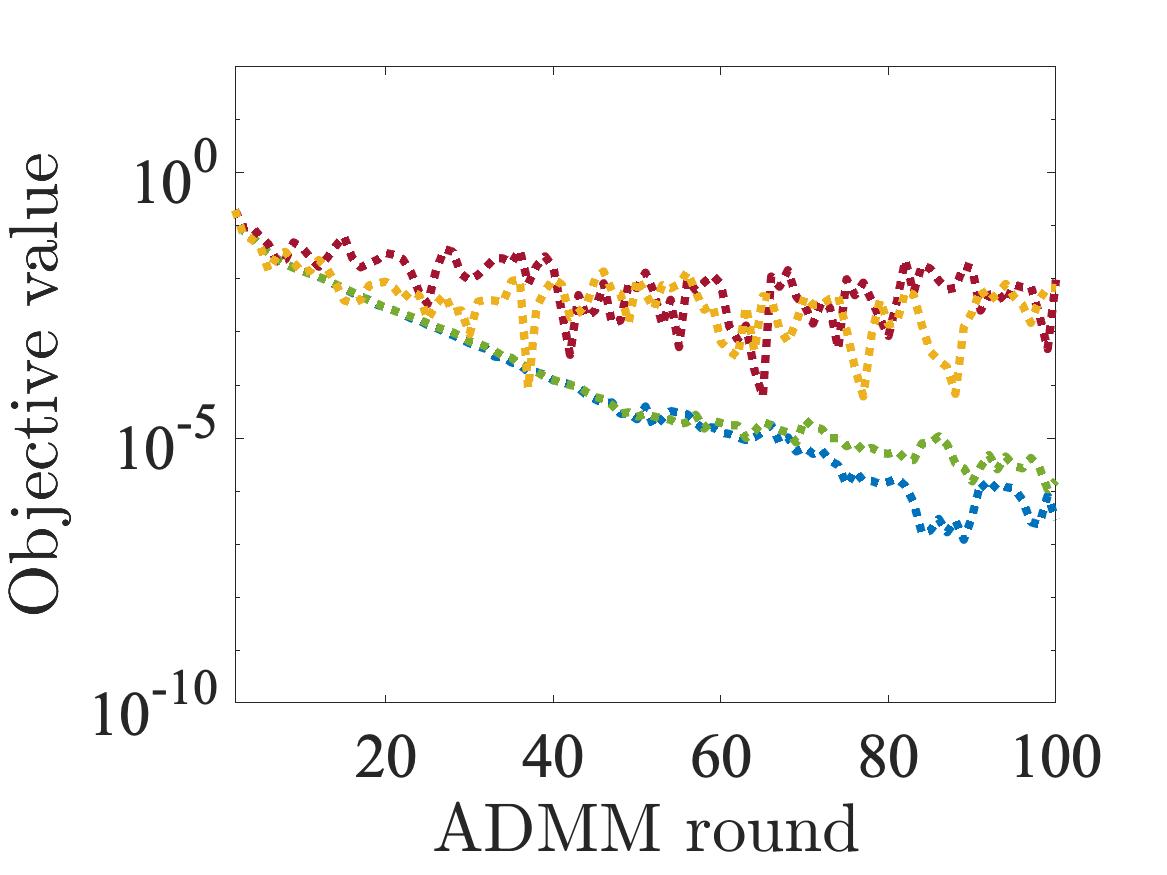}   
\includegraphics[width=\textwidth]{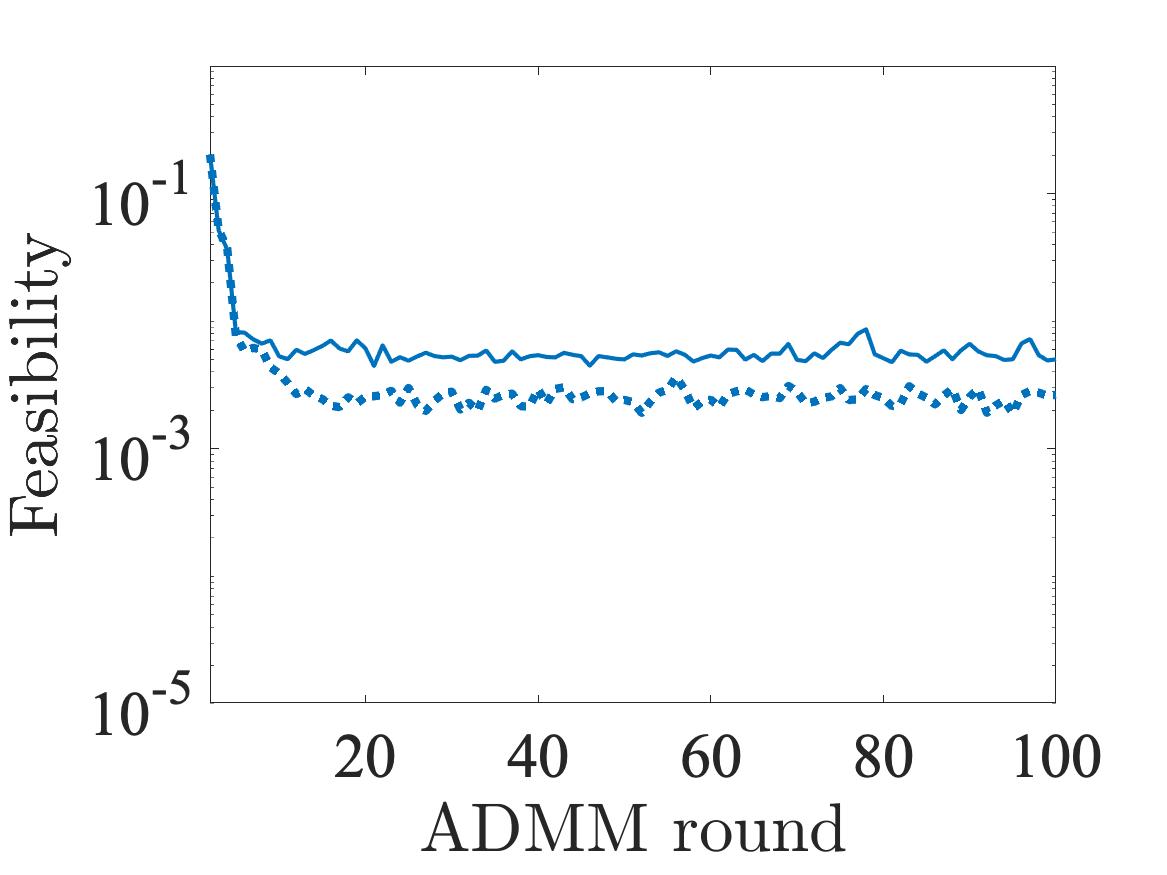}  
\caption{$\bar{\epsilon}=0.1$}
\end{subfigure}
\begin{subfigure}[b]{0.23\textwidth}
\centering      
\includegraphics[width=\textwidth]{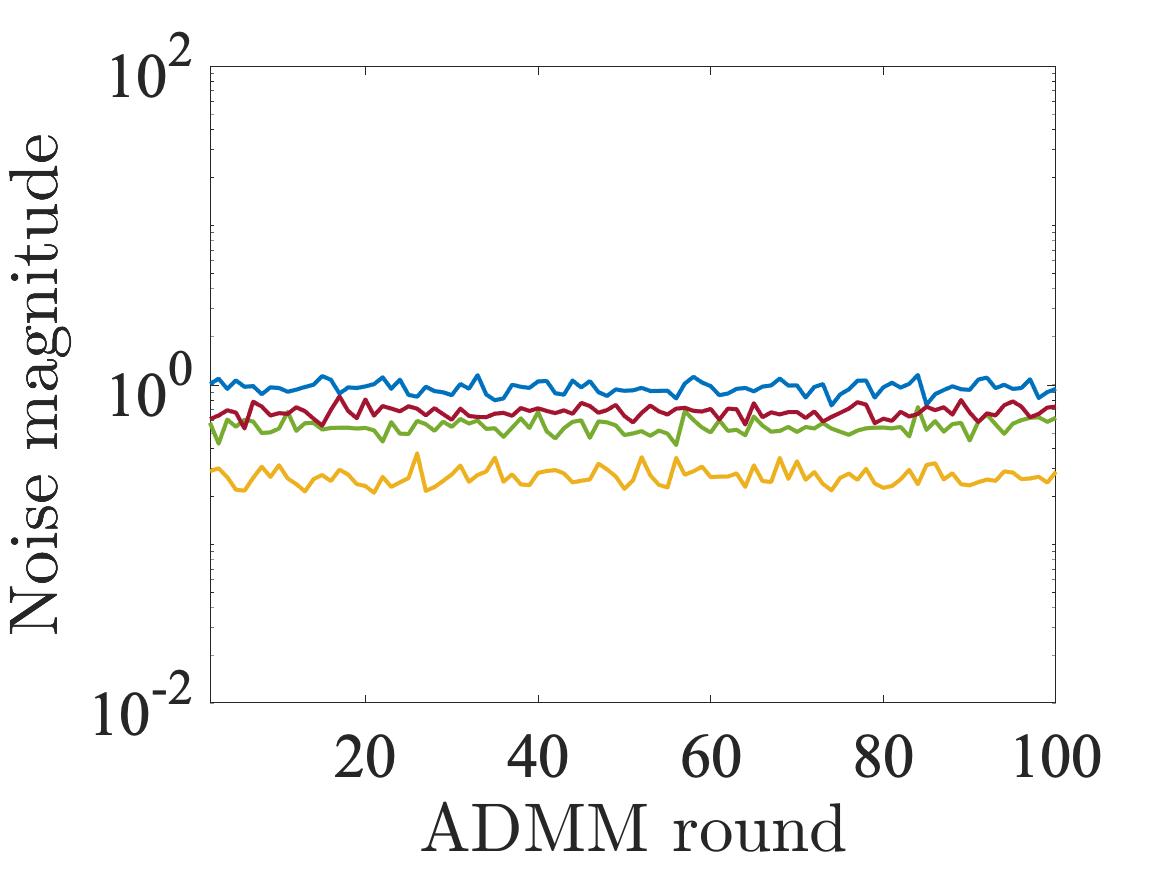}  
\includegraphics[width=\textwidth]{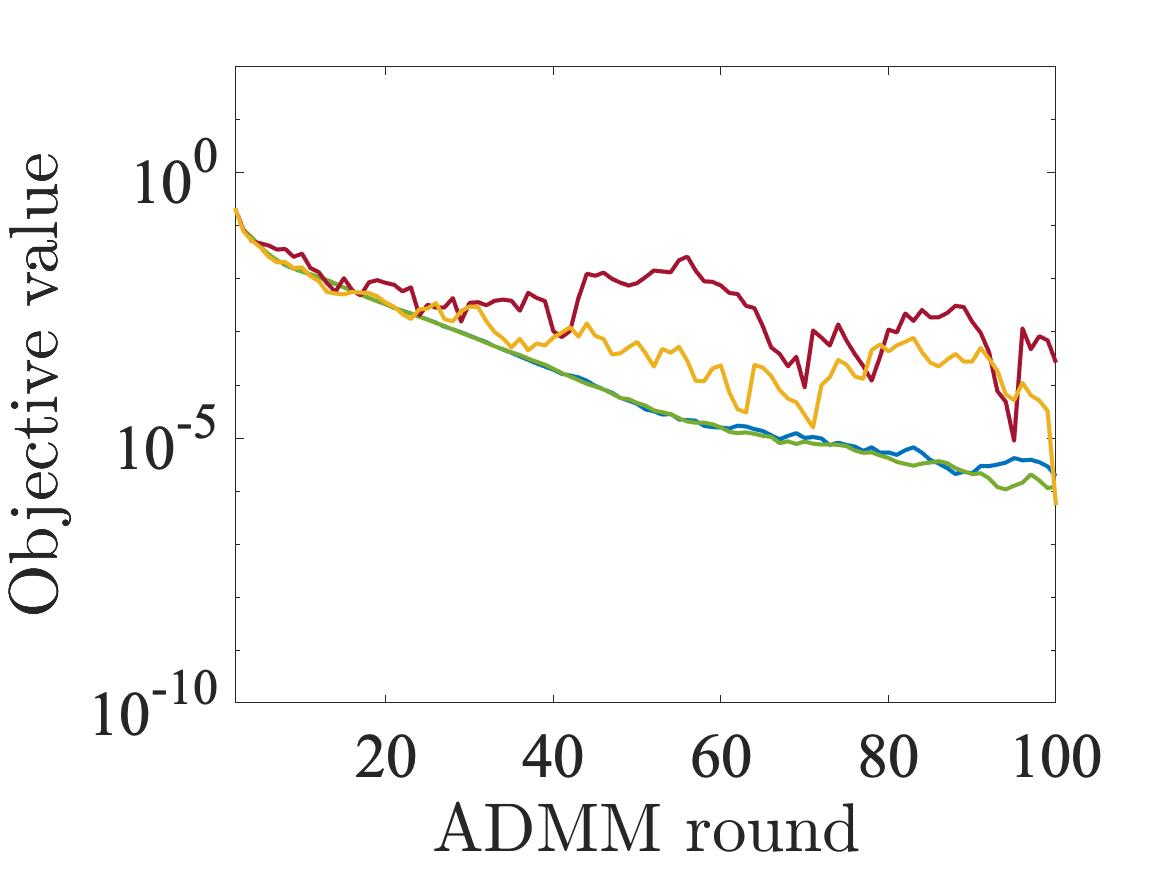}   
\includegraphics[width=\textwidth]{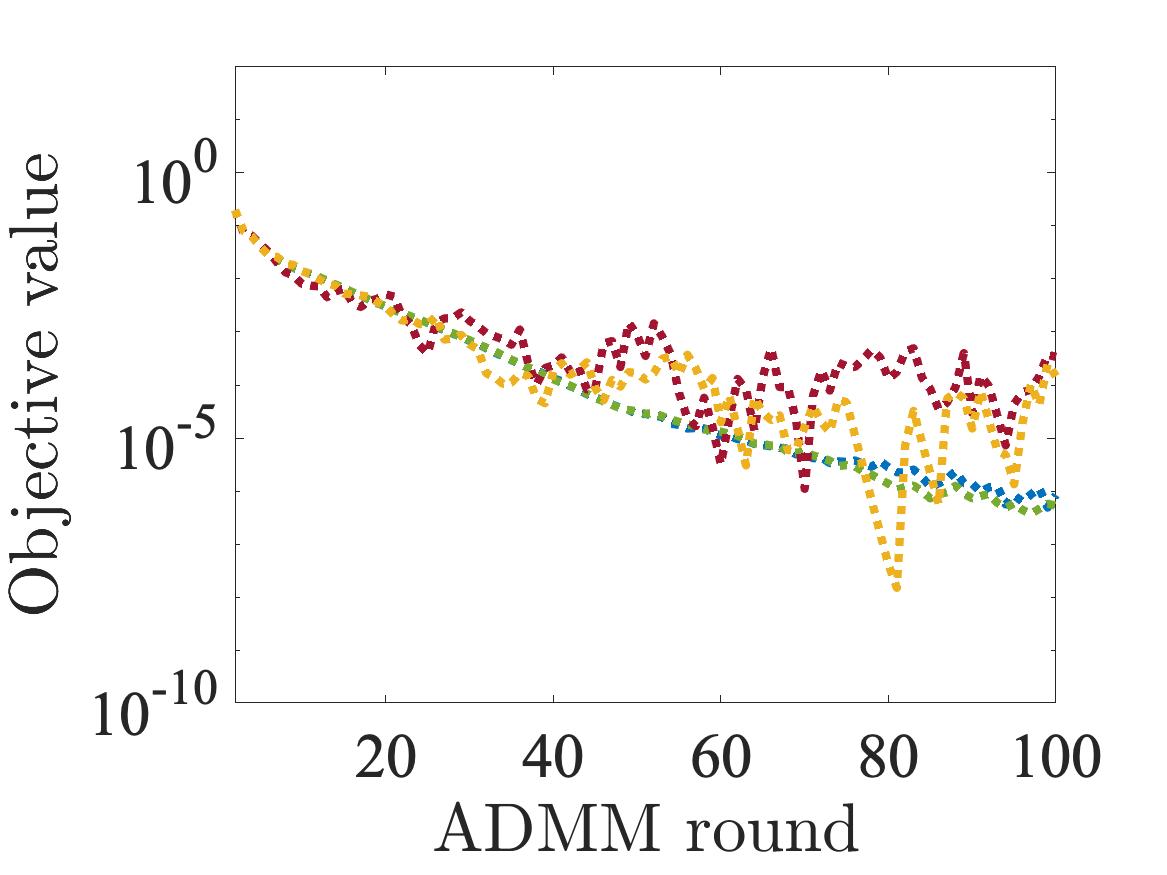}   
\includegraphics[width=\textwidth]{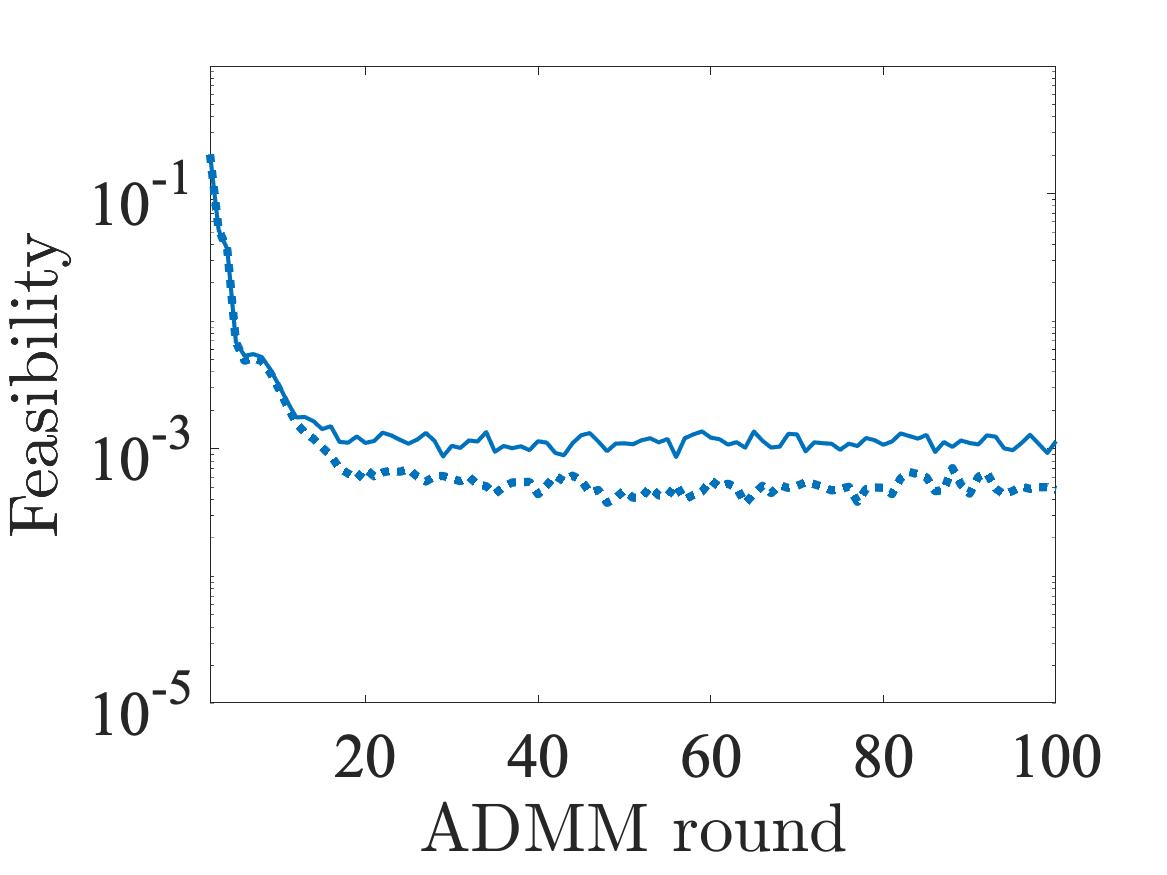}  
\caption{$\bar{\epsilon}=0.5$}
\end{subfigure}     
\begin{subfigure}[b]{0.23\textwidth}
\centering      
\includegraphics[width=\textwidth]{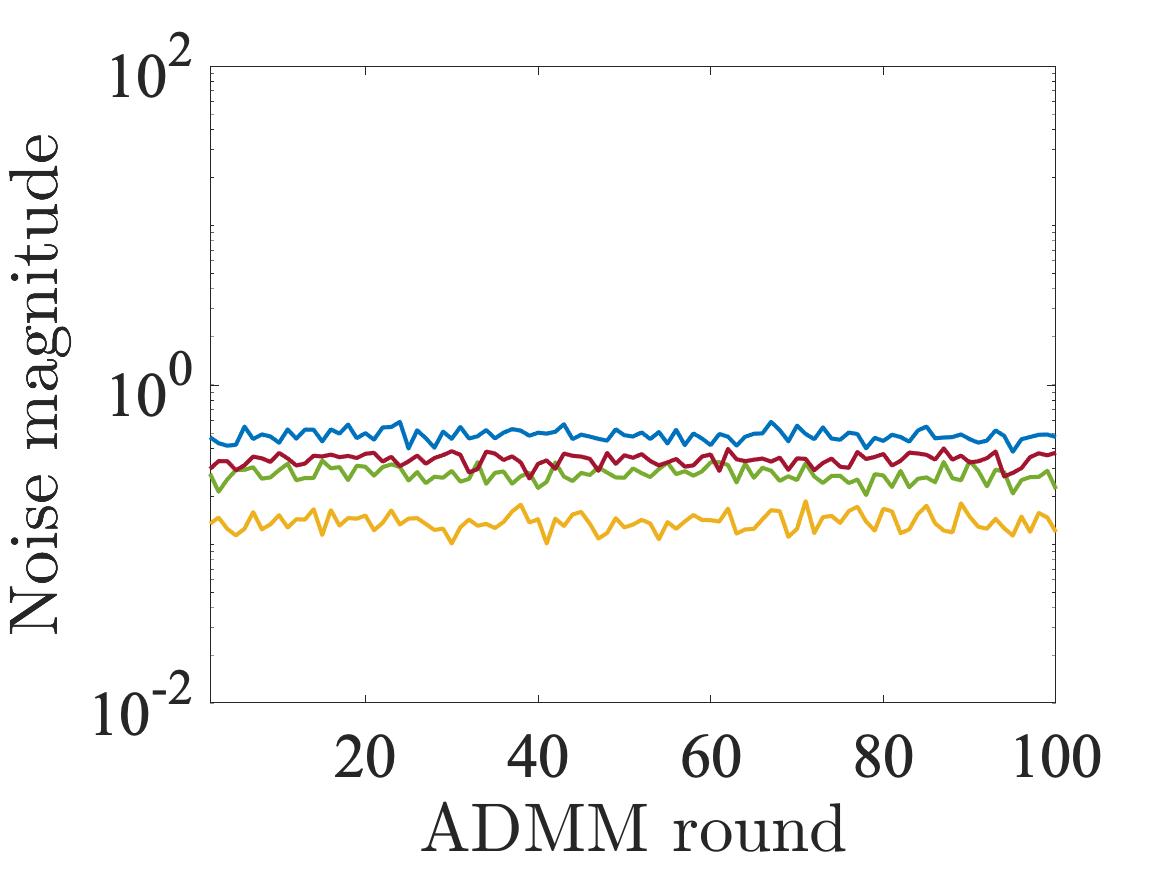}  
\includegraphics[width=\textwidth]{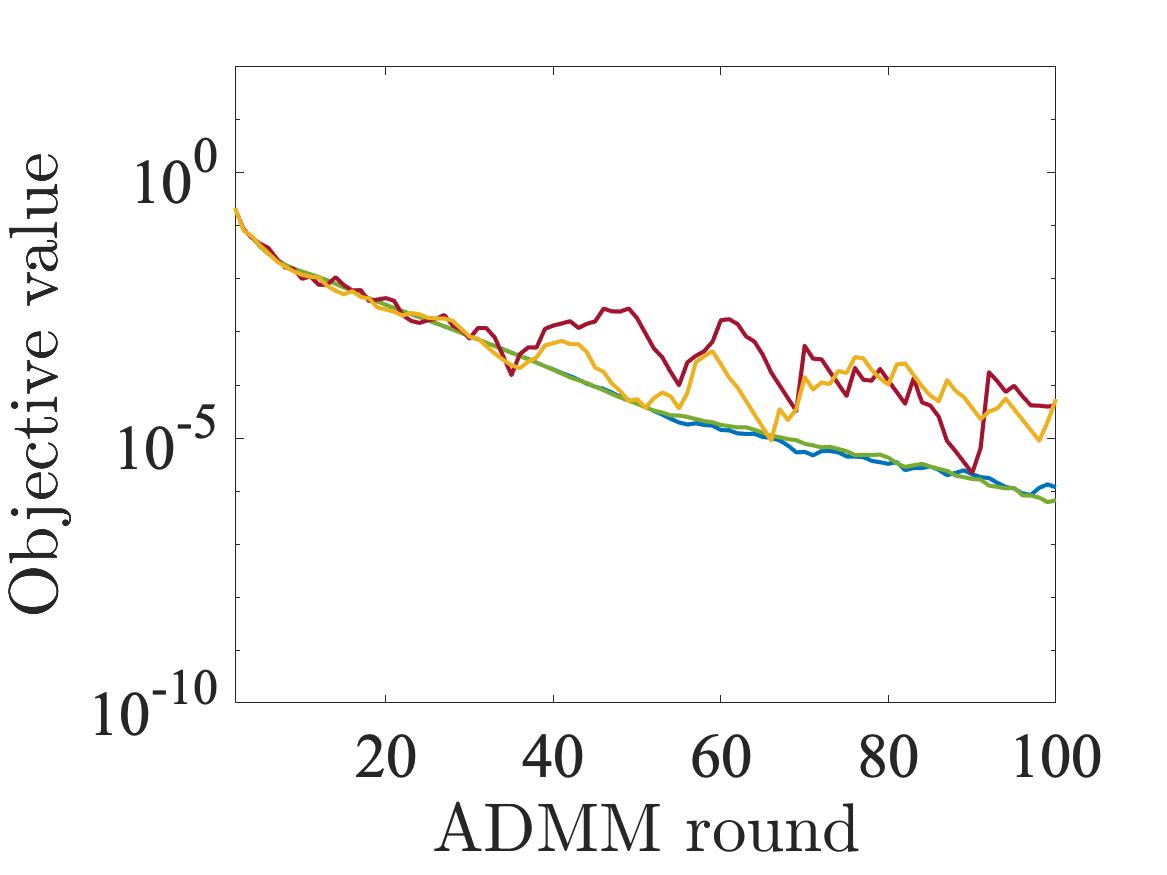}   
\includegraphics[width=\textwidth]{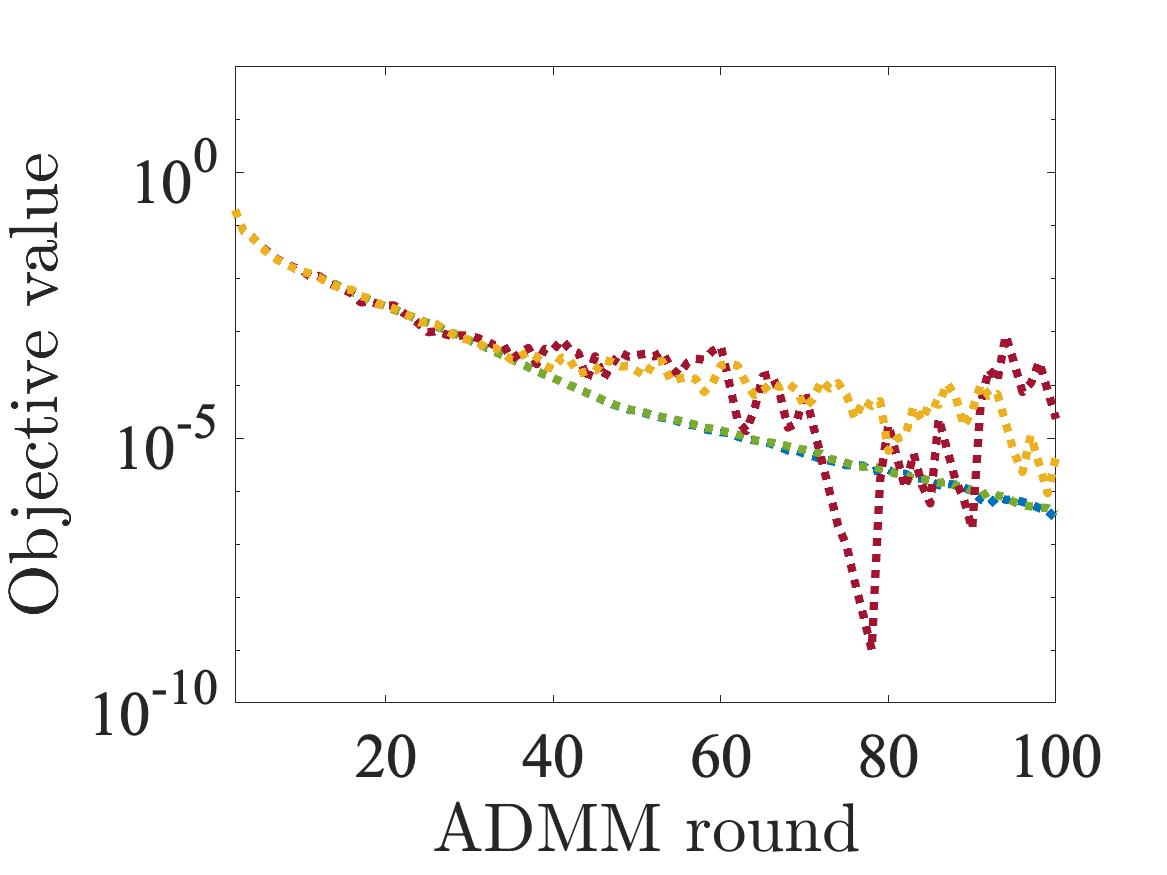}   
\includegraphics[width=\textwidth]{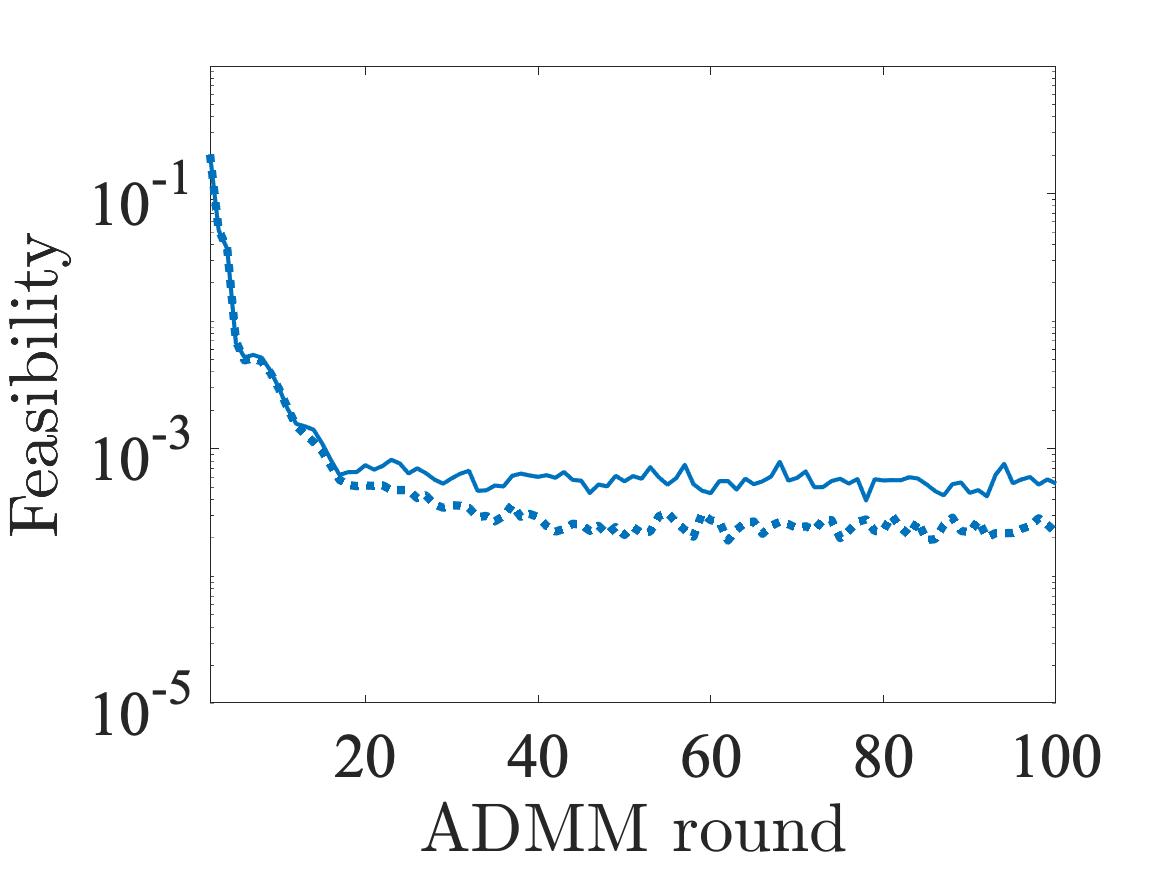}  
\caption{$\bar{\epsilon}=1$}
\end{subfigure}       
\caption{Comparison of \texttt{ObjG}, \texttt{ObjL}, \texttt{OutG}, \texttt{OutL} in terms of noise magnitude ($1^{\text{st}}$ row), objective value ($2^{\text{nd}}$ and $3^{\text{rd}}$ rows), and feasibility ($4^{\text{th}}$ row) by using case 14.}
\label{fig:case14}
\end{figure*}

\begin{figure*}[!ht]  
\centering
\begin{subfigure}[b]{0.23\textwidth}
\centering        
\includegraphics[width=\textwidth]{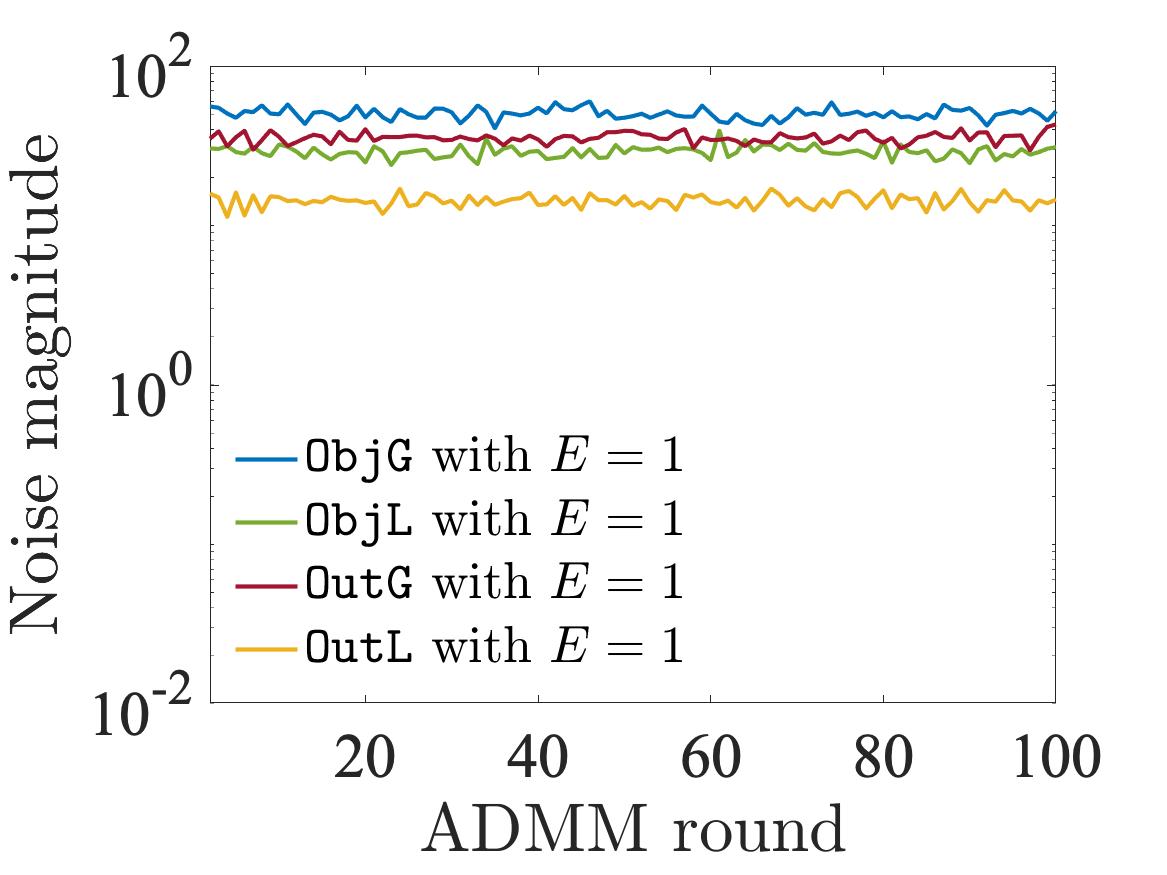}  
\includegraphics[width=\textwidth]{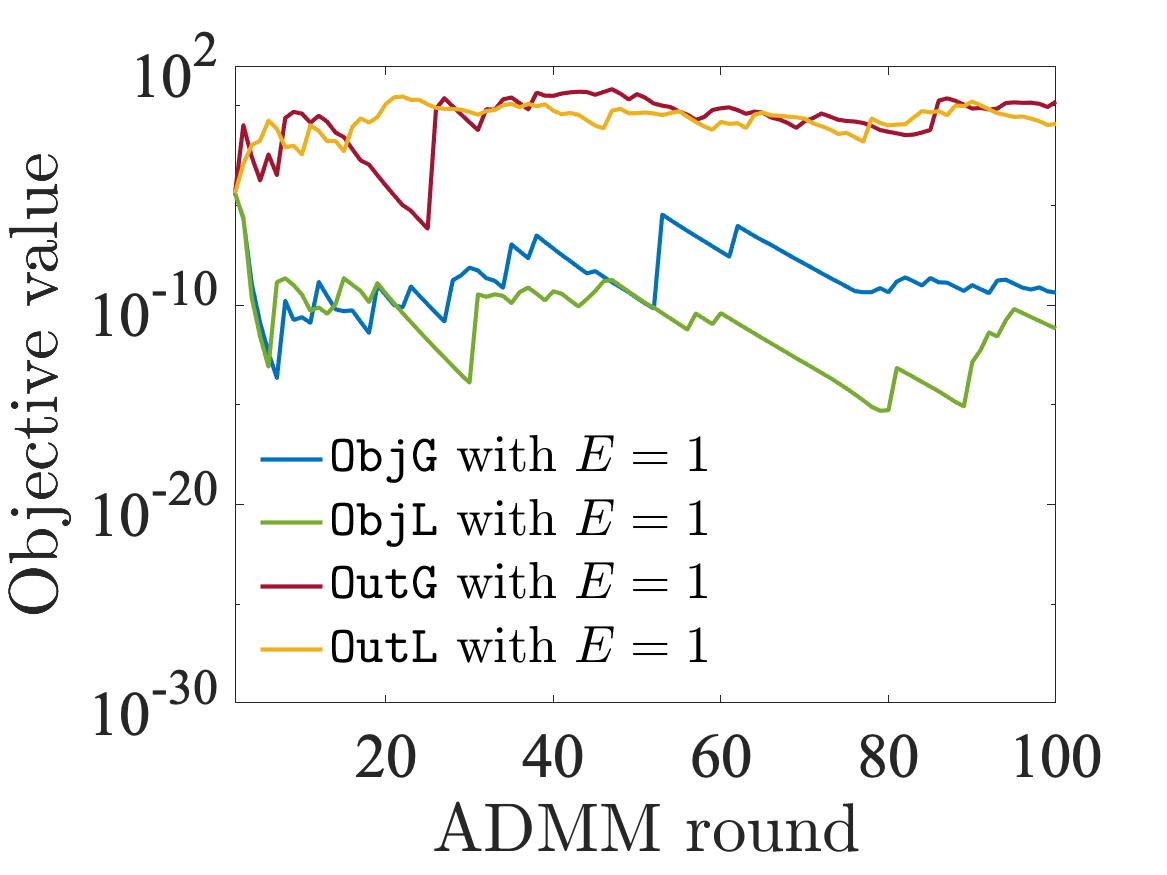}   
\includegraphics[width=\textwidth]{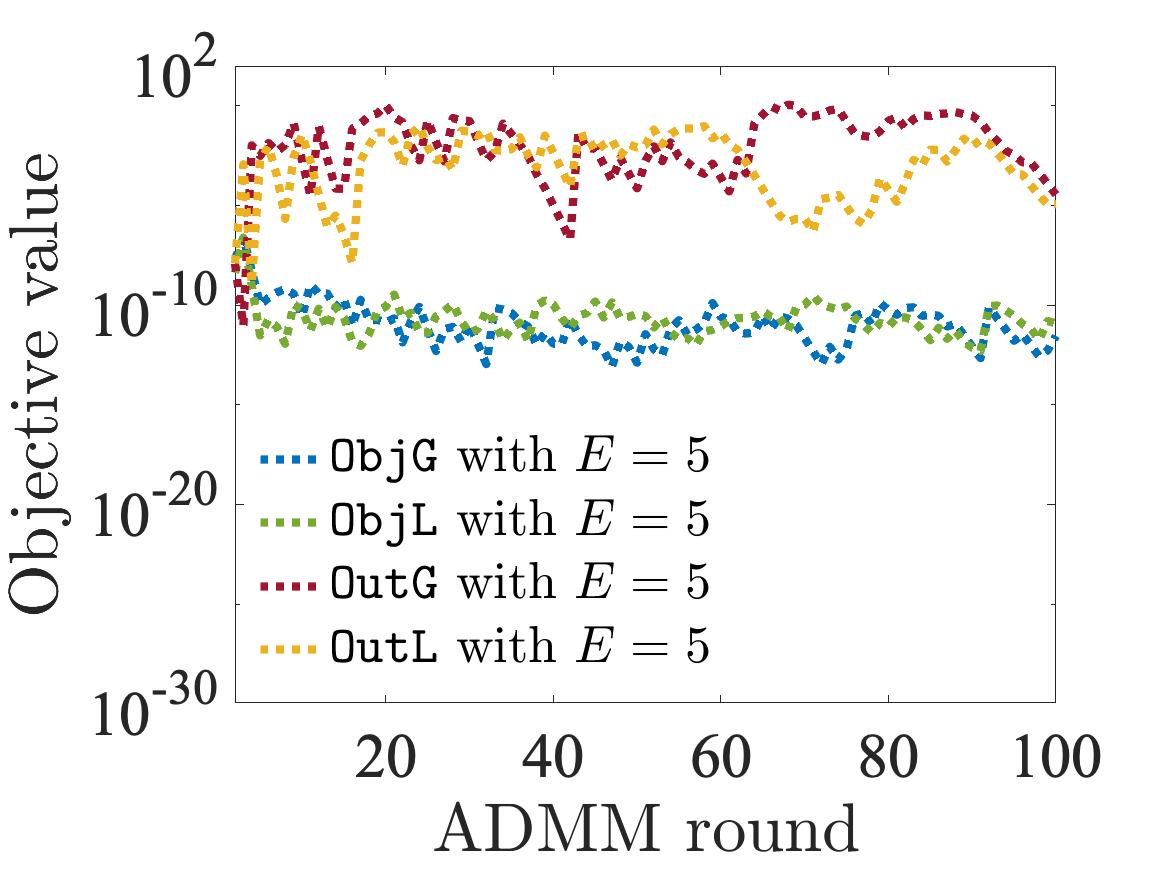}   
\includegraphics[width=\textwidth]{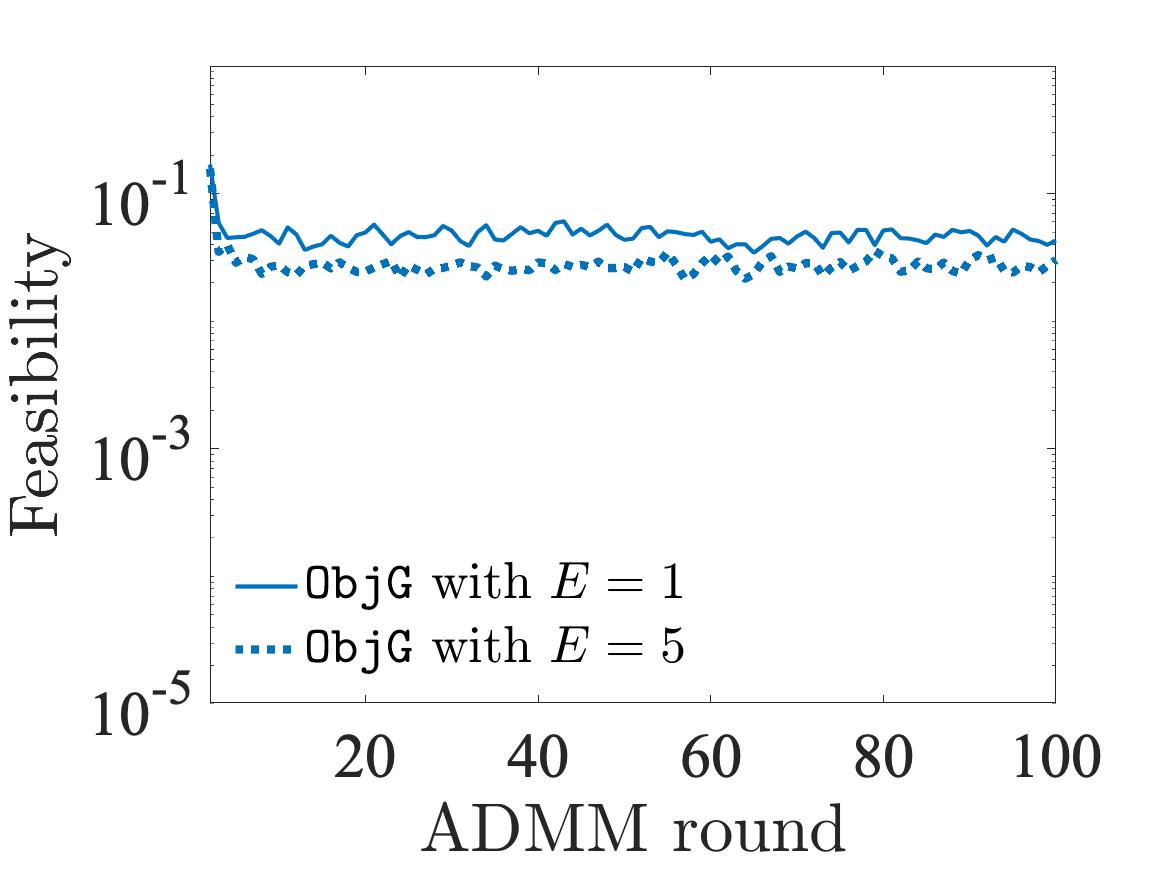}  
\caption{$\bar{\epsilon}=0.05$}
\end{subfigure}
\begin{subfigure}[b]{0.23\textwidth}
\centering        
\includegraphics[width=\textwidth]{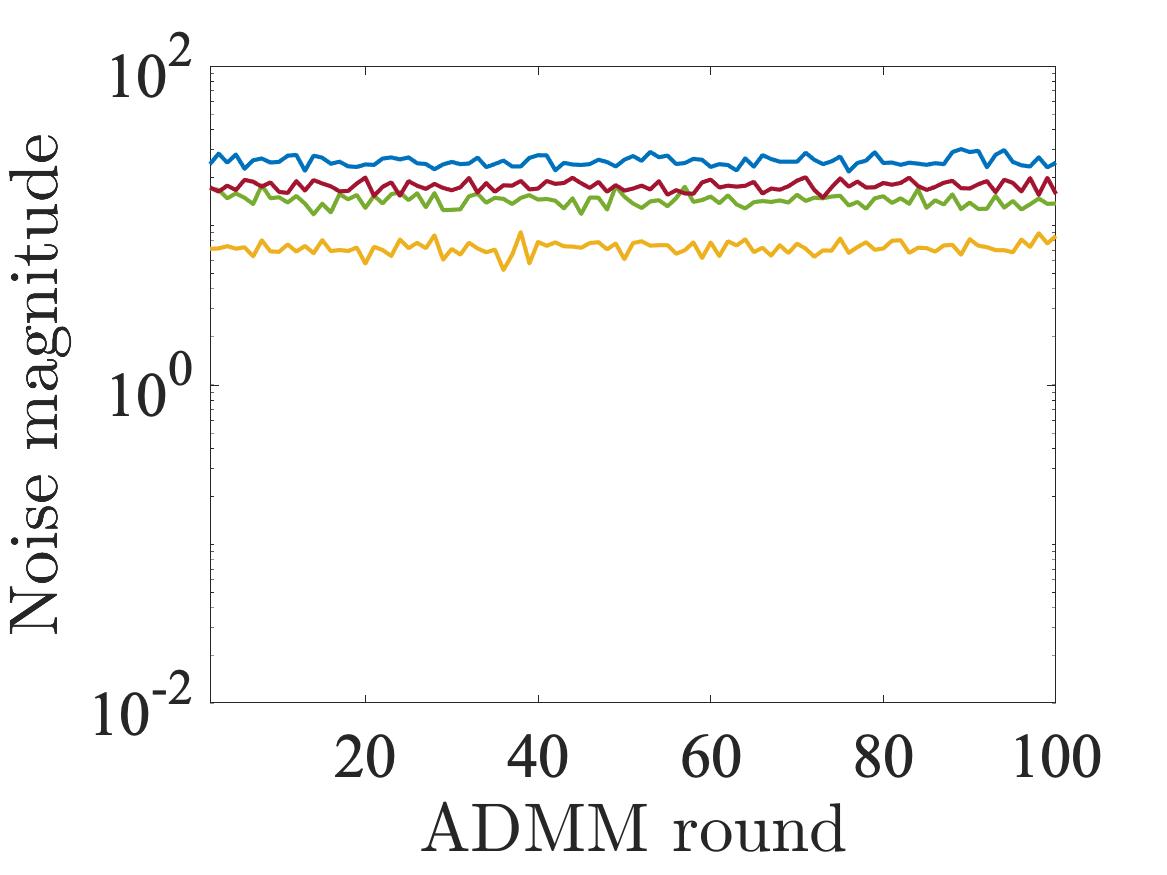}  
\includegraphics[width=\textwidth]{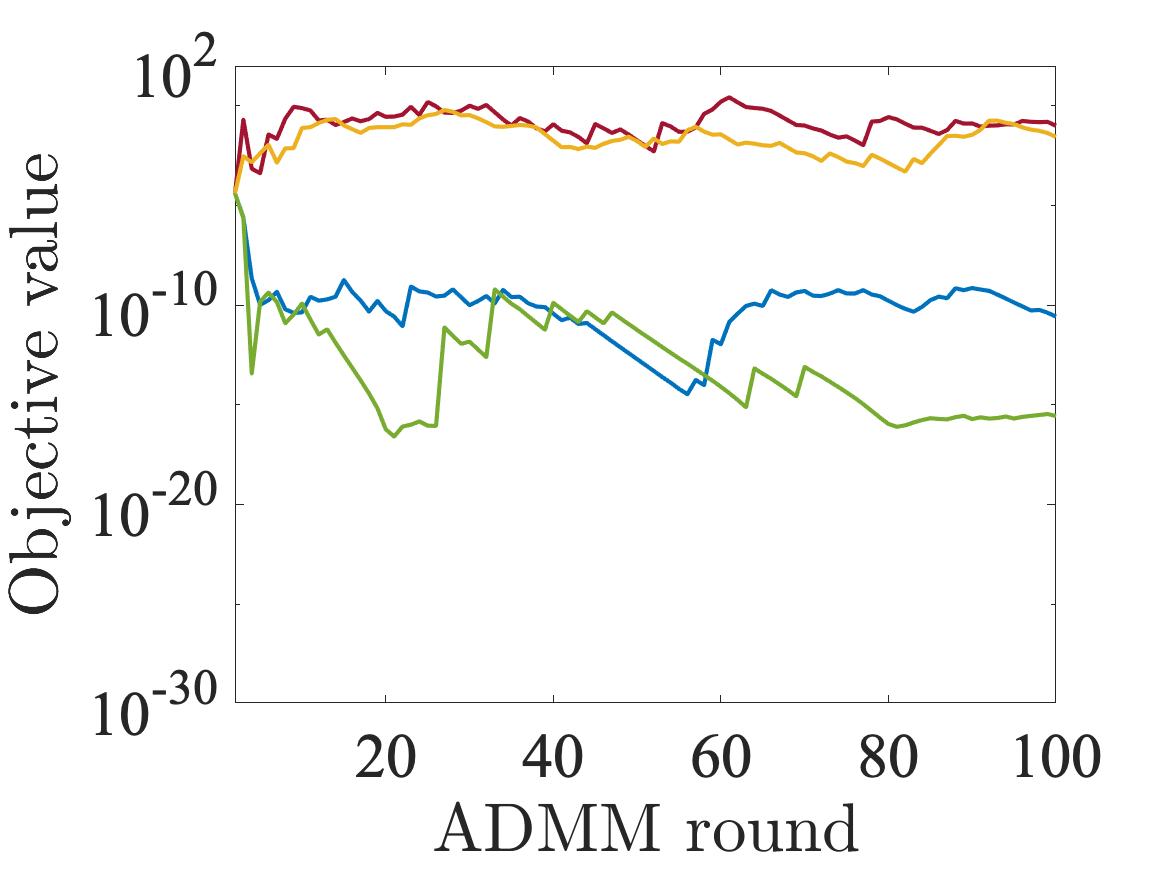}   
\includegraphics[width=\textwidth]{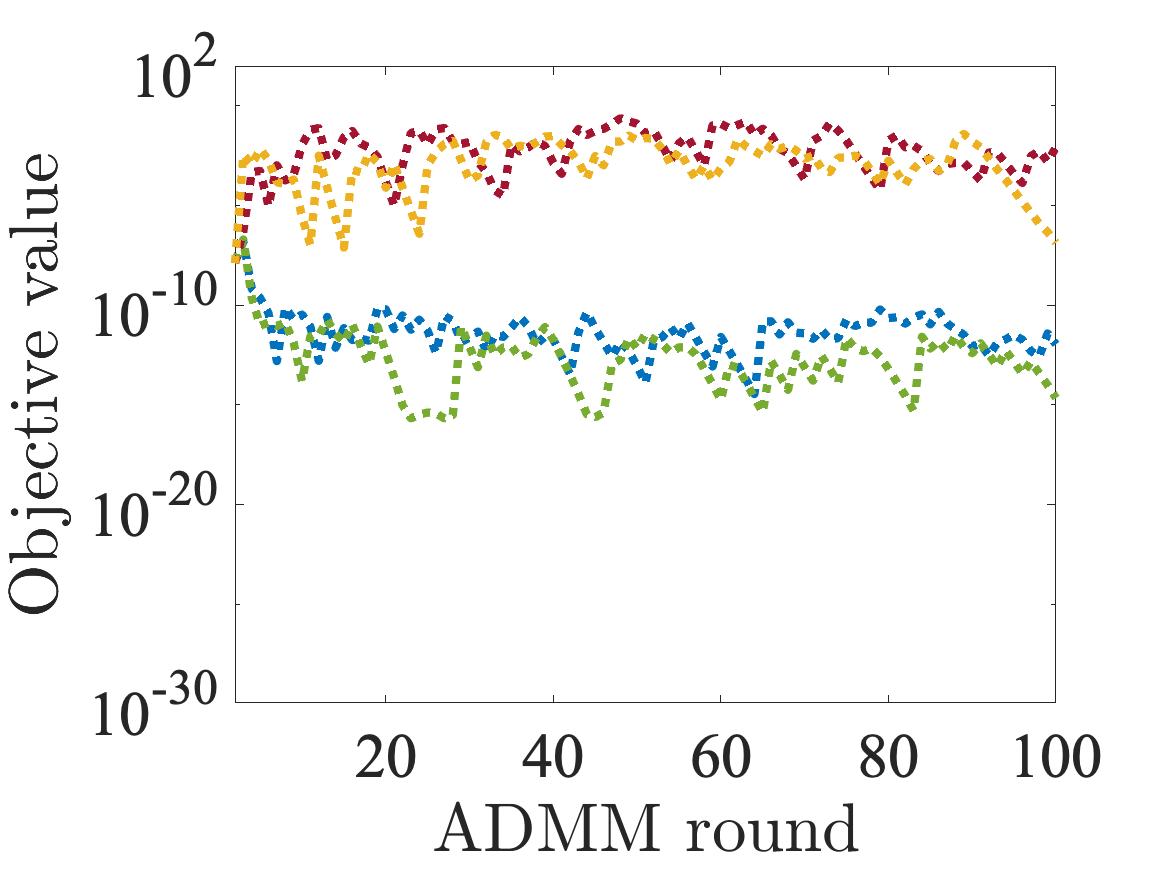}   
\includegraphics[width=\textwidth]{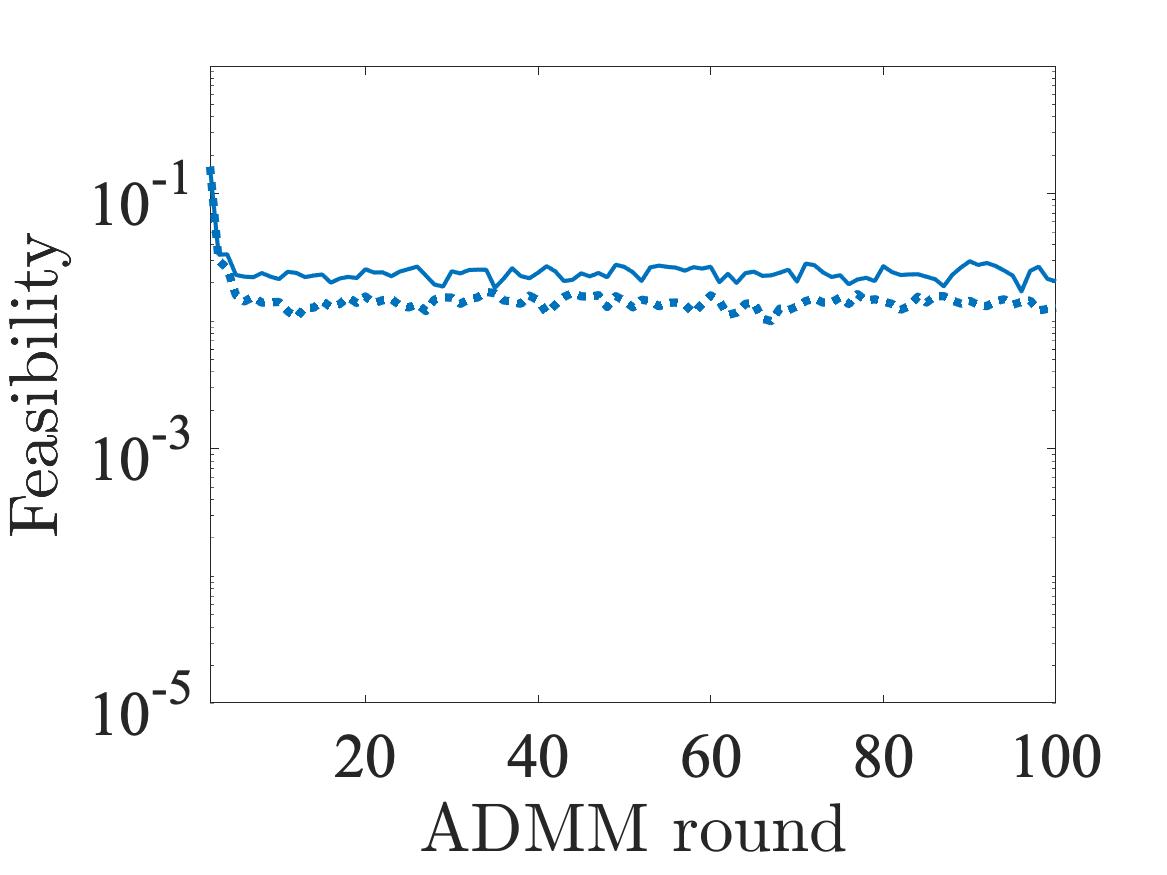}  
\caption{$\bar{\epsilon}=0.1$}
\end{subfigure}
\begin{subfigure}[b]{0.23\textwidth}
\centering      
\includegraphics[width=\textwidth]{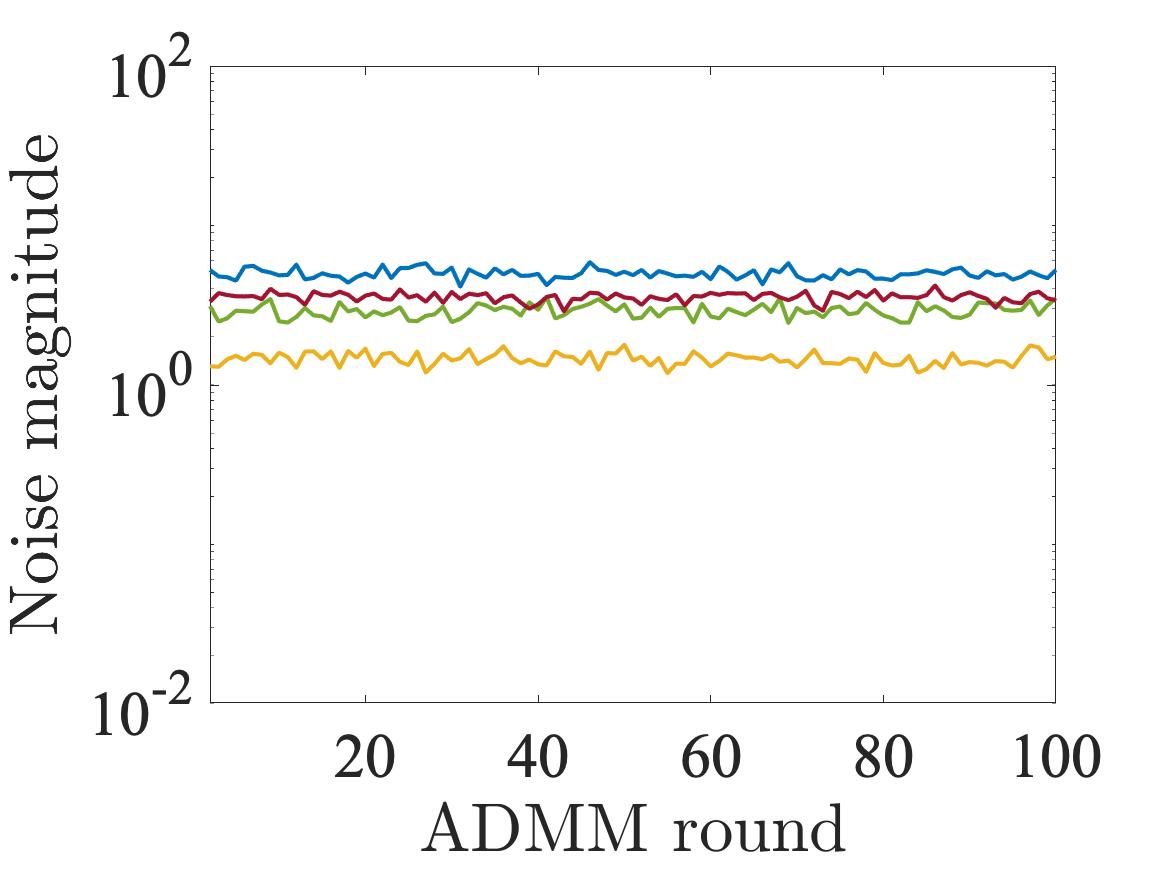}  
\includegraphics[width=\textwidth]{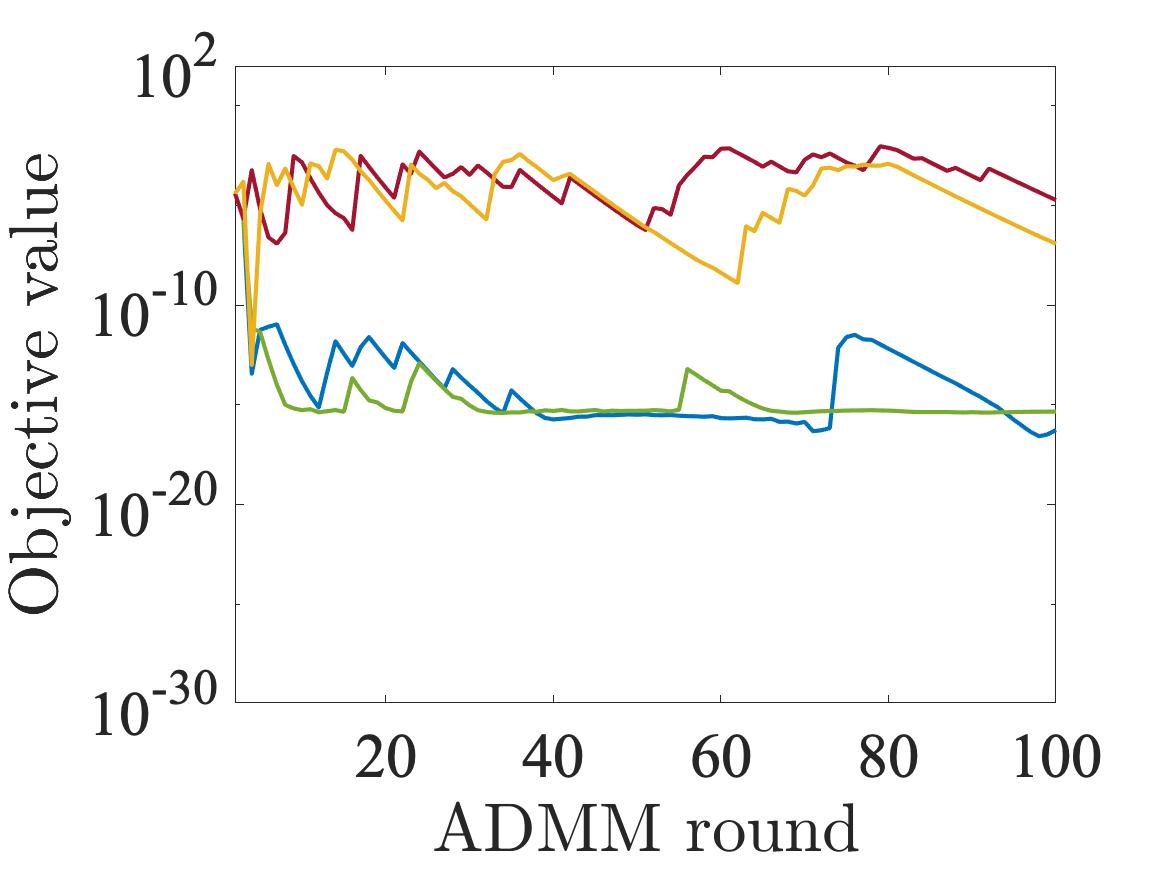}   
\includegraphics[width=\textwidth]{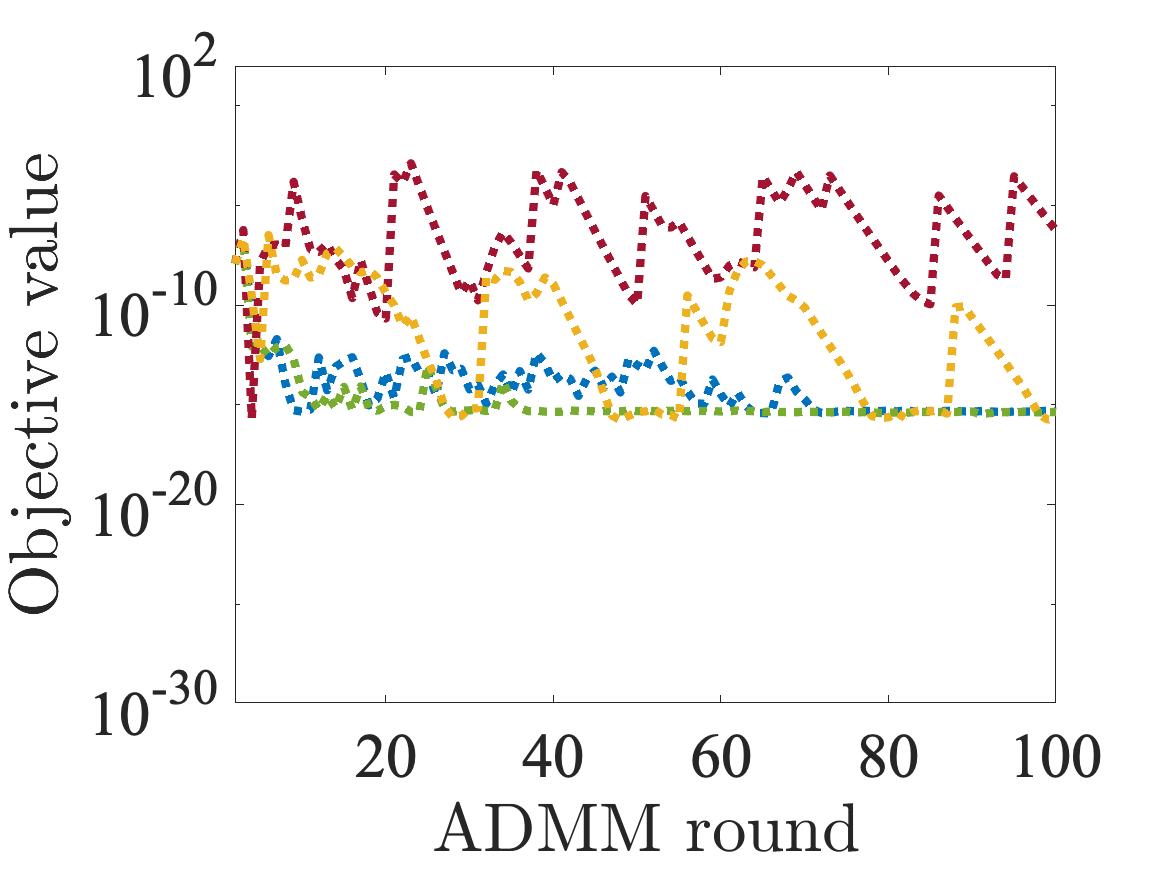}   
\includegraphics[width=\textwidth]{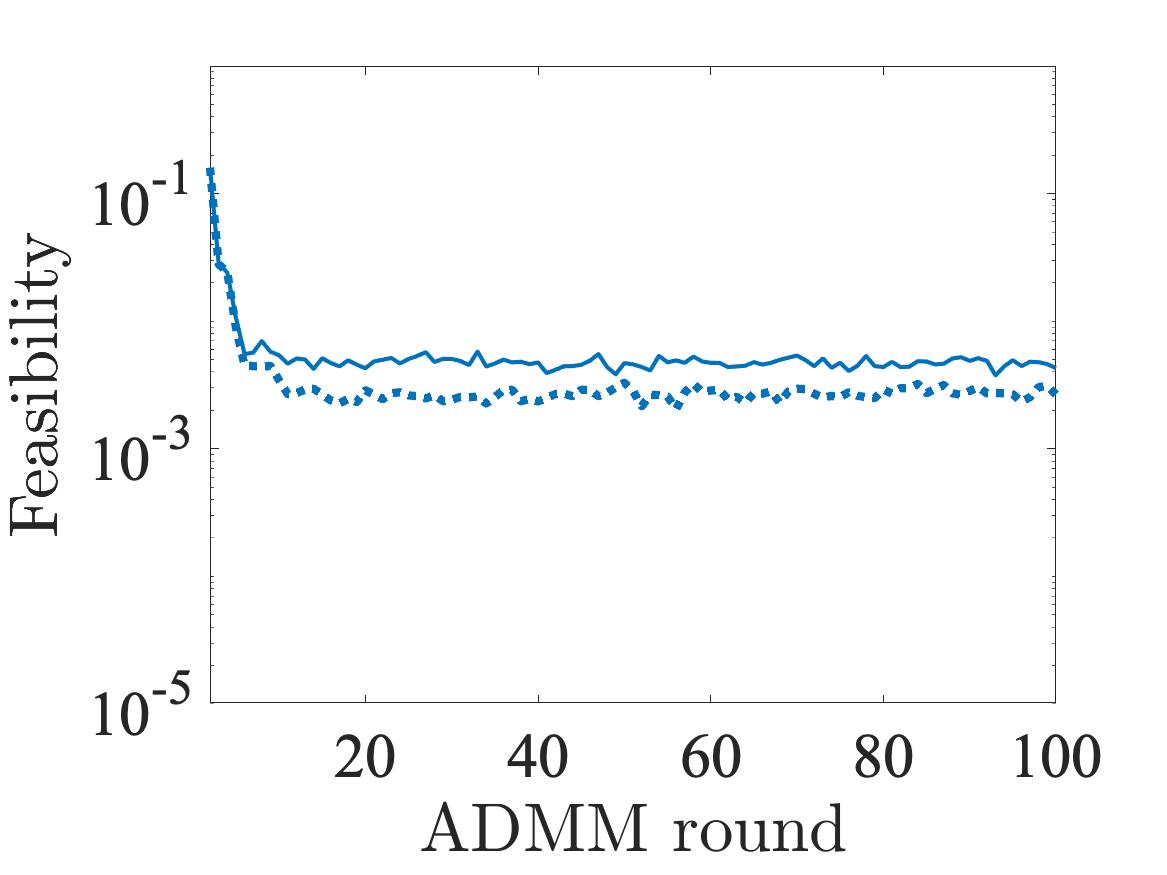}  
\caption{$\bar{\epsilon}=0.5$}
\end{subfigure}     
\begin{subfigure}[b]{0.23\textwidth}
\centering      
\includegraphics[width=\textwidth]{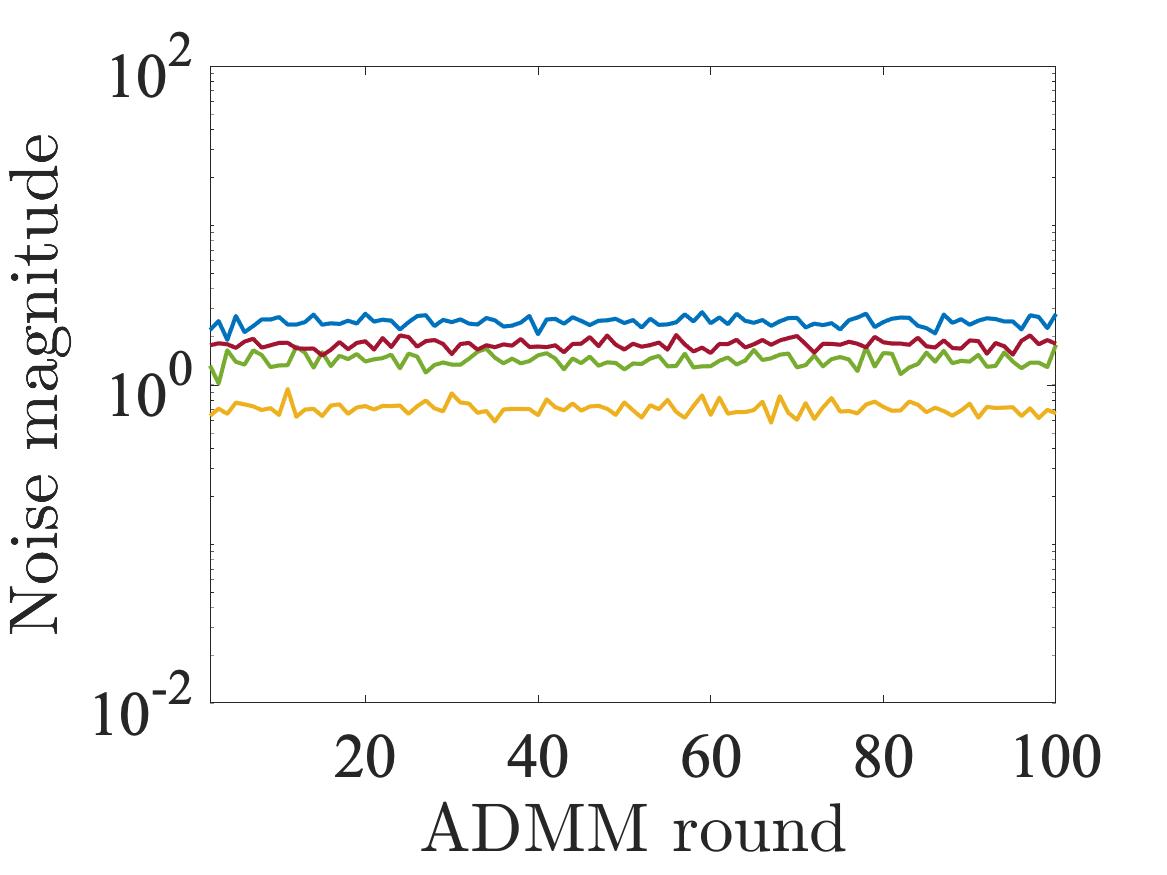}  
\includegraphics[width=\textwidth]{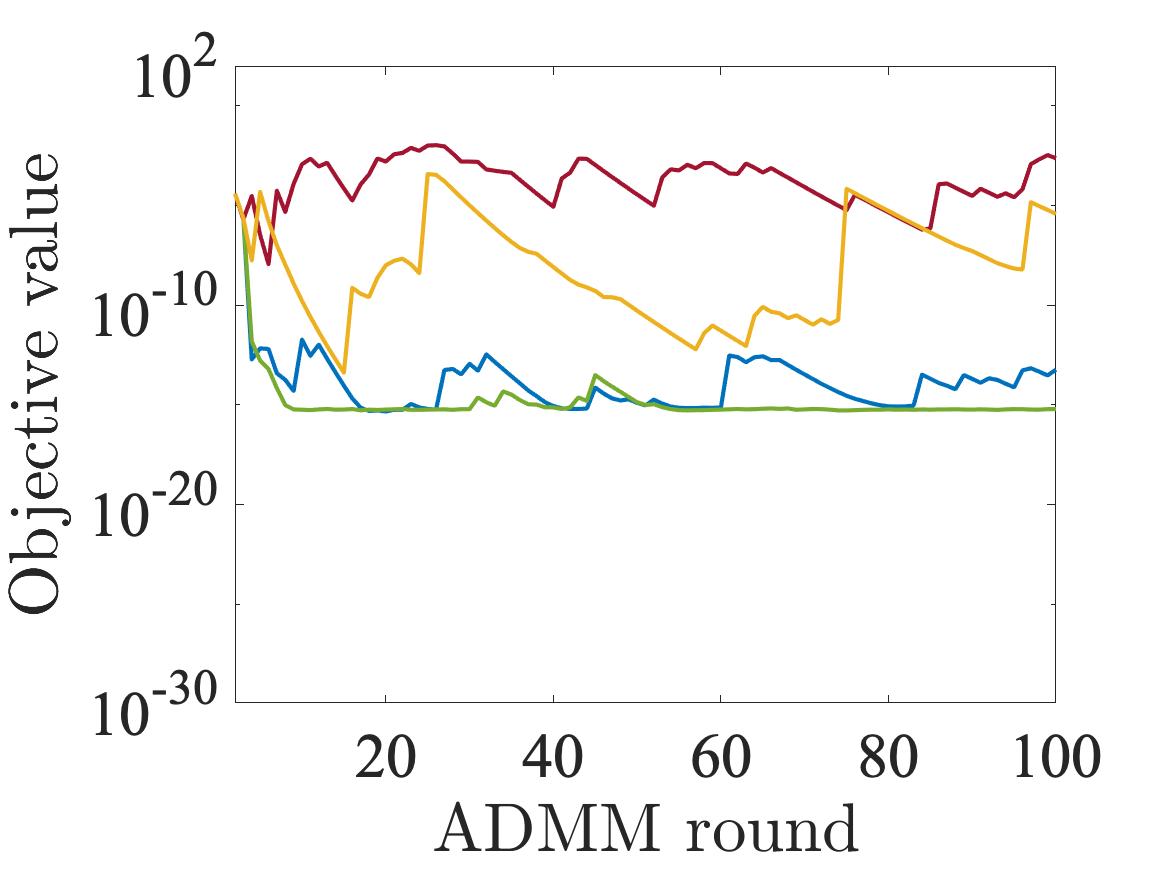}   
\includegraphics[width=\textwidth]{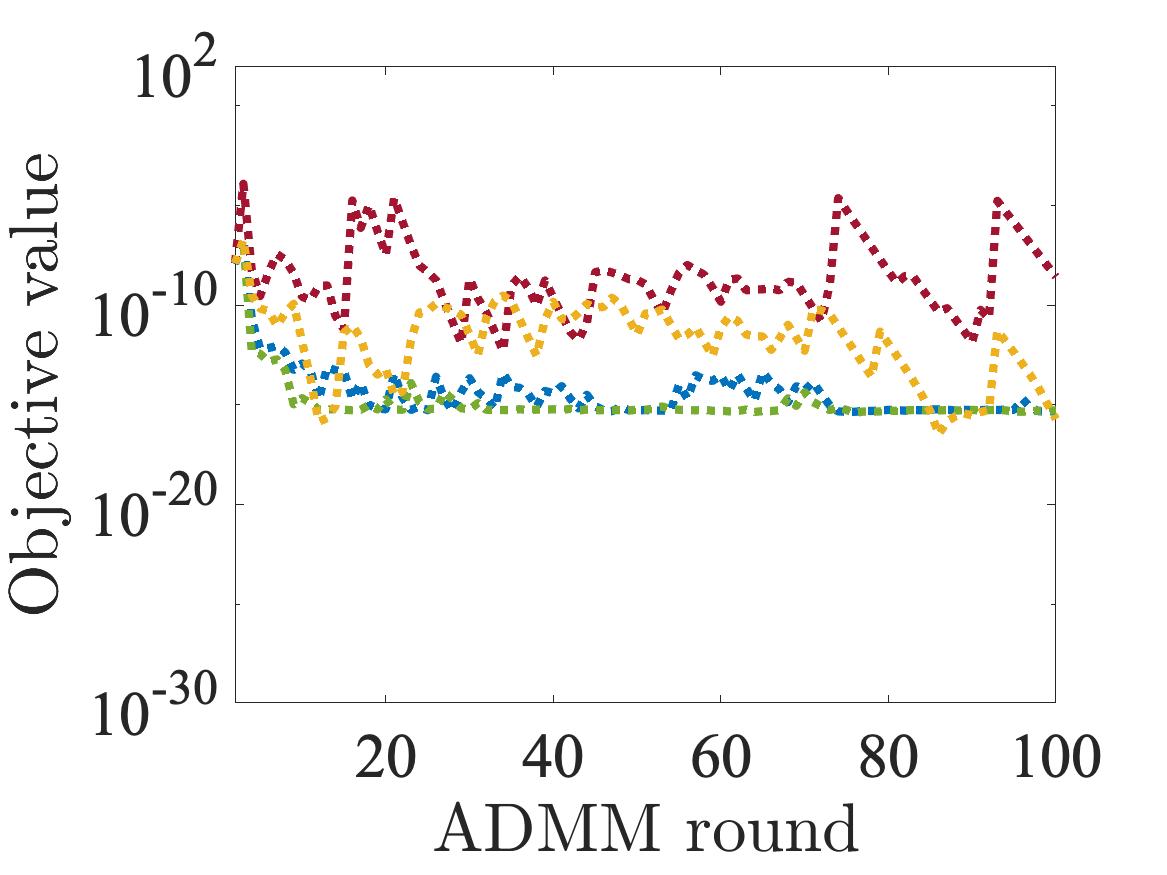}   
\includegraphics[width=\textwidth]{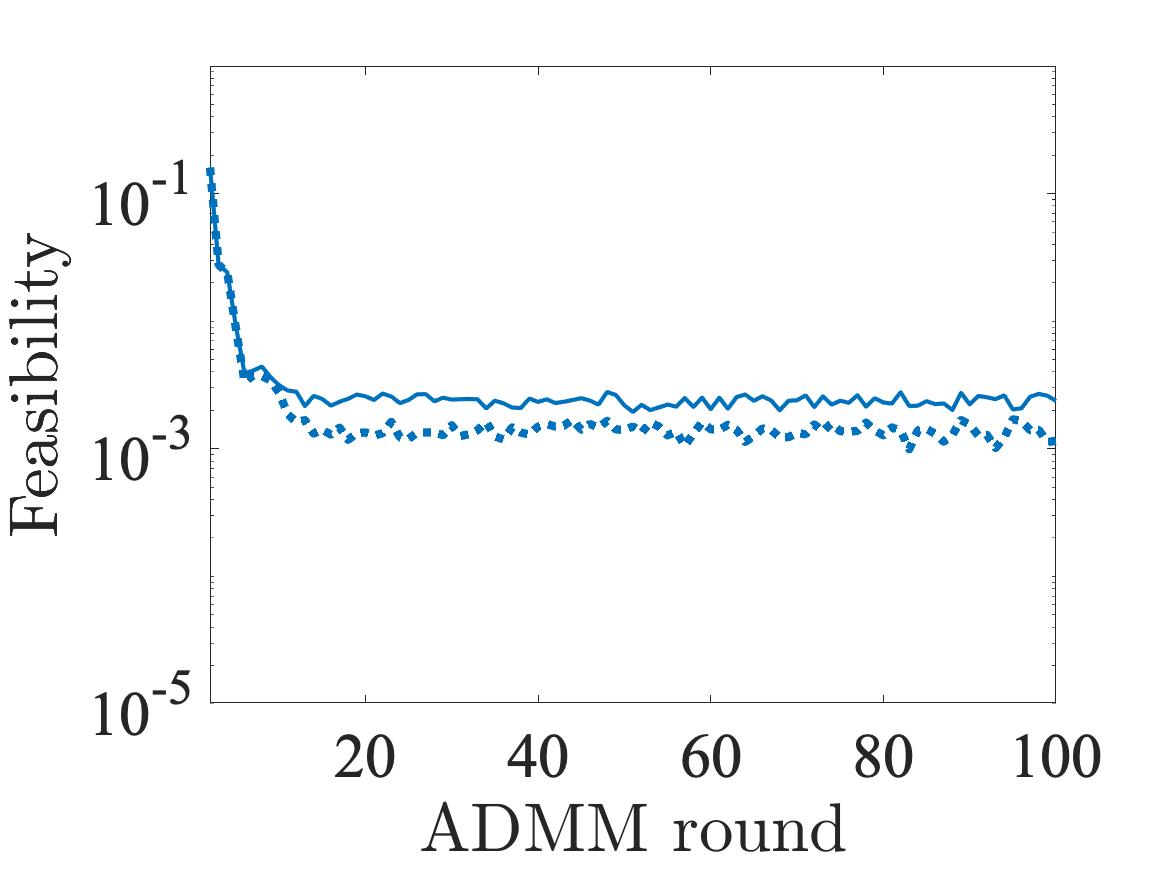}  
\caption{$\bar{\epsilon}=1$}
\end{subfigure}       
\caption{Comparison of \texttt{ObjG}, \texttt{ObjL}, \texttt{OutG}, \texttt{OutL} in terms of noise magnitude ($1^{\text{st}}$ row), objective value ($2^{\text{nd}}$ and $3^{\text{rd}}$ rows), and feasibility ($4^{\text{th}}$ row) by using case 118.}
\label{fig:case118}
\end{figure*}  

\subsubsection{Performance comparison}
We report the performance of \texttt{ObjG}, \texttt{ObjL}, \texttt{OutG}, and \texttt{OutL} for case 14 and case 118 instances in Figure \ref{fig:case14} and \ref{fig:case118}, respectively. 

First, we compare the four algorithms with respect to the magnitude of noise introduced during $T=100$ ADMM rounds, which is defined as
\begin{align}
  \sum_{p=1}^P  \sum_{e=1}^E  \sum_{i=1}^n | \tilde{\xi}^{t,e}_{p i} |, \ \forall t \in [T], \label{form:noise_mag}
\end{align}
where $\tilde{\xi}^{t,e}_{pi} \in \mathbb{R}$ is sampled from either Laplace or Gaussian distribution.
The noise magnitude \eqref{form:noise_mag} under various $\bar{\epsilon}$ is reported in figures from the first row of Figure \ref{fig:case14} and \ref{fig:case118}.
From these figures, we make the following observations:
\begin{enumerate}
  \item The noise magnitude increases as $\bar{\epsilon}$ decreases.
  \item For fixed $\bar{\epsilon}$, the proposed algorithm \texttt{ObjG} (resp., \texttt{ObjL}) requires more noise compared with \texttt{OutG} (resp., \texttt{OutL}).
  \item For fixed $\bar{\epsilon}$, algorithms with the Gaussian mechanism (i.e., \texttt{ObjG} and \texttt{OutG}) require more noise compared with algorithms with the Laplace mechanism (i.e., \texttt{ObjL} and \texttt{OutL}).
\end{enumerate} 
To explain these observations, we report in Table \ref{table:var} the variance of distributions used for sampling each element of $\tilde{\xi}^{t,e}_p$.
For all algorithms considered in the table, decreasing $\bar{\epsilon}$ increases the variance, which can explain the first observation.
From the table we note that \texttt{ObjG} has the largest variance, followed by \texttt{OutG}, \texttt{ObjL}, and \texttt{OutL}. 
These explain the second and third observations.

\begin{table}[H]
\centering
\caption{Variance of distributions used for sampling noise in the four algorithms.}
\begin{tabular}{ c c c c}
\hline
\texttt{ObjG} & \texttt{ObjL}  & \texttt{OutG}  & \texttt{OutL} \\ \hline    
$2 \ln ( \frac{1.25}{\bar{\delta}}) (\frac{2\beta}{\bar{\epsilon}})^2$ 
& $2  (\frac{2\beta}{\bar{\epsilon}})^2$ 
& $2 \ln ( \frac{1.25}{\bar{\delta}}) (\frac{\beta}{\bar{\epsilon}})^2$ 
& $2  (\frac{\beta}{\bar{\epsilon}})^2$   \\
\hline
\end{tabular}
\label{table:var}
\end{table}

Second, we report the objective function values (which should be zero at optimality in a non-private setting) produced by the four algorithms under various $\bar{\epsilon}$ in the figures located at the second row of Figures \ref{fig:case14} and \ref{fig:case118}, when $E=1$, namely, single local update.
We notice that the proposed algorithms \texttt{ObjG} and \texttt{ObjL} produce lower objective function values compared with \texttt{OutG} and \texttt{OutL} even though our algorithm requires more noise to be introduced for guaranteeing the same level of data privacy. 
This shows the effectiveness of our approach, which ensures feasibility of the intermediate solutions randomized for DP.
We also note that \texttt{ObjL} performs better than \texttt{ObjG}. The reason is  that \texttt{ObjL} requires less noise in this experiment.

Third, we set the number of local updates to $E=5$ and report the objective values in the figures located at the third row of Figures \ref{fig:case14} and \ref{fig:case118}. For larger $\bar{\epsilon}$, \texttt{OutG} and \texttt{OutL} sometimes produce better objective function values. Clearly, however,  our algorithms provide better objective values for smaller $\bar{\epsilon}$ (i.e., stronger data privacy).
The reason is mainly that there is no $\bar{\epsilon}$-dependent error bound in our convergence rate presented in Section \ref{sec:convergence} while the error bound exists for the existing DP algorithms \cite{huang2020differentially} based on the output perturbation, which increases as $\bar{\epsilon}$ decreases. 

Fourth, we report feasibility, namely, $\|Aw^t+Bz^t\|_2$, where $w$ is a global solution, $z$ is a local solution, and $Aw+Bz$ represents the consensus constraints (which should be zero at optimality in a non-private setting), produced by \texttt{ObjG} for $E \in \{1,5\}$ in figures located at the last row of Figures \ref{fig:case14} and \ref{fig:case118}. 
We observe that increasing the number of local updates provides a solution with better quality. 
Similar behavior is observed in other algorithms \texttt{ObjL}, \texttt{OutG}, and \texttt{OutL},  numerically demonstrating the effectiveness of introducing the multiple local update.

\subsection{Federated learning} \label{sec:numerical_fl}
In ML, transferring data into a central server for training can be limited because of, for example, low bandwidth or data privacy issues. Such limitations have motivated the use of distributed optimization, also known as federated learning. 
In most FL literature, a distributed empirical risk minimization model is utilized that is a form of \eqref{model:dist}.
Specifically, there are multiple $P$ agents that have their own data and machine for local training. 
These agents cooperate to find a global model parameter $w^*$ by iteratively training using the local objective function $f_p (w) := (\frac{I_p}{I}) \frac{1}{I_p}\sum_{i \in \mathcal{I}_p}  l (w; x_i, y_i)$, where $\mathcal{I}_p$ is an index set of local data samples, $I_p := |\mathcal{I}_p|$ is the number of local data samples, $I:=\sum_{p=1}^P I_p$ is the total number of data samples, 
$l$ represents a loss function, 
and $x_i$ and $y_i$ represent data features and data labels, respectively.
We note that most FL literature assumes $\mathcal{W} := \mathbb{R}^n$, but more recently the importance of imposing hard constraints instead of soft constraints (i.e., penalizing in the objective function) has been discussed in the ML community \cite{gallego2022controlled, marquez2017imposing}. 
 
In particular, we consider a multiclass logistic regression model for classifying image data.
Such problems can be formulated as the form of \eqref{model:dist_1} with the following local objective function:
\begin{align}
f_p (z_p) :=  - \frac{1}{I} \sum_{i=1}^{I_p} \sum_{k=1}^K y_{pik} \ln \Big(\frac{\exp(\sum_{j=1}^J  x_{pij} z_{pjk}) }{\sum_{k'=1}^K \exp(\sum_{j=1}^J  x_{pij} z_{pjk'}) } \Big), \ \forall p \in [P], \label{FL_local_obj}
\end{align}
where $z_p \in \mathbb{R}^{J \times K}$ represents local model parameters and $x_p \in \mathbb{R}^{I_p \times J}$ and $y_p \in \mathbb{R}^{I_p \times K}$ are data features and labels, respectively. Note that for each local data sample $i \in [I_p]$, there are $J$ data features and $K$ data labels.
In this experiment we consider a simple local constraint $\mathcal{W}_p :=  [-u, u]$ for some $u \in \mathbb{R}^{J \times K}_{++}$. 

In this setting agents are not required to share their local image data with other agents. 
 One can, however,  reconstruct the law data if local model parameters are exposed.
To protect data, one can utilize a DP algorithms if one can admit the inevitable trade-off between data privacy and learning performance. 
In this experiment we aim to show that our algorithm \texttt{ObjG} outperforms the existing algorithm \texttt{OutG} while ensuring the same level of data privacy under the existence of convex constraints. 
We note that \texttt{OutG} proposed in \cite{huang2019dp, huang2020differentially} demonstrated better learning performance compared with other existing DP algorithms, such as DP-SGD in \cite{abadi2016deep} and DP-ADMM in \cite{zhang2016dynamic}, under the same level of data privacy, when there are no constraints to satisfy.

\subsubsection{Experimental settings} \label{sec:numerical_fl_exp}
We consider two publicly available datasets for image classification: MNIST \cite{lecun1998mnist} and FEMNIST \cite{caldas2018leaf}.
For the MNIST dataset, we split the 60,000 training data points over $P=10$ agents, each of which is assigned to have the same number of IID datasets.
For the FEMNIST dataset, we follow the preprocess procedure\footnote{https://github.com/TalwalkarLab/leaf/tree/master/data/femnist} to sample 5\% of the entire 805,263 data points in a non-IID manner, resulting in 36,708 training samples distributed over $P=195$ agents.
For local constraints $\mathcal{W}_p :=  [-u, u]$, we set $u=0.1$ for MNIST and $u=1.0$ for FEMNIST.
We set $\bar{\delta}=10^{-6}$ and consider various $\bar{\epsilon} \in \{0.05, 0.1, 0.5, 1\}$ where smaller $\bar{\epsilon}$ ensures stronger data privacy.
For \texttt{OutG}, we compute the $L_2$ sensitivity of the form \eqref{L2_sensitivity} as discussed in \cite{huang2019dp}.
We note that the $L_2$ sensitivity in \cite{huang2019dp} is derived from the optimality condition when $\mathcal{W}_p:=\mathbb{R}^n$.
It can be utilized for the simple box constraint considered in this experiment, but it cannot be utilized straightforwardly for a general convex $\mathcal{W}_p$. 
For \texttt{ObjG}, we compute the $L_2$ sensitivity of the form \eqref{redefined_L2sensitivity}.
Specifically, we add a data sample $(x_{pi^*}, y_{pi^*})$ to the existing dataset $\mathcal{D}_p$ to construct a neighboring dataset $\mathcal{D}_p'$. 
Then we can compute the $L_2$ sensitivity \eqref{redefined_L2sensitivity} based on \eqref{FL_local_obj} as follows:
\begin{align}
\bar{\Delta}^{t,e}_{p,2} = \max_{ \forall (\mathcal{D}_p, \mathcal{D}_p')}  \Big[ \sum_{j=1}^J \sum_{k=1}^K \Big\{ \frac{1}{I+1}   x_{pi^* j} \Big(-y_{pi^*k} + \frac{  \exp(\sum_{j'=1}^J  x_{pi^*j} z^{t,e}_{p j' k})}{\sum_{k'=1}^K \exp(\sum_{j'=1}^J  x_{pi^*j'} z^{t,e}_{pjk'})}  \Big) \Big\}^2 \Big]^{0.5}, \label{fl_sensitivity}
\end{align} 
which is computed by finding $(x_{pi^*}, y_{pi^*})$ from the universe of datasets that maximizes \eqref{fl_sensitivity} for given $z^{t,e}_p$.

A parameter $\rho^t$ is chosen after conducting some tuning process discussed in Appendix \ref{apx:hyper}. We note that the chosen parameter $\hat{\rho}^t$ is nondecreasing and bounded above, thus satisfying Assumption \ref{assump:convergence} (i).

We implemented the algorithms in Python and ran the experiments on Swing, a 6-node GPU computing cluster at Argonne National Laboratory. Each node of Swing has 8 NVIDIA A100 40 GB GPUs, as well as 128 CPU cores.

\begin{figure*}[!ht]  
\centering
\begin{subfigure}[b]{0.23\textwidth}
\centering        
\includegraphics[width=\textwidth]{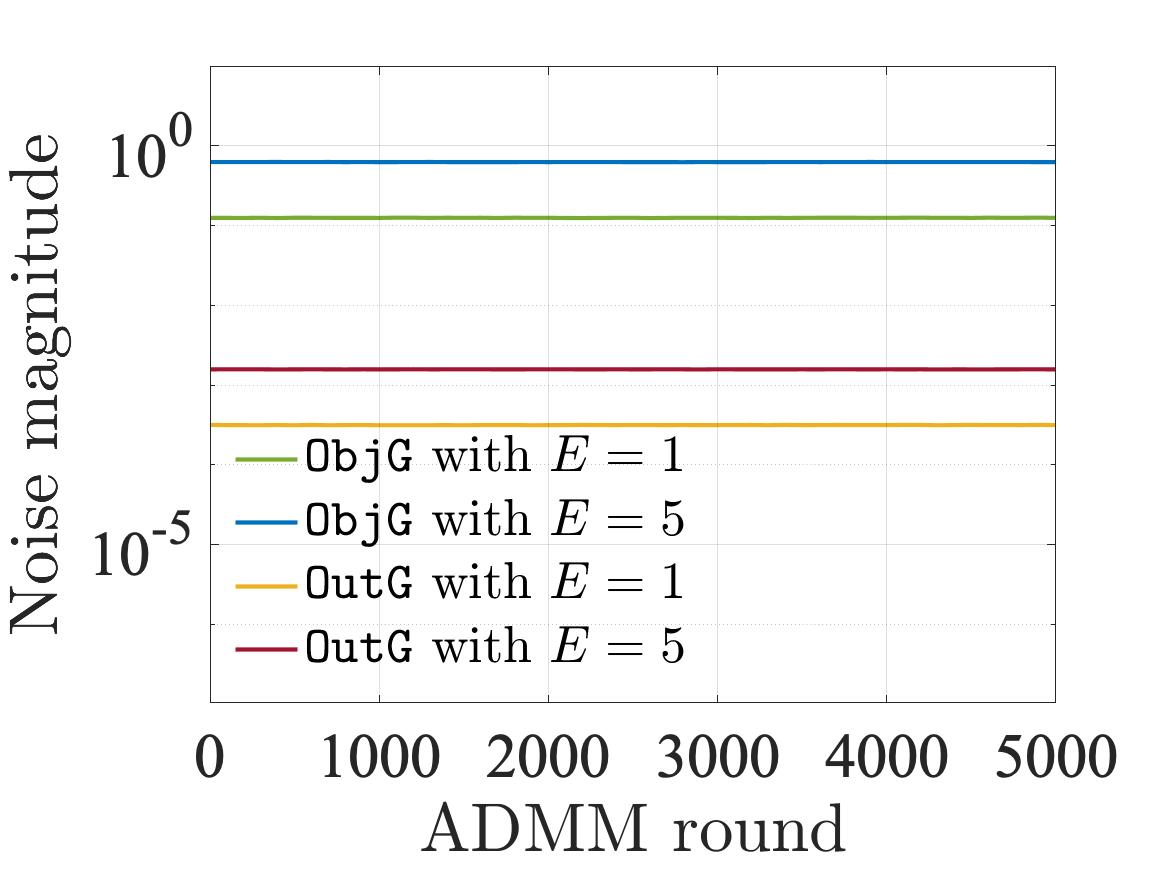}
\includegraphics[width=\textwidth]{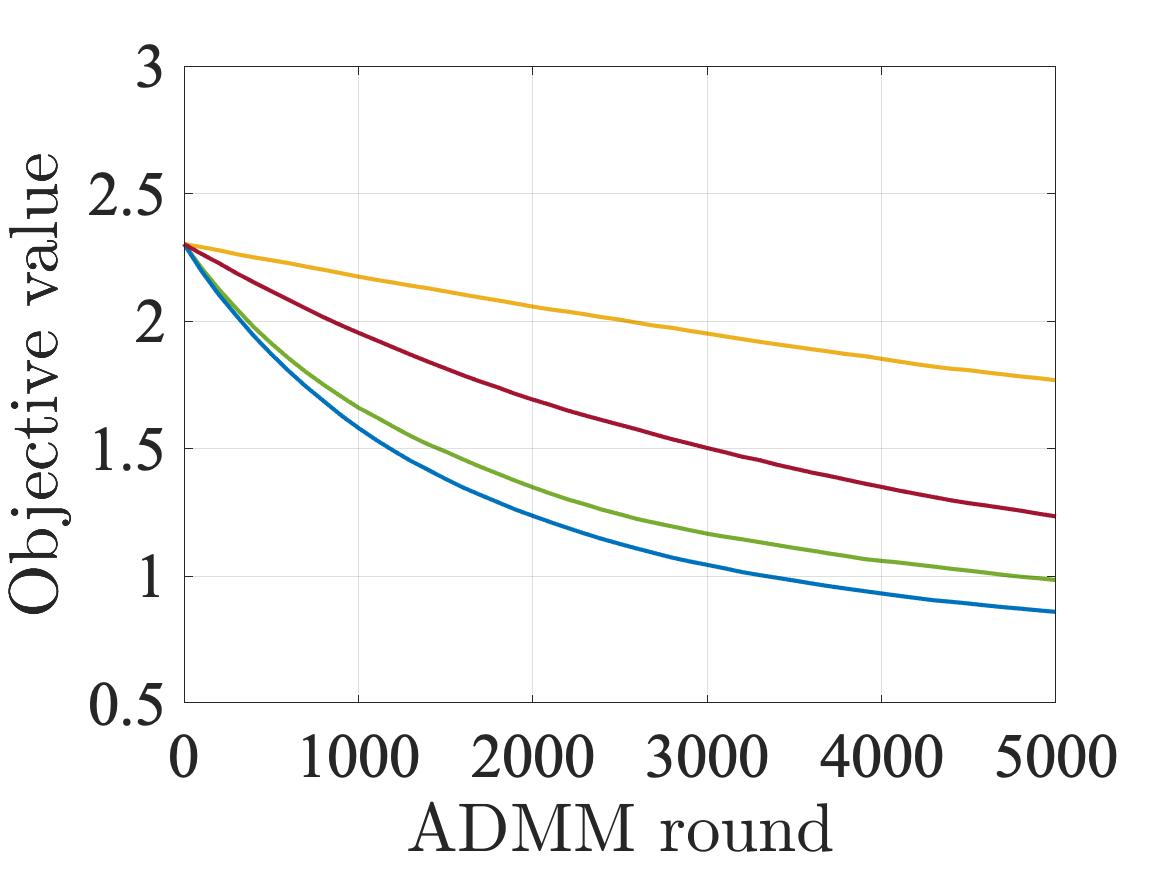}  
\includegraphics[width=\textwidth]{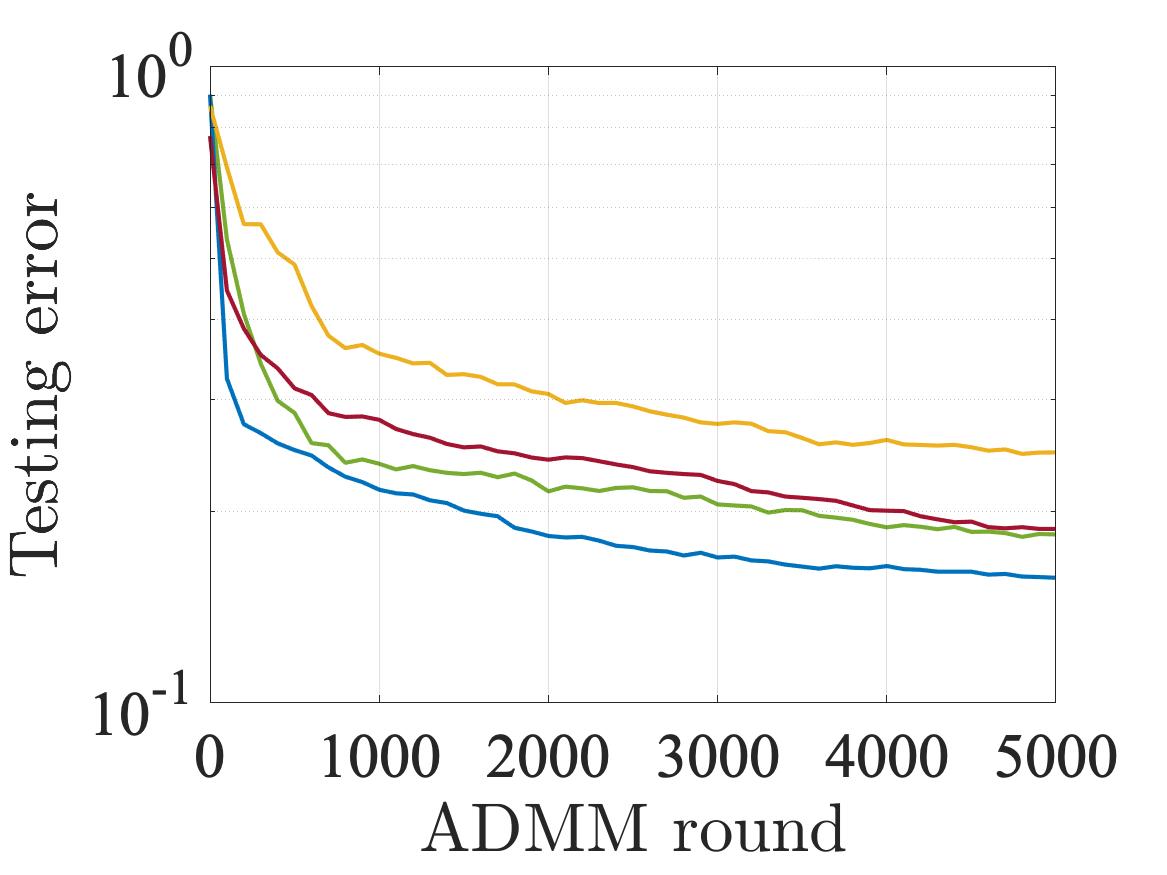}  
\caption{$\bar{\epsilon}=0.05$}
\end{subfigure}
\begin{subfigure}[b]{0.23\textwidth}
\centering        
\includegraphics[width=\textwidth]{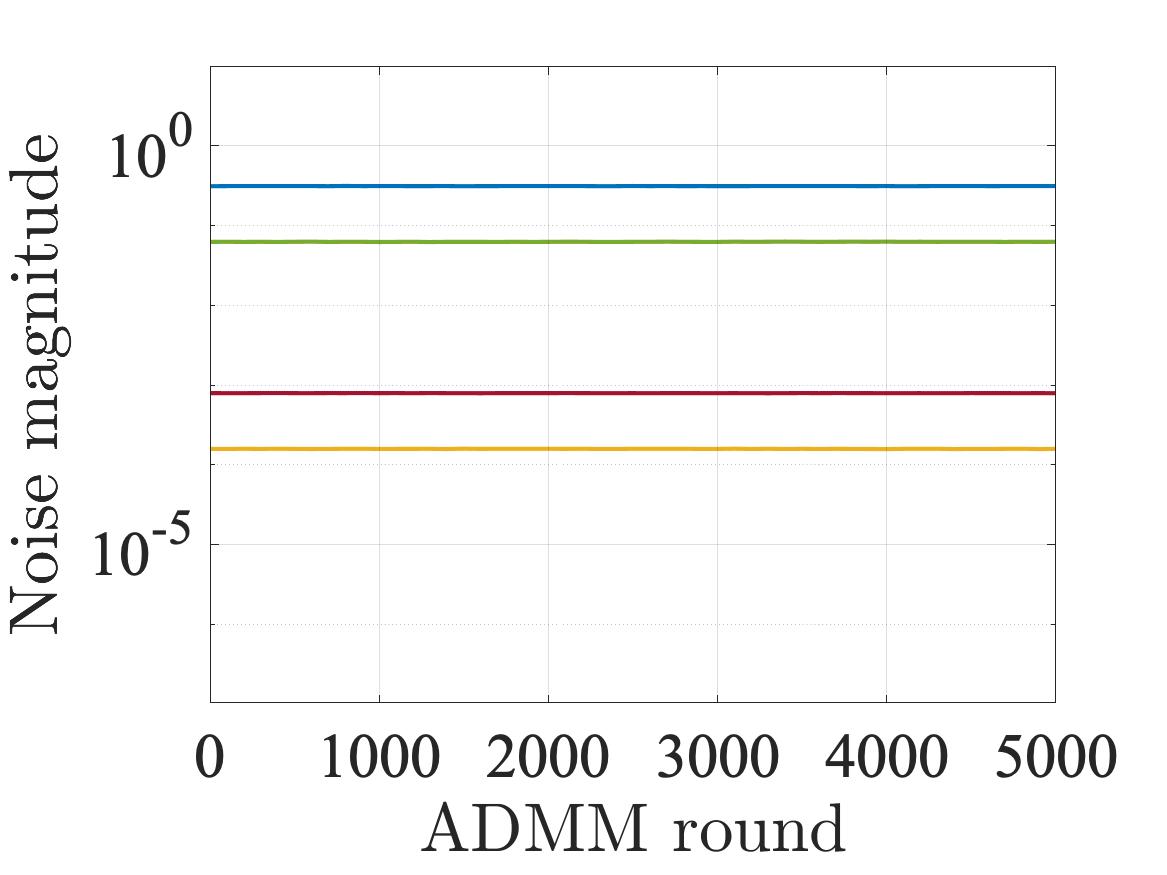}
\includegraphics[width=\textwidth]{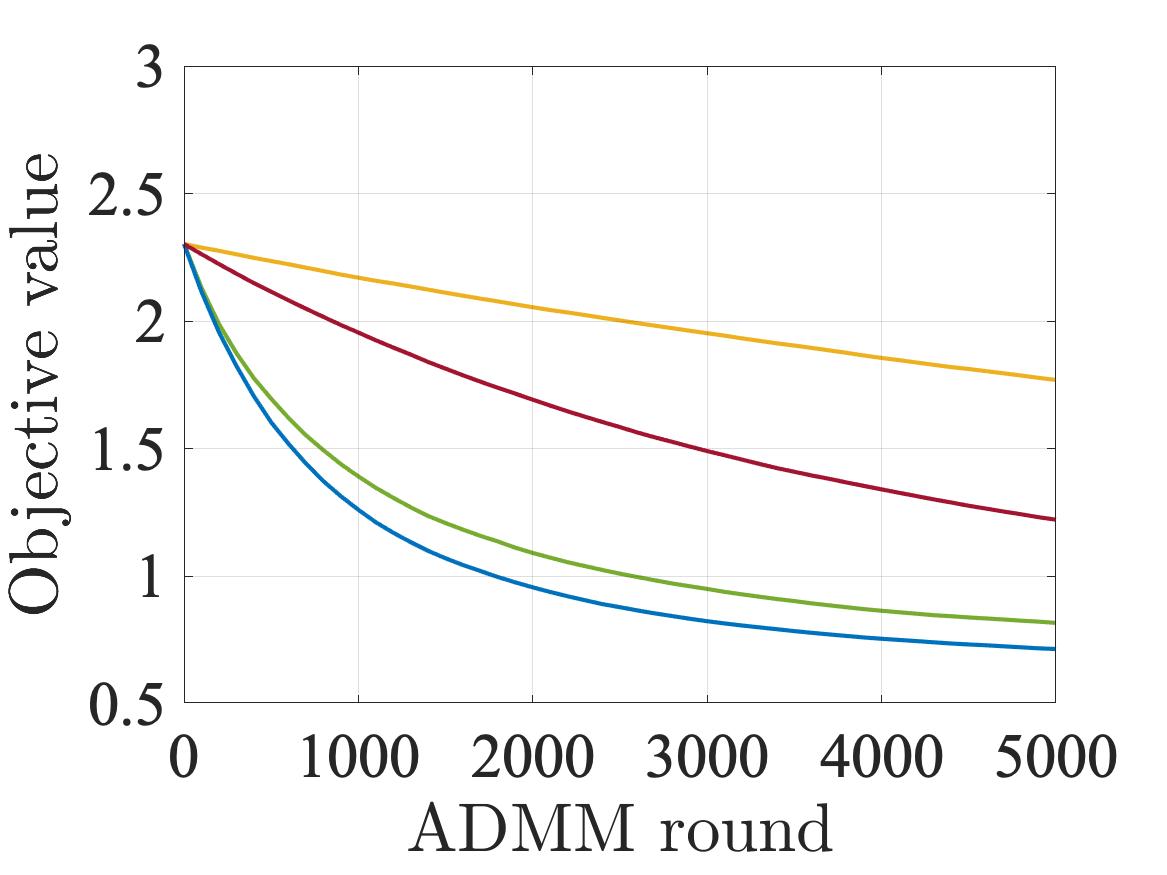}  
\includegraphics[width=\textwidth]{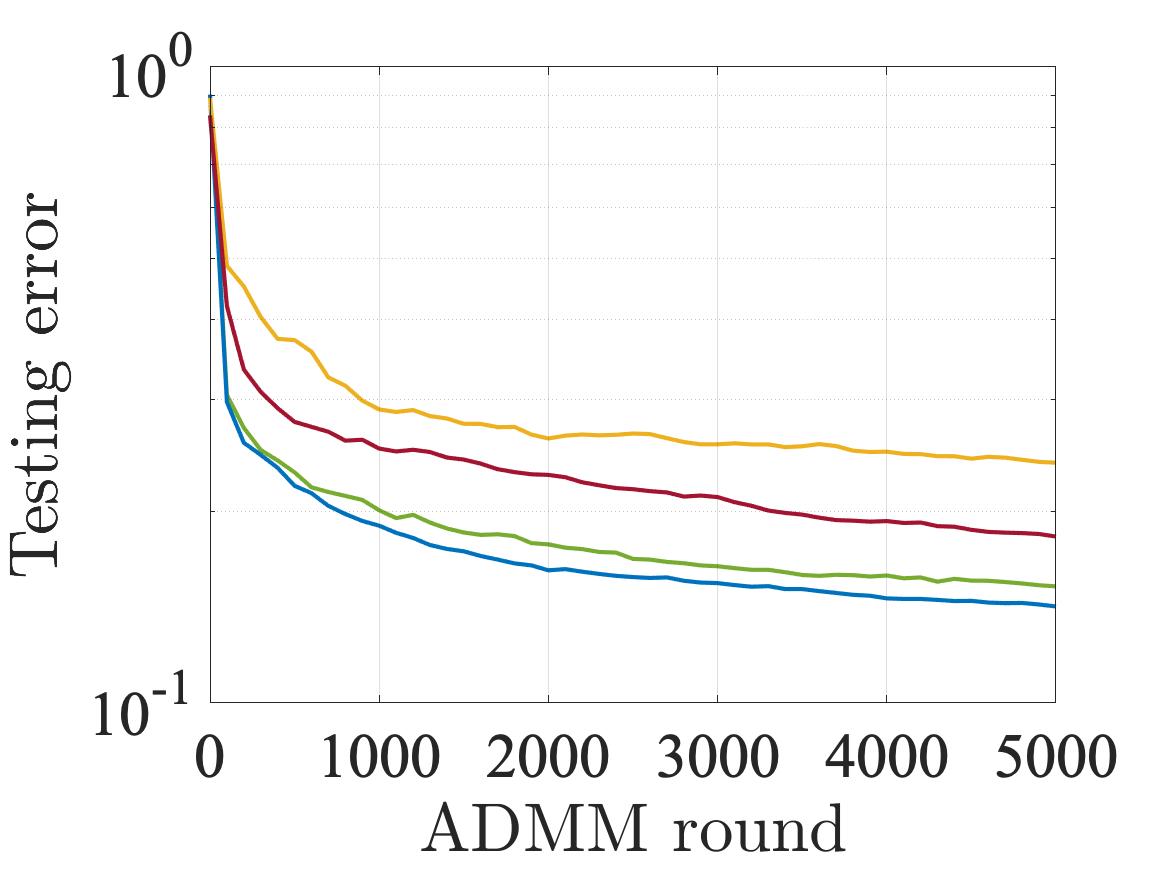}  
\caption{$\bar{\epsilon}=0.1$}
\end{subfigure}
\begin{subfigure}[b]{0.23\textwidth}
\centering      
\includegraphics[width=\textwidth]{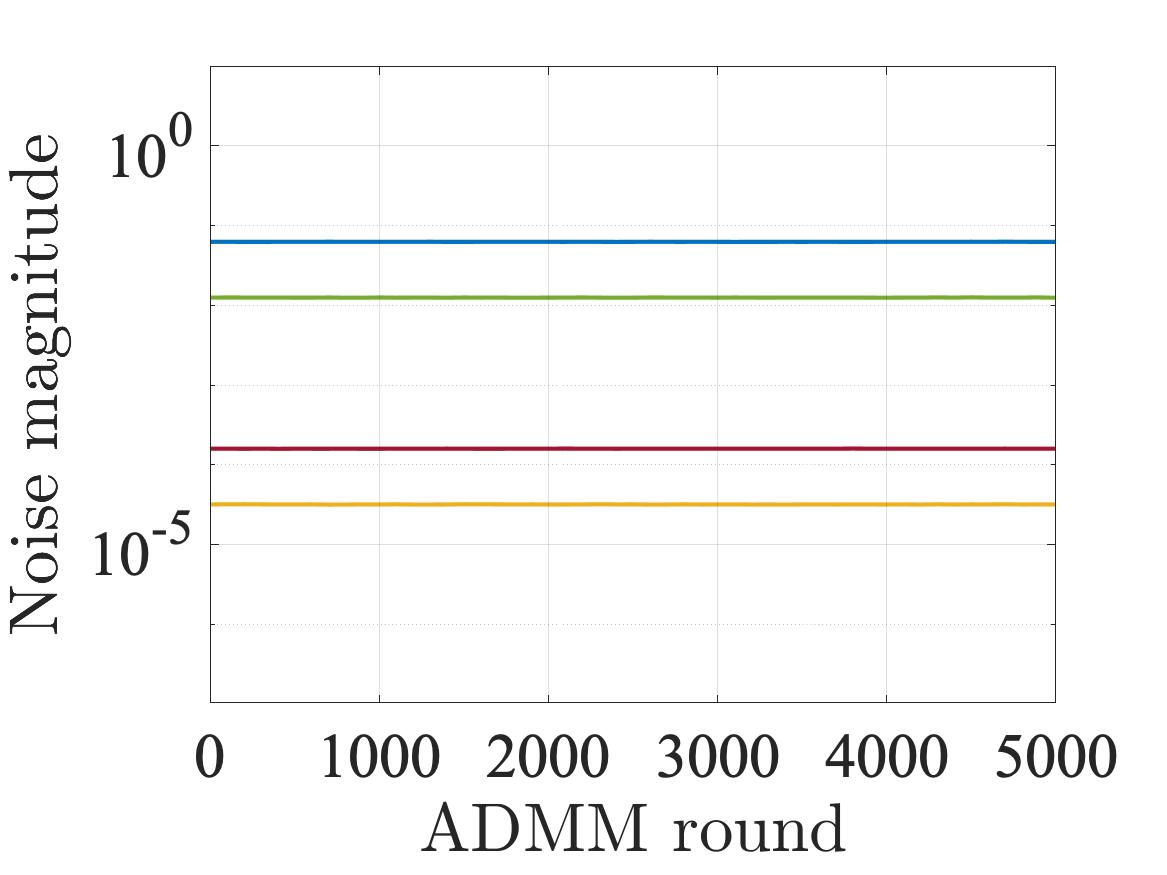}
\includegraphics[width=\textwidth]{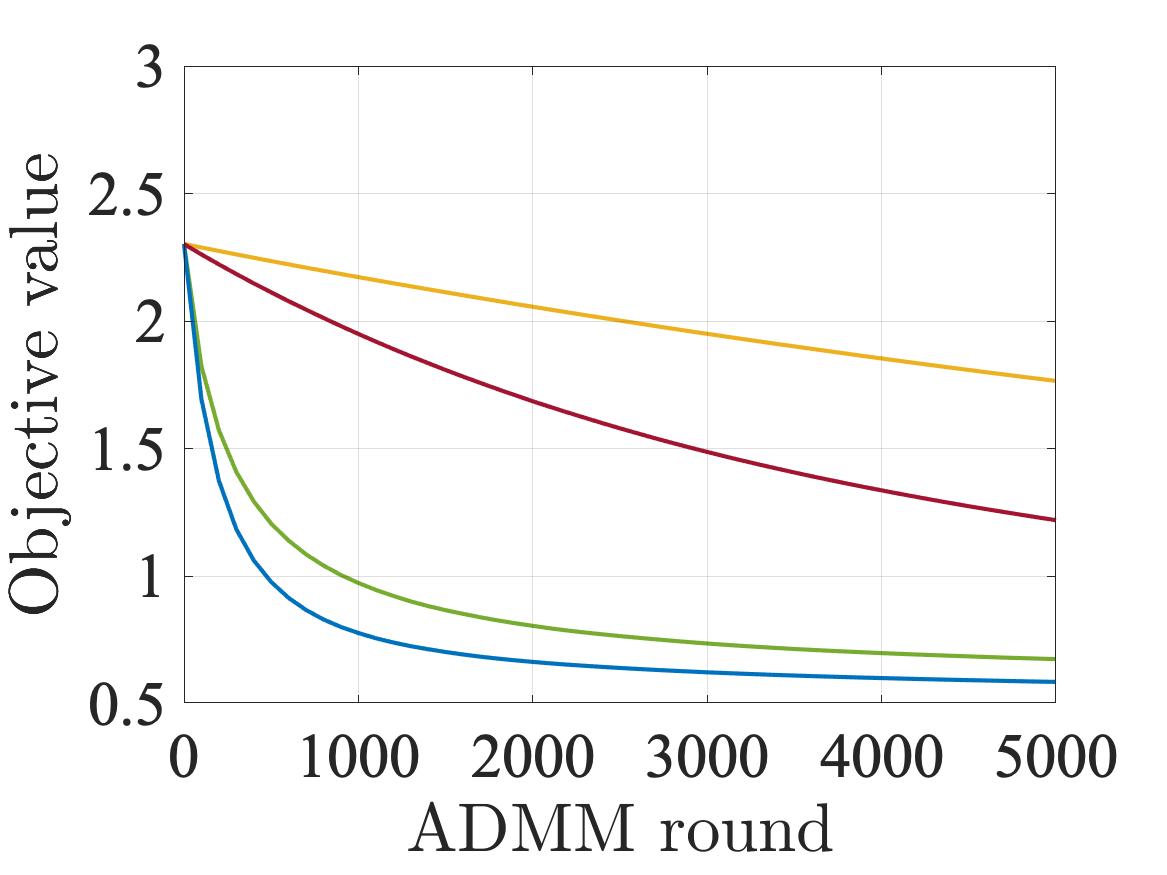}  
\includegraphics[width=\textwidth]{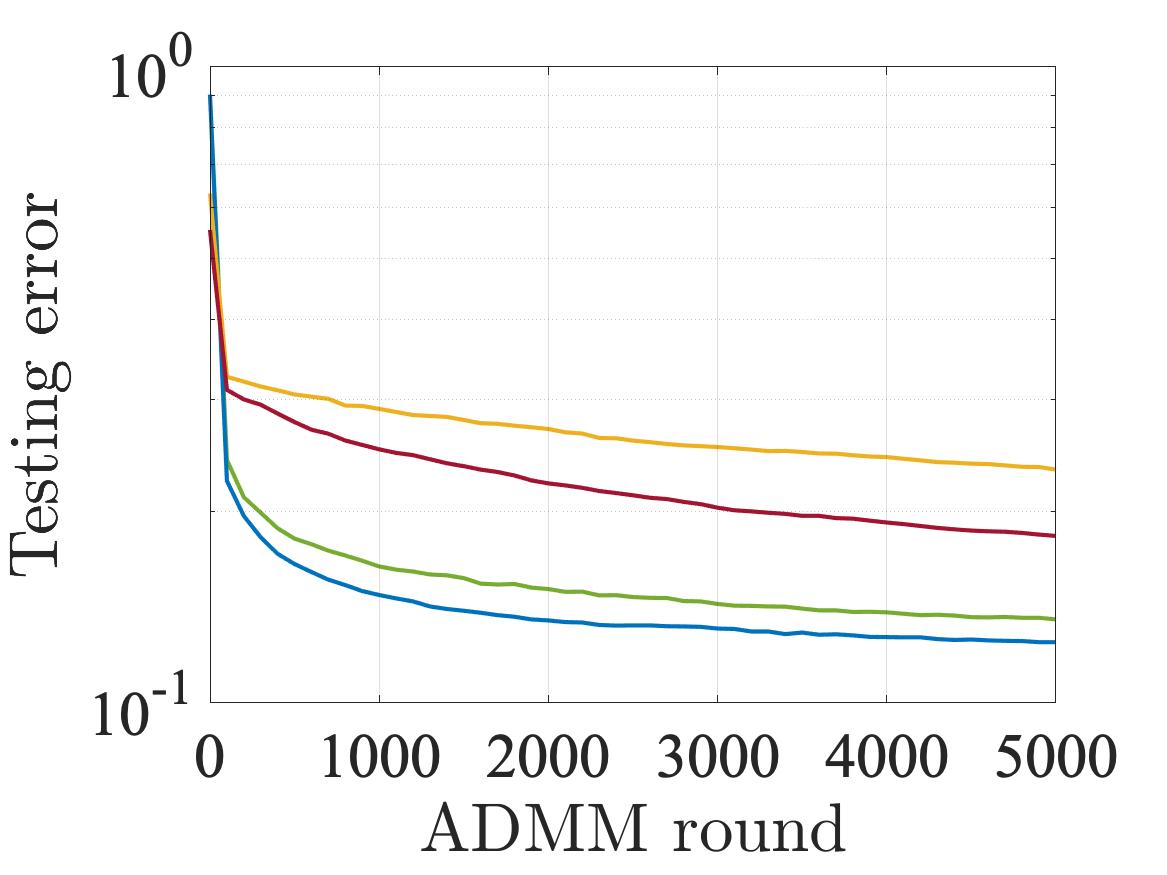}  
\caption{$\bar{\epsilon}=0.5$}
\end{subfigure}     
\begin{subfigure}[b]{0.23\textwidth}
\centering      
\includegraphics[width=\textwidth]{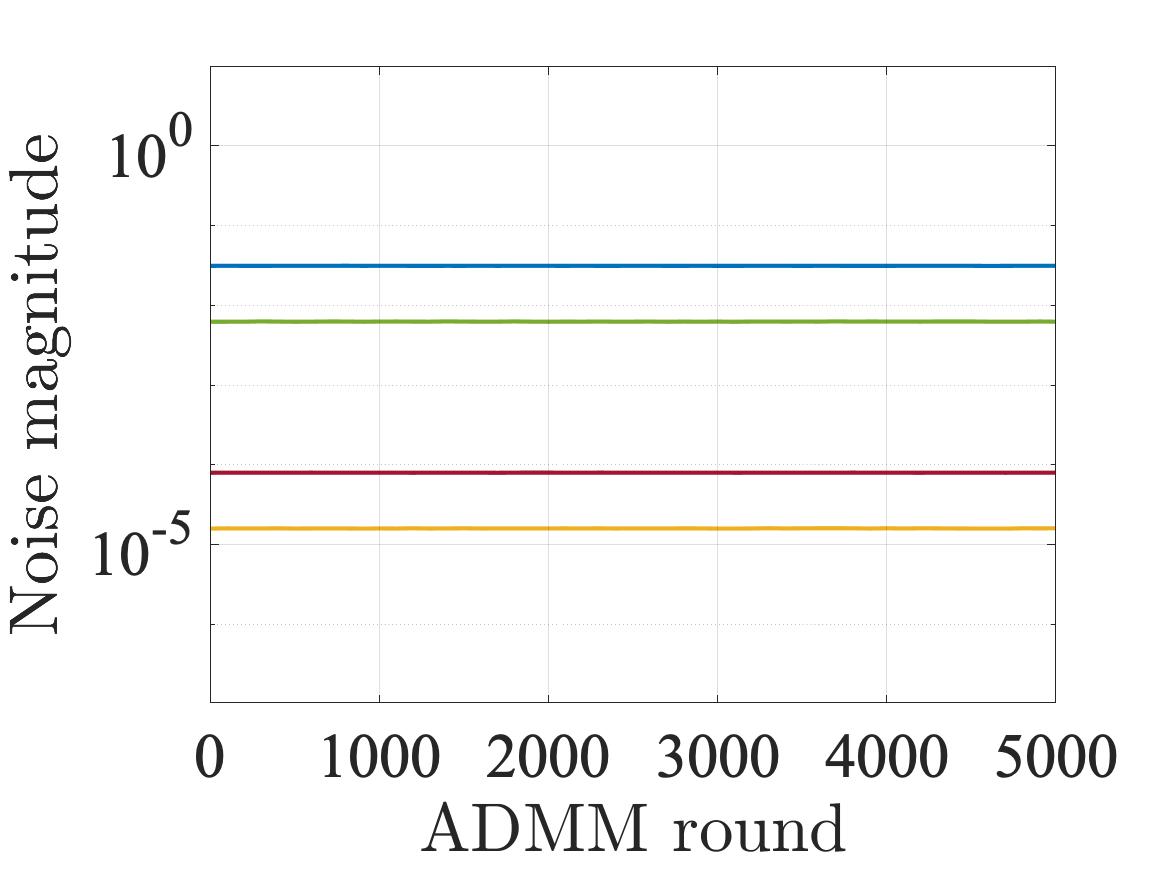}
\includegraphics[width=\textwidth]{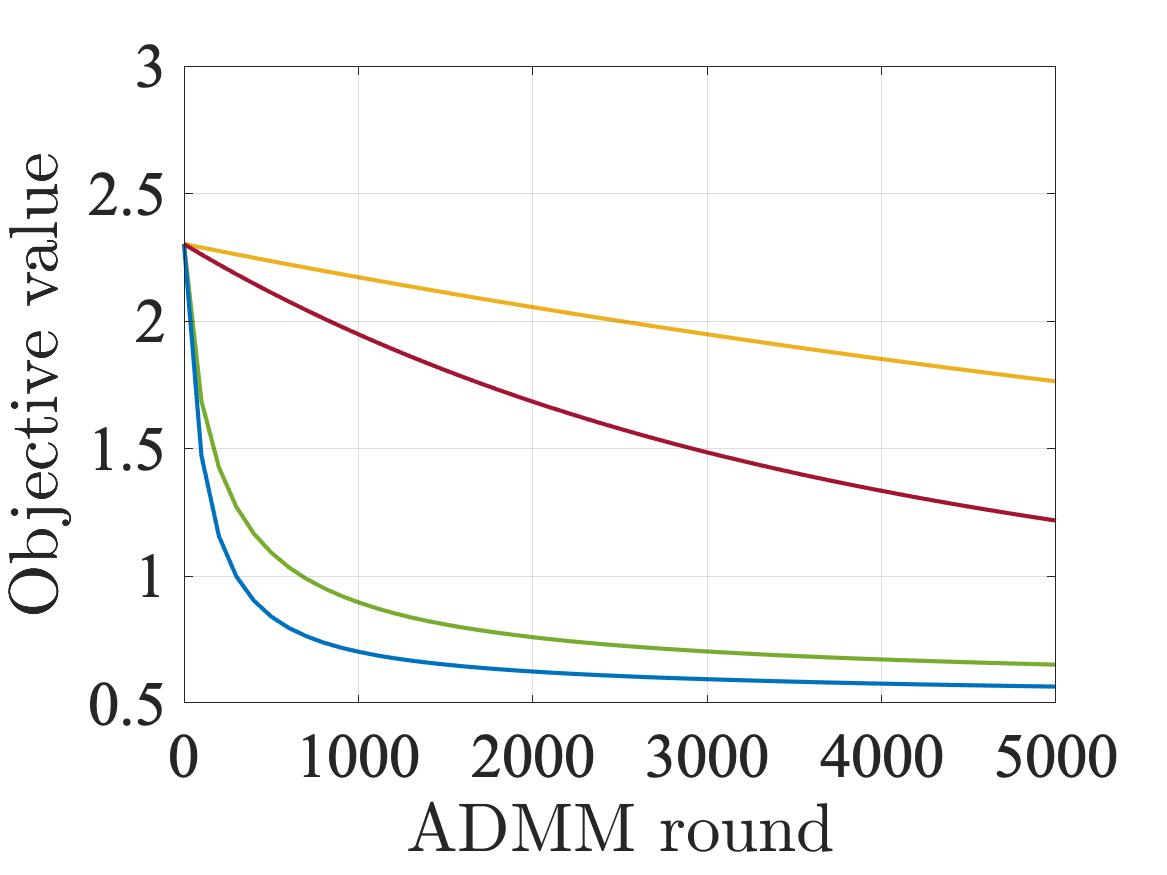}  
\includegraphics[width=\textwidth]{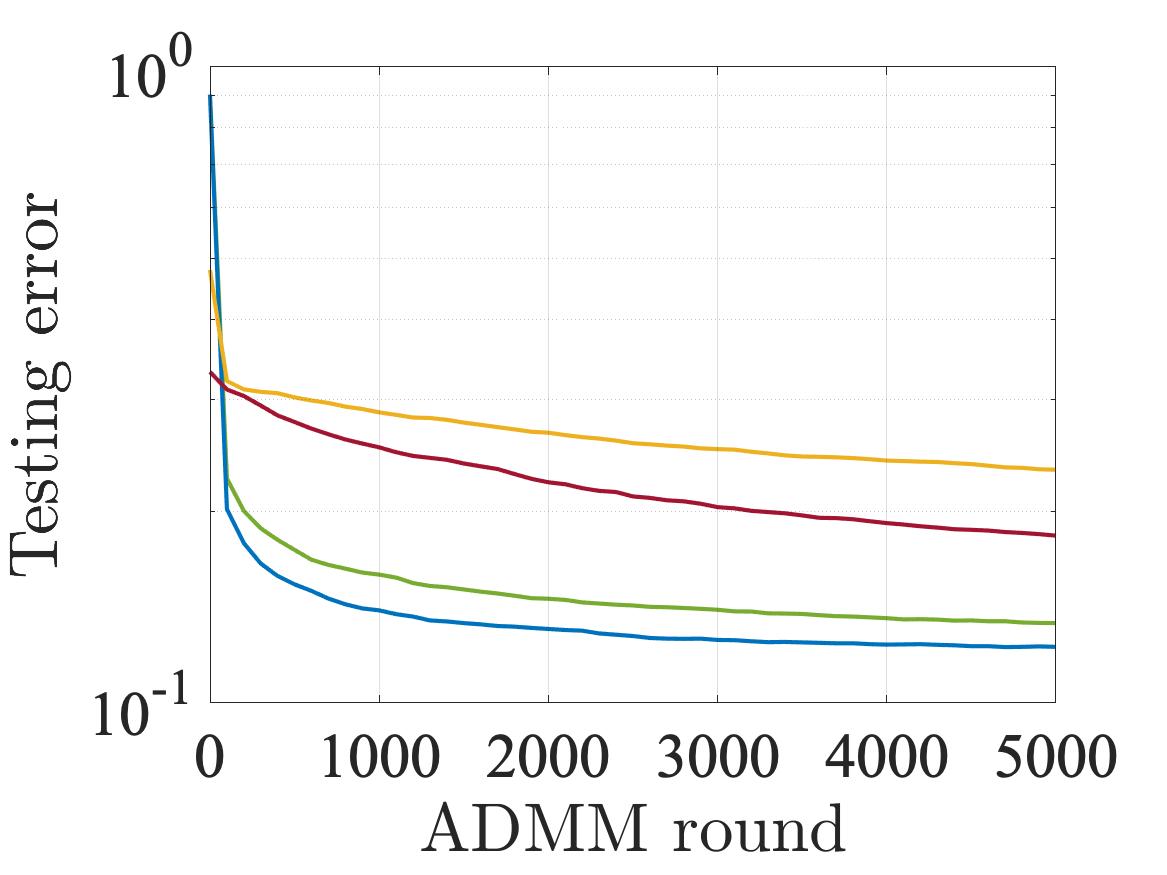}  
\caption{$\bar{\epsilon}=1$}
\end{subfigure}           
\caption{Comparison of \texttt{ObjG} and \texttt{OutG} with respect to the noise magnitude ($1^{\text{st}}$ row), objective value ($2^{\text{nd}}$ row), and testing error ($3^{\text{rd}}$ row) by using MNIST.}
\label{fig:MNIST}
\end{figure*}

\begin{figure*}[!ht]  
  \centering
  \begin{subfigure}[b]{0.23\textwidth}
  \centering        
  \includegraphics[width=\textwidth]{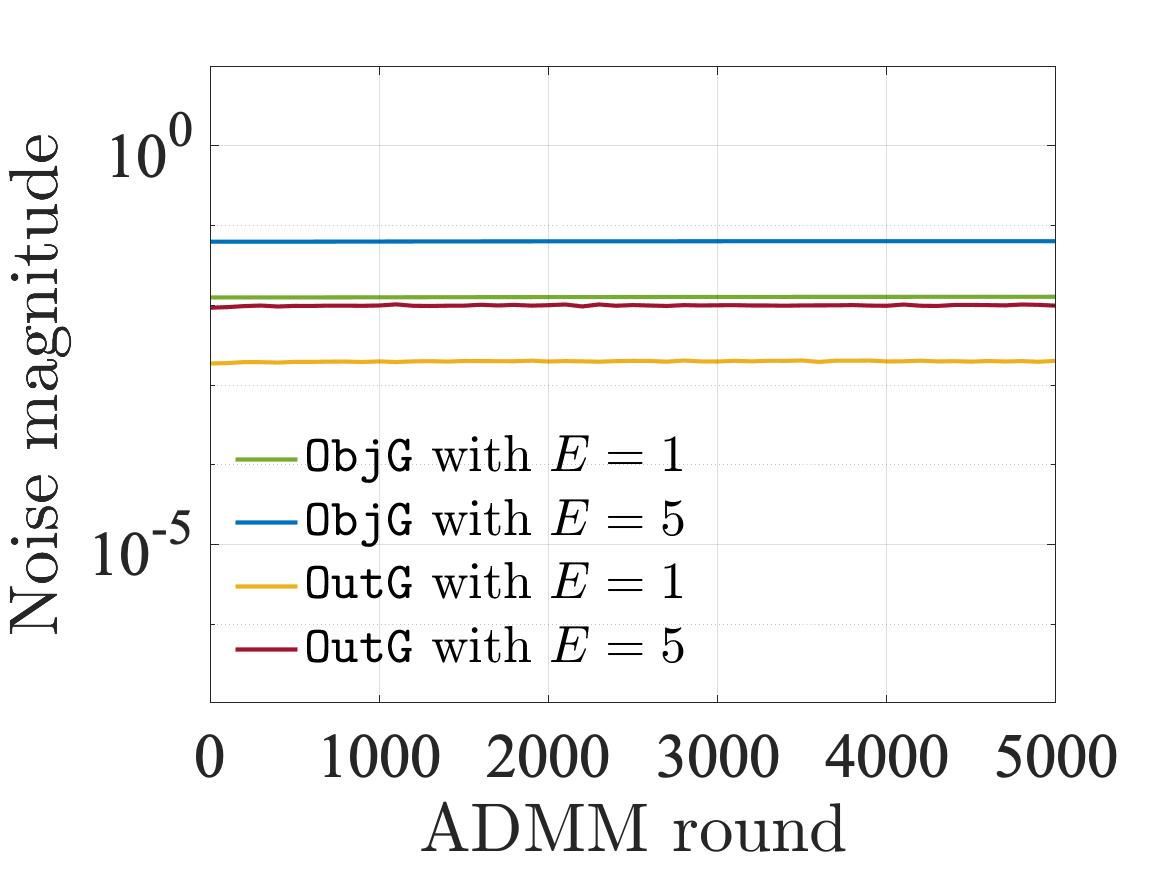}  
  \includegraphics[width=\textwidth]{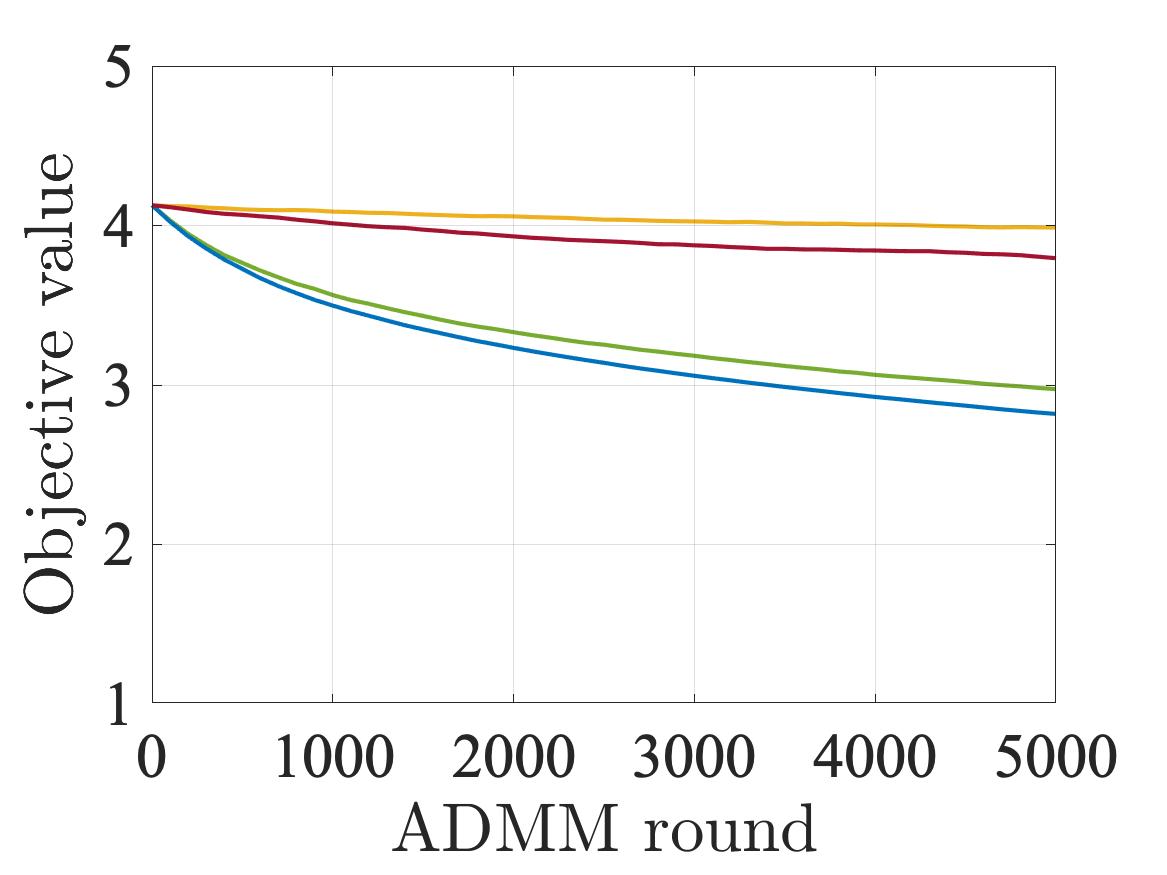}  
  \includegraphics[width=\textwidth]{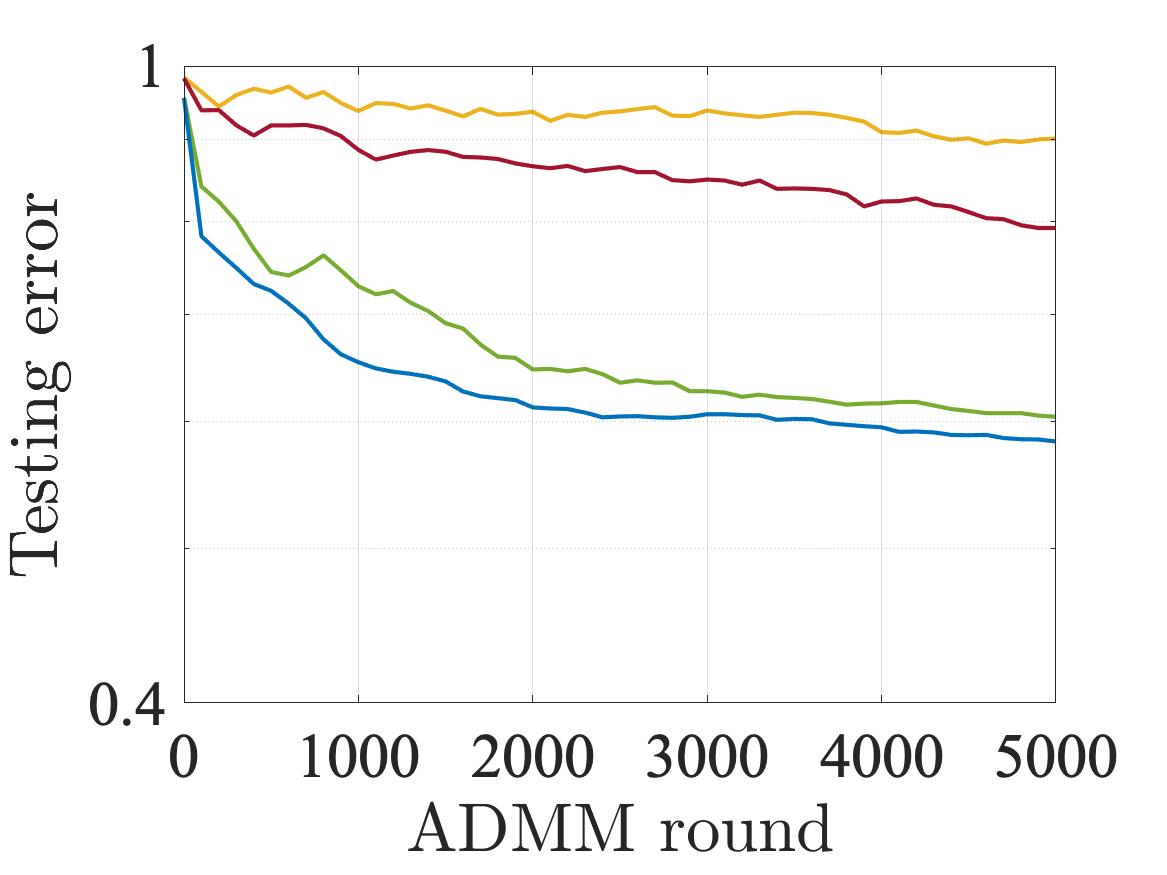}  
  \caption{$\bar{\epsilon}=0.05$}
  \end{subfigure}
  \begin{subfigure}[b]{0.23\textwidth}
  \centering        
  \includegraphics[width=\textwidth]{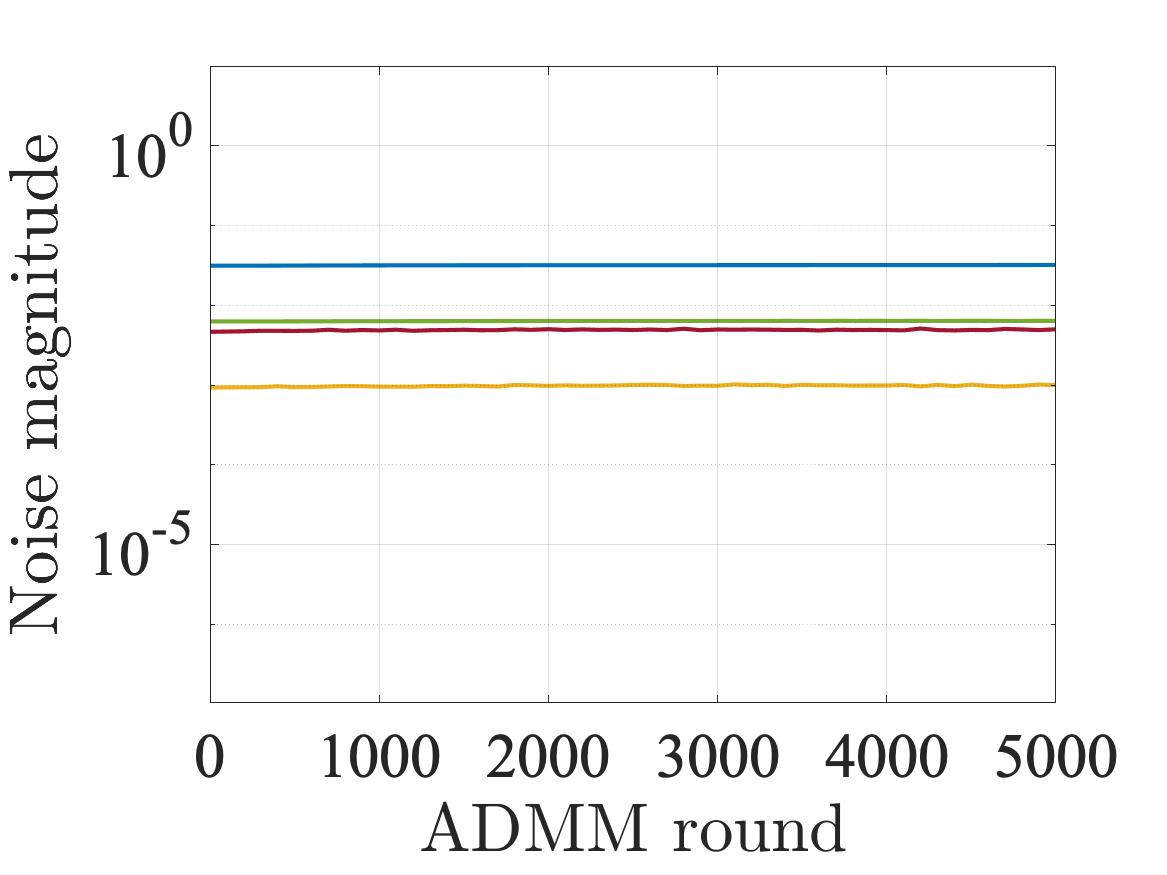}  
  \includegraphics[width=\textwidth]{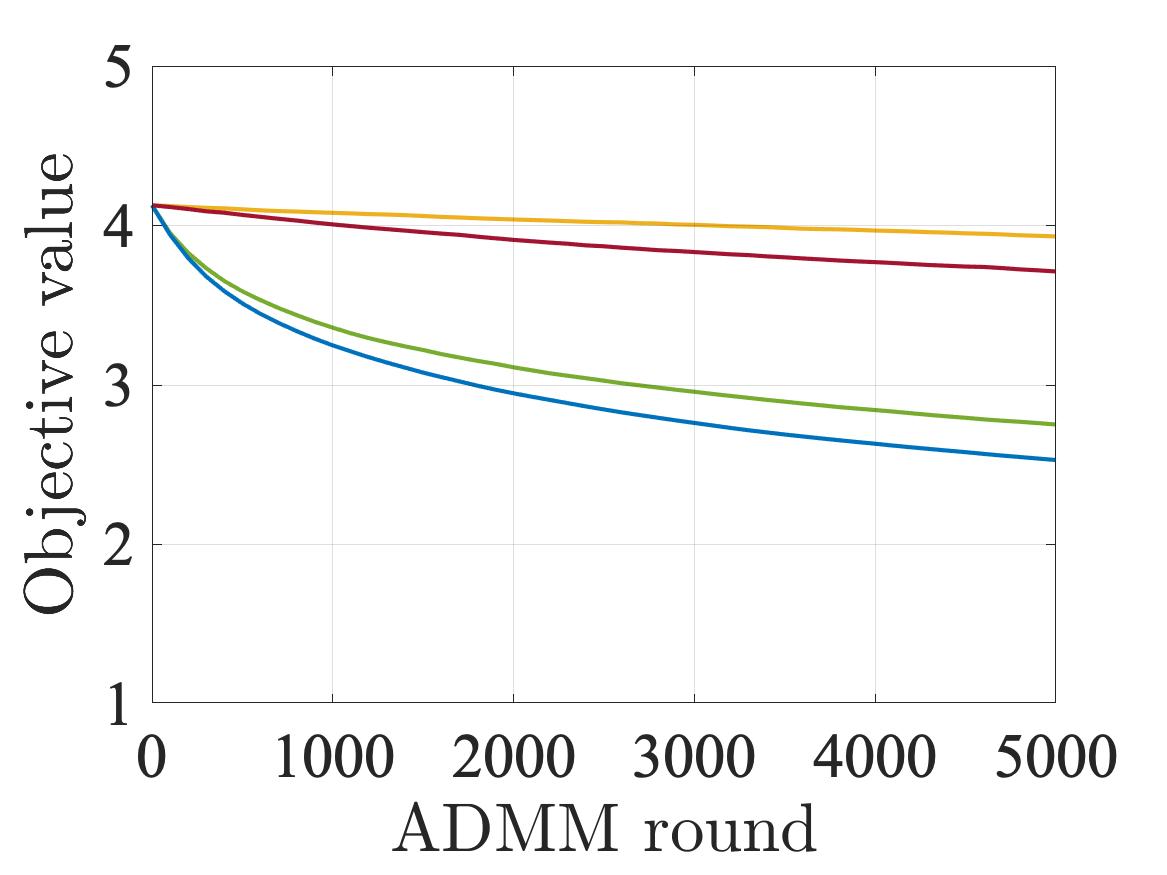}  
  \includegraphics[width=\textwidth]{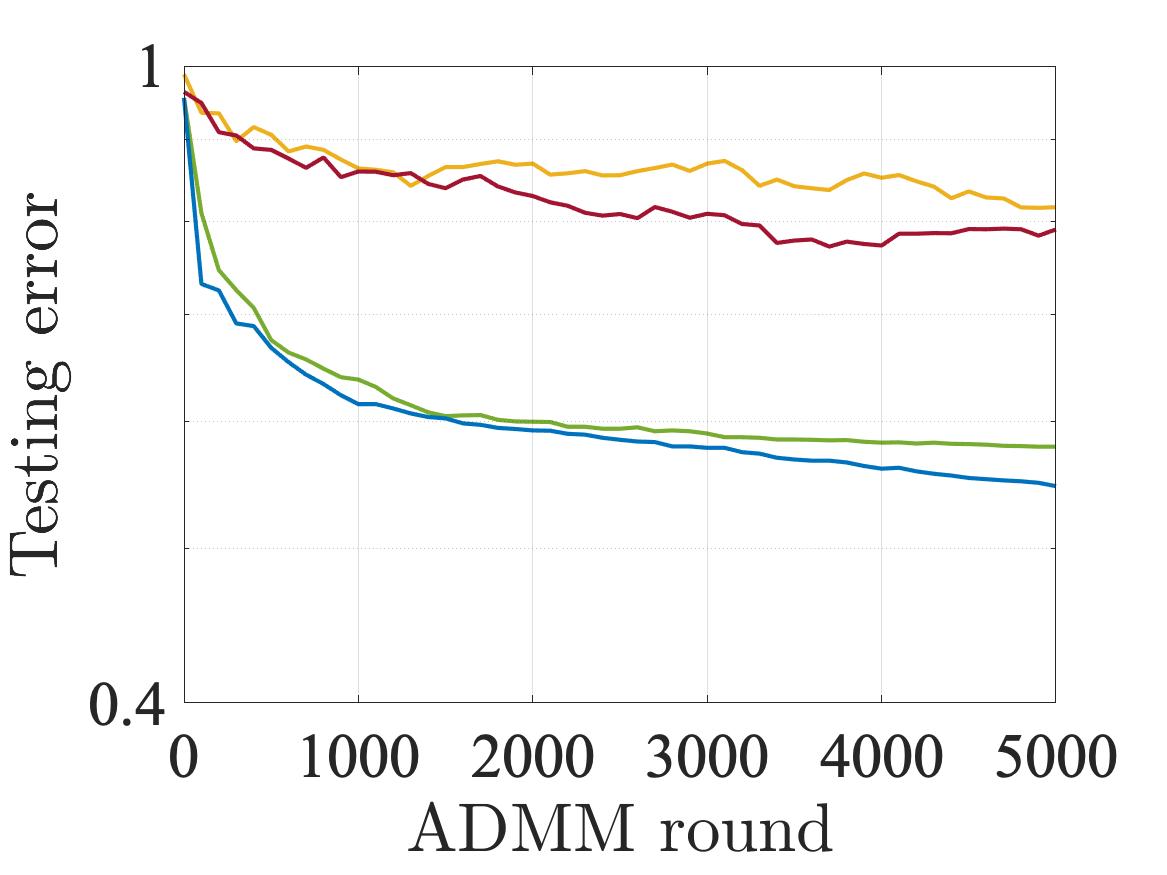}    
  \caption{$\bar{\epsilon}=0.1$}
  \end{subfigure}
  \begin{subfigure}[b]{0.23\textwidth}
  \centering      
  \includegraphics[width=\textwidth]{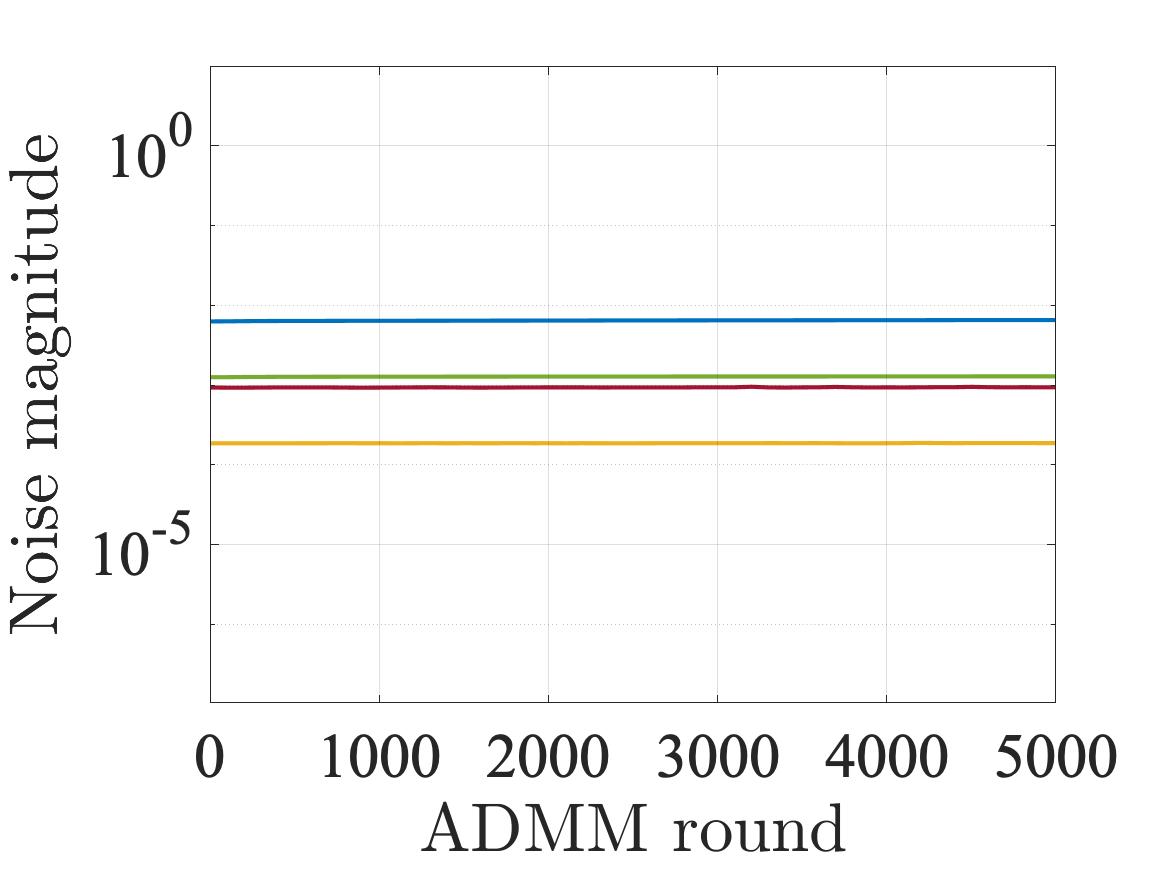}  
  \includegraphics[width=\textwidth]{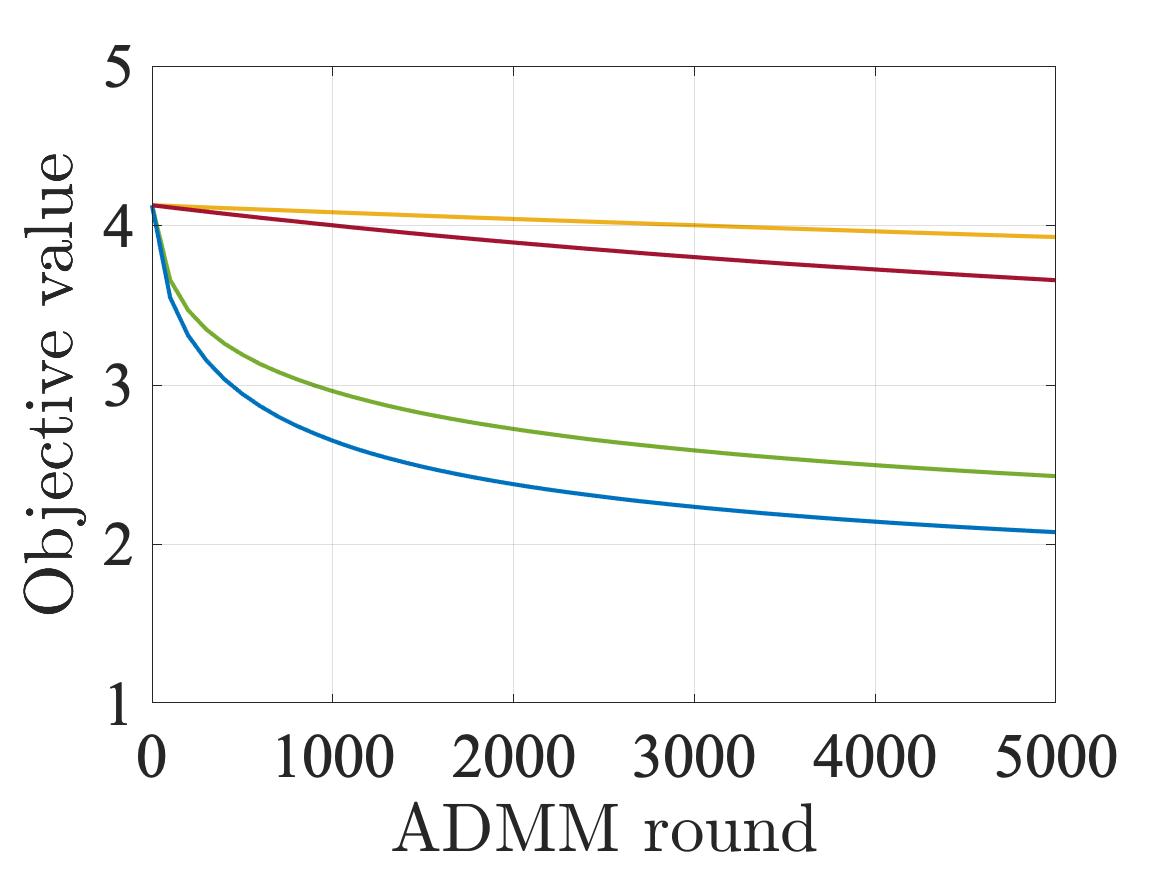}  
  \includegraphics[width=\textwidth]{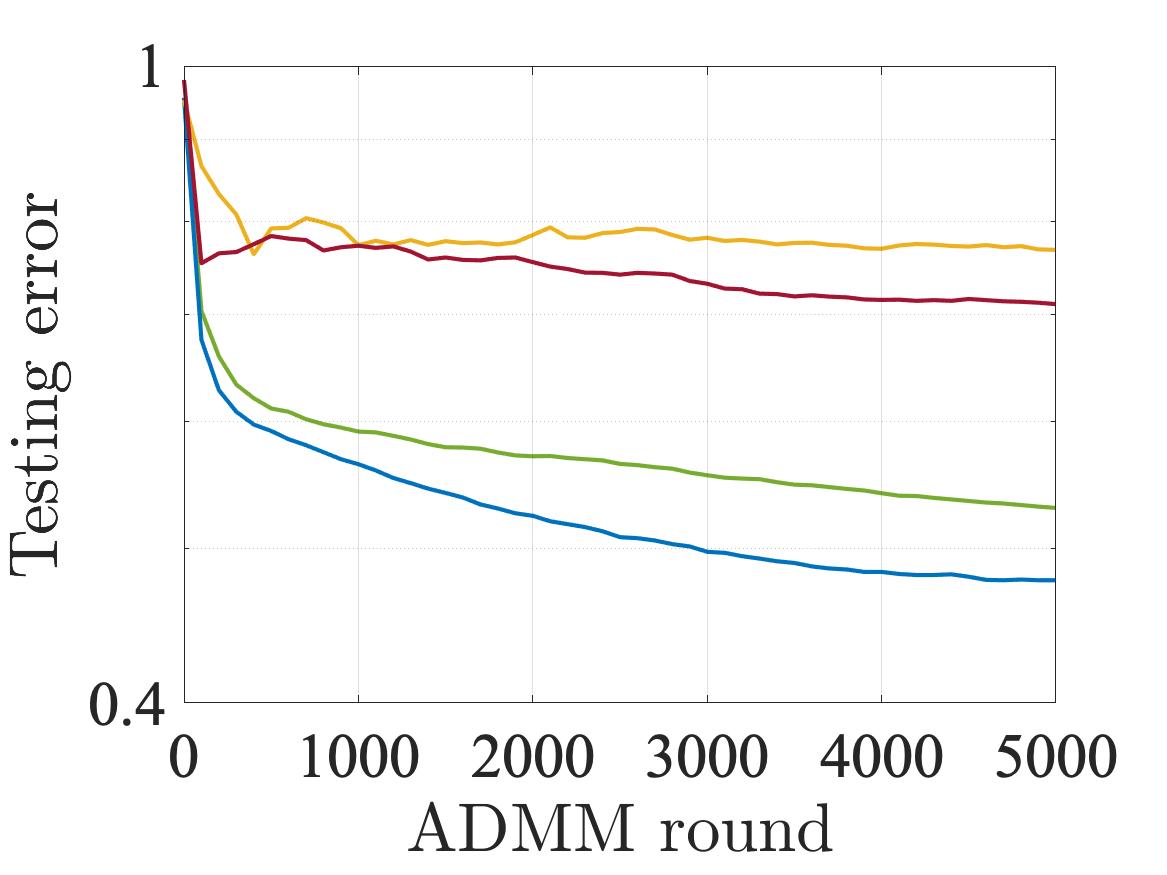}  
  \caption{$\bar{\epsilon}=0.5$}
  \end{subfigure}     
  \begin{subfigure}[b]{0.23\textwidth}
  \centering      
  \includegraphics[width=\textwidth]{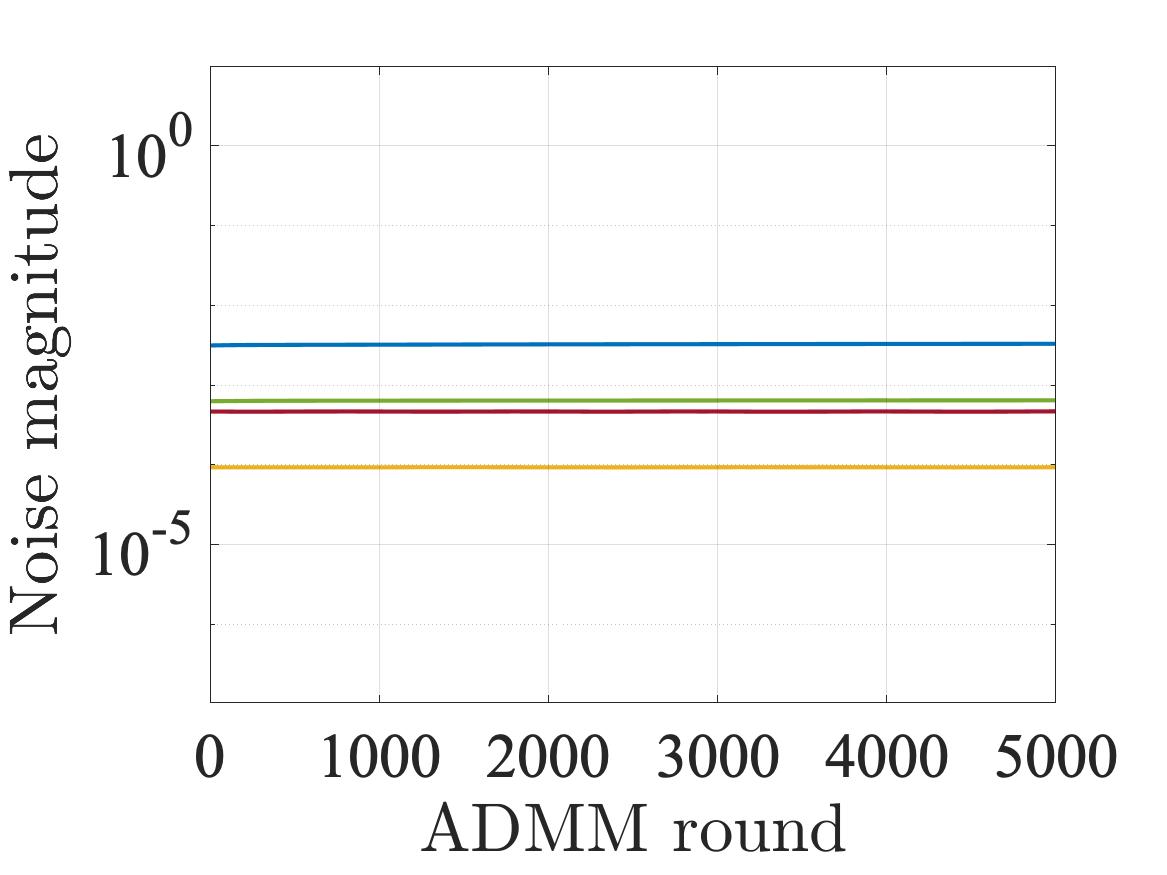}  
  \includegraphics[width=\textwidth]{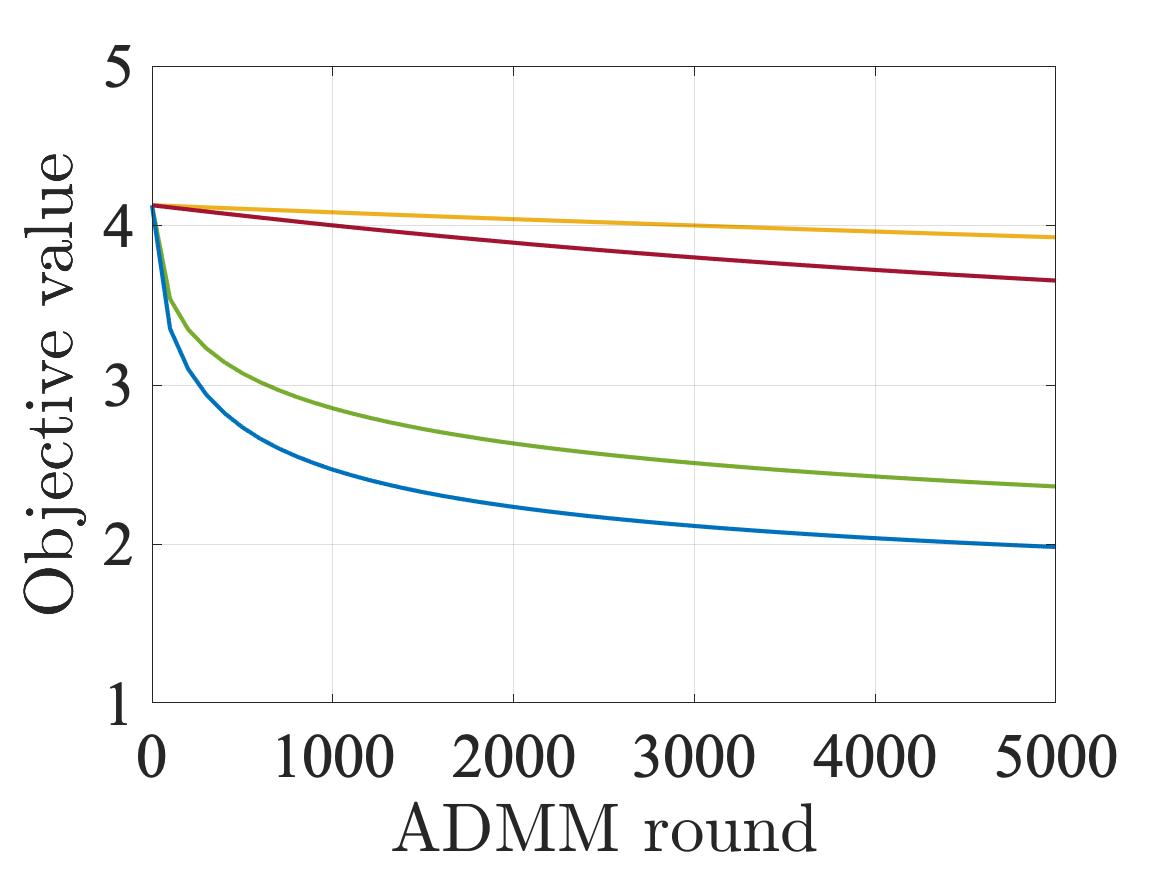}  
  \includegraphics[width=\textwidth]{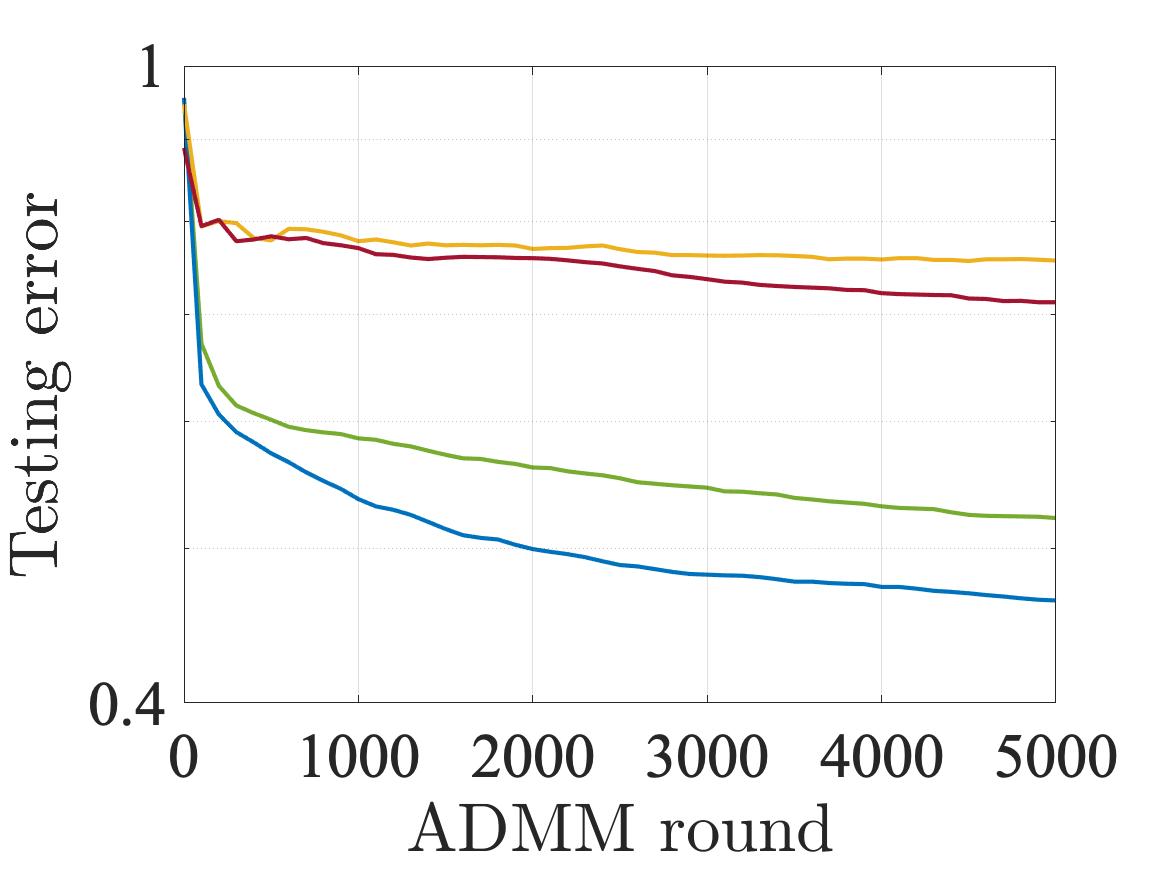}  
  \caption{$\bar{\epsilon}=1$}
  \end{subfigure}           
  \caption{Comparison of \texttt{ObjG} and \texttt{OutG} with respect to the noise magnitude ($1^{\text{st}}$ row), objective value ($2^{\text{nd}}$ row), and testing error ($3^{\text{rd}}$ row) by using FEMNIST.}
  \label{fig:FEMNIST}
  \end{figure*}

\subsubsection{Performance comparison} 
We report the performance of the proposed algorithm \texttt{ObjG} and the existing algorithm \texttt{OutG} for MNIST and FEMNIST image classification in Figure \ref{fig:MNIST} and \ref{fig:FEMNIST}, respectively.

First, we report in the first row of Figures \ref{fig:MNIST} and \ref{fig:FEMNIST} the maginitude of noise introduced during $T=5000$ ADMM rounds, which is defined as
\begin{align}
  \sum_{p=1}^P \sum_{e=1}^E \Big(  \frac{1}{n}  \sum_{i=1}^n | \tilde{\xi}^{t,e}_{p i} | \Big), \ \forall t \in [T],
\end{align}
where $\tilde{\xi}^{t,e}_{pi} \in \mathbb{R}$ is sampled from a Gaussian distribution with zero mean and variance $\sigma^2 = 2 \ln (1.25/\bar{\delta}) ( \bar{\Delta}^{t,e}_{p,2} / \bar{\epsilon})^2$, where $\bar{\Delta}^{t,e}_{p,2}$ is the $L_2$ sensitivity computed as described in Section \ref{sec:numerical_fl_exp}.
We observe that the noise magnitude for \texttt{ObjG} is greater than that for \texttt{OutG} mainly because the sensitivity of \texttt{ObjG} is greater.

Second, we report the objective values (i.e., training costs) produced by \texttt{ObjG} and \texttt{OutG} in the second row of Figures \ref{fig:MNIST} and \ref{fig:FEMNIST}, respectively.
We observe that the objective value increases as $\bar{\epsilon}$ decreases, which indicates the trade-off between data privacy and learning performance, well known in the literature on DP algorithms \cite{dwork2014algorithmic}.
Also, the objective values of \texttt{ObjG} are lower than those of \texttt{OutG}.
Furthermore, increasing $E$ lowers the objective values, showing the efficacy of the multiple local updates.
  
Third, we report the testing errors produced by \texttt{ObjG} and \texttt{OutG} in the third row of Figures \ref{fig:MNIST} and \ref{fig:FEMNIST}, respectively.
We observe that the testing errors of \texttt{ObjP} with multiple local updates significantly outperforms the other algorithms. 
This result implies that our algorithm can mitigate the trade-off between data privacy and learning performance.   

\section{Conclusion} \label{sec:conclusion}
In this paper we proposed a differentially private distributed optimization algorithm for general convex optimization problems. 
Our algorithm addresses a data privacy concern in distributed optimization by communicating randomized outputs of the local optimization models, which prevents reconstructing data stored at the local machine.
Different from the existing DP algorithms in the literature, most of which have been developed for unconstrained machine learning models, our algorithm is developed for constrained convex optimization models, which can be beneficial because constraints are necessary in most optimization control problems and which become more important in some machine learning applications.

Specifically, our algorithm generalizes the existing linearized ADMM by introducing multiple local updates to reduce communication costs and by incorporating the objective perturbation method to the local optimization models so that the resulting randomized outputs are always feasible (i.e., satisfying the constraints).
This is different from the popular output perturbation method used in the existing DP algorithms where the resulting randomized outputs constructed by adding noise to the true output, which could be infeasible. 
We presented two different mechanisms, namely, Laplace and Gaussian mechanisms, that ensure differential privacy for every iteration of the proposed algorithm.
These are similar to the Laplace and Gaussian mechanisms used in the output perturbation method except for the sensitivity computation. 
We also provided a convergence analysis that showed that the rate of convergence in expectation is sublinear without any error bound.  
We discussed two different application areas, namely, distributed control of power flow and federated learning, where the proposed algorithms can be utilized.
From the experiments, we numerically demonstrate the outperformance of the proposed algorithm over the existing DP algorithms with the output perturbation method, when there are constraints to satisfy.

We note that the privacy and convergence analyses in this paper are limited to the convex setting and thus not directly applicable to the nonconvex setting.

\newpage

\appendix 
 
\section{Proof of Theorem \ref{thm:smooth}} \label{apx-thm:smooth} 
\subsection{Preliminaries}
\subsubsection{Some constants} \label{sec:constants}
First, $U_1$ is well defined under Assumption \ref{assump:convergence}-(iii).
The necessary and sufficient condition of Assumption \ref{assump:convergence}-(iii) is that, for all $u \in \mathcal{W}_p$  and $v \in \partial f_p(u)$, $\| v \|_{\star} \leq H$, where $\|\cdot\|_{\star}$ is the dual norm.
Since the dual norm of the Euclidean norm is the Euclidean norm, we have $\| f'_p(u) \| \leq H$.
Since the objective function, which is a maximum of finite continuous functions, is continuous and $\mathcal{W}_p$ is compact, $U_1$ is well defined. 

Second, $U_2$ is well defined because the objective function $\|u-v\|$ is continuous and the feasible region $\mathcal{W}_p$ is compact. 

Third, $U_3$ is well defined because of the existence of $U_1$.

\subsubsection{Basic equations.}
We note that for any symmetric matrix $A$, 
\begin{align}
(a-b)^{\top} A (c-d) = \frac{1}{2} \{ \| a-d \|^2_A  - \| a-c \|^2_A + \| c-b \|^2_A - \| d - b \|^2_A \},  \label{rule}
\end{align}
where $a$, $b$, $c$, and $d$ are vectors of the same size.

We define $\tilde{\lambda}^t_p := \lambda^t_p + \rho^t (w^{t+1} - z^t_p)$ for fixed $t \in [T]$ and $p \in [P]$.
From the optimality condition of \eqref{ADMM-1}, namely, 
$\sum_{p=1}^P \lambda^t_p + \rho^t(w^{t+1}-z^t_p)= \sum_{p=1}^P \tilde{\lambda}^t_p =0$, we have
\begin{align}
&  \sum_{p=1}^P \langle  \tilde{\lambda}^t_p, w^{t+1}-w \rangle = 0, \ \forall w \in \mathbb{R}^n. \label{Block1:optimality_condition}
\end{align}

\subsection{Inequality derivation for a fixed iteration $t$ and $e$}
For a given $p \in [P]$, the optimality condition of \eqref{DPADMM-2-Prox} is given by 
\begin{align*}
& \langle \nabla f_p(z^{t,e}_p) - \{ \underbrace{ \lambda^t_p  + \rho^t(w^{t+1}-z^{t,e+1}_p) }_{ \text{``\texttt{A}''} } \} + \tilde{\xi}^{t,e}_p, z^{t,e+1}_p - z_p \rangle  \nonumber \\
\leq  & \frac{1}{\eta^{t}} \underbrace{\langle z_p - z^{t,e+1}_p, z^{t,e+1}_p - z_p^{t,e} \rangle}_{\text{``\texttt{B}''}}, \ \forall z_p \in \mathcal{W}_p.
\end{align*}
By defining $\lambda^{t,e+1}_p  := \lambda^t_p + \rho^t (w^{t+1} - z^{t,e+1}_p)$ for the ``\texttt{A}'' term and applying \eqref{rule} on the ``\texttt{B}'' term from the above inequalities, we have
\begin{align}
& \langle \nabla f_p(z^{t,e}_p) - \lambda^{t,e+1}_p + \tilde{\xi}^{t,e}_p, z^{t,e+1}_p - z_p \rangle \label{Block2:smooth_inequality_0} \\
\leq & \frac{1}{2\eta^{t}} \Big( \| z_p - z_p^{t,e}\|^2 - \| z_p - z_p^{t,e+1}\|^2 - \| z_p^{t,e} - z_p^{t,e+1} \|^2 \Big). \nonumber
\end{align}
By adding a term $\langle \lambda^{t,e+1}_p, z^{t,e+1}_p-z_p \rangle$ to the gradient inequality $f_p(z^{t,e}_p) - f_p(z_p)  \leq \langle \nabla f_p(z^{t,e}_p), z^{t,e}_p - z_p \rangle$ for all $z_p$, we derive
\begin{align}
& f_p(z^{t,e}_p) - f_p(z_p) - \underbrace{\langle \lambda^{t,e+1}_p, z^{t,e+1}_p-z_p \rangle}_{\text{``\texttt{C}''}}  \nonumber \\
\leq &  \langle \nabla f_p(z^{t,e}_p), z^{t,e}_p-z^{t,e+1}_p \rangle + \langle \nabla f_p(z^{t,e}_p) - \lambda^{t,e+1}_p, z^{t,e+1}_p - z_p \rangle \nonumber \\
= &  \langle  \nabla f_p(z^{t,e}_p) + \tilde{\xi}^{t,e}_p, z^{t,e}_p-z^{t,e+1}_p \rangle + \langle  \tilde{\xi}^{t,e}_p, z_p - z^{t,e}_p \rangle \nonumber \\
& + \underbrace{\langle \nabla f_p(z^{t,e}_p) - \lambda^{t,e+1}_p + \tilde{\xi}^{t,e}_p, z^{t,e+1}_p-z_p \rangle}_{\text{applying } \eqref{Block2:smooth_inequality_0}}  \nonumber \\
\leq & \langle  \nabla f_p(z^{t,e}_p) + \tilde{\xi}^{t,e}_p, z^{t,e}_p-z^{t,e+1}_p \rangle  + \langle  \tilde{\xi}^{t,e}_p, z_p - z^{t,e}_p \rangle  \nonumber \\
& + \frac{1}{2\eta^{t}} \Big( \| z_p - z_p^{t,e}\|^2 - \| z_p - z_p^{t,e+1}\|^2 - \| z_p^{t,e} - z_p^{t,e+1} \|^2 \Big).\nonumber 
\end{align}
Since the ``\texttt{C}'' term from the above inequalities can be written as
\begin{align}
& \langle \lambda^{t,e+1}_p, z^{t,e+1}_p-z_p \rangle = \langle \lambda^{t+1}_p, z^{t,e+1}_p-z_p \rangle + \langle \underbrace{\lambda^{t,e+1}_p - \lambda^{t+1}_p}_{=\rho^t(z^{t+1}_p - z^{t,e+1}_p ) }, z^{t,e+1}_p-z_p \rangle, \nonumber
\end{align}
we obtain
\begin{align}
& f_p(z^{t,e}_p) - f_p(z_p) - \langle \lambda^{t+1}_p, z^{t,e+1}_p-z_p \rangle \label{Block2:smooth_inequality_1} \\
\leq & \rho^t \langle z^{t+1}_p - z^{t,e+1}_p, z^{t,e+1}_p-z_p \rangle + \underbrace{\langle \nabla f_p(z^{t,e}_p) + \tilde{\xi}^{t,e}_p, z^{t,e}_p-z^{t,e+1}_p \rangle}_{\text{``\texttt{D}''}} + \langle  \tilde{\xi}^{t,e}_p, z_p - z^{t,e}_p \rangle   \nonumber \\
& + \frac{1}{2\eta^{t}} \Big( \| z_p - z_p^{t,e}\|^2 - \| z_p - z_p^{t,e+1}\|^2 - \| z_p^{t,e} - z_p^{t,e+1} \|^2 \Big), \ \forall z_p \in \mathcal{W}_p.   \nonumber
\end{align}
Next, we derive from the ``\texttt{D}'' term in \eqref{Block2:smooth_inequality_1} that
\begin{align}
& \langle \nabla f_p(z^{t,e}_p) + \tilde{\xi}^{t,e}_p, z^{t,e}_p-z^{t,e+1}_p \rangle\\
= & \underbrace{\langle \tilde{\xi}^{t,e}_p, z^{t,e}_p-z^{t,e+1}_p \rangle}_{ \text{applying Young's inequality} }  + \underbrace{\langle \nabla f_p(z^{t,e}_p), z^{t,e}_p-z^{t,e+1}_p \rangle}_{\text{applying the L-smoothness of $f_p$}}  \nonumber \\
\leq & \Big\{
\frac{1}{2(1/\eta^{t} - L)} \| \tilde{\xi}^{t,e}_p \|^2 +  \frac{1/\eta^{t} - L} {2} \| z_p^{t,e+1} - z^{t,e}_p\|^2  \Big\}  \nonumber \\
& + \Big\{ f_p(z_p^{t,e}) - f_p(z_p^{t,e+1}) + \frac{L}{2}\| z_p^{t,e+1} - z_p^{t,e} \|^2 \Big\}, \nonumber 
\end{align}
where $1/\eta^t - L > 0$ by the construction of $\eta^t$ in \eqref{smooth_proximity}.
Therefore, we derive from \eqref{Block2:smooth_inequality_1} the following inequalities:
\begin{align}
  & f_p(z^{t,e+1}_p) - f_p(z_p) - \langle  \lambda^{t+1}_p, z^{t,e+1}_p-z_p \rangle  \label{Block2:smooth_inequality_2} \\
  \leq & \rho^t \langle z^{t+1}_p - z^{t,e+1}_p, z^{t,e+1}_p-z_p \rangle +   \frac{1}{2(1/\eta^{t} - L)} \| \tilde{\xi}^{t,e}_p \|^2 \nonumber \\
  &  + \frac{1}{2\eta^{t}} \Big( \|z_p-z^{t,e}_p\|^2- \|z_p-z^{t,e+1}_p\|^2   \Big) +  \langle \tilde{\xi}^{t,e}_p, z_p - z^{t,e}_p \rangle, \ \forall z_p \in \mathcal{W}_p. \nonumber 
\end{align}

\subsection{Inequality derivation for a fixed iteration $t$}
Summing \eqref{Block2:smooth_inequality_2} over all $e \in [E]$ and dividing the resulting inequalities by $E$, for all $z_p \in \mathcal{W}_p$, we obtain 
\begin{align}
  & \frac{1}{E} \sum_{e=1}^E f_p(z^{t,e+1}_p) - f_p(z_p) -  \langle \lambda^{t+1}_p, \underbrace{\frac{1}{E} \sum_{e=1}^E z^{t,e+1}_p}_{\substack{=z^{t+1}_p \text{ from line 16} \\ \text{in Algorithm \ref{algo:DP-IADMM-Prox}} }}-z_p \rangle \label{Block2:smooth_inequality_3}\\
  \leq & \rho^t  \underbrace{ \frac{1}{E} \sum_{e=1}^E \langle z^{t+1}_p - z_p^{t,e+1}, z_p^{t,e+1} - z_p \rangle}_{\text{``\texttt{E}''}}  \nonumber \\
  & +  \frac{1}{E} \sum_{e=1}^E \Big\{ \frac{\| \tilde{\xi}^{t,e}_p \|^2 }{2(1/\eta^{t} - L)} + \frac{1}{2\eta^{t}} \Big( \|z_p-z^{t,e}_p\|^2- \|z_p-z^{t,e+1}_p\|^2 \Big) + \langle  \tilde{\xi}^{t,e}_p, z_p - z^{t,e}_p \rangle \Big\}. \nonumber 
\end{align}
The ``\texttt{E}'' term from \eqref{Block2:smooth_inequality_3} is nonpositive because
\begin{align}
  & \frac{1}{E}\sum_{e=1}^E \langle z_p^{t+1} - z_p^{t,e+1}, z_p^{t,e+1} - z_p \rangle  \label{non_positive_trick} \\
  = & \frac{1}{E^2}\sum_{e=1}^E \sum_{e'=1}^E \langle z_p^{t,e'+1} - z_p^{t,e+1}, z_p^{t,e+1} - z_p \rangle \nonumber \\
  = & \frac{1}{E^2} \sum_{e=1}^E \sum_{e'=1:e'>e}^E \left( \langle z_p^{t,e'+1} - z_p^{t,e+1}, z_p^{t,e+1} - z_p \rangle + \langle z_p^{t,e+1} - z_p^{t,e'+1}, z_p^{t,e'+1} - z_p \rangle \right) \nonumber \\
  = & \frac{1}{E^2} \sum_{e=1}^E \sum_{e'=1:e'>e}^E \langle z_p^{t,e'+1} - z_p^{t,e+1}, - z_p^{t,e'+1} + z_p^{t,e+1} \rangle  \nonumber \\
  \leq & \frac{1}{E^2} \sum_{e=1}^E \sum_{e'=1:e'>e}^E - \| z_p^{t,e'+1} - z_p^{t,e+1}\|^2 \leq 0. \nonumber
\end{align}
Summing the inequalities resulting from \eqref{Block2:smooth_inequality_3} and \eqref{non_positive_trick} over $p \in [P]$, we have
\begin{align}
  & \sum_{p=1}^P \Big[ \frac{1}{E} \sum_{e=1}^E f_p(z^{t,e+1}_p) - f_p(z_p) -  \langle \lambda^{t+1}_p, z^{t+1}_p - z_p \rangle \Big] \label{Block2:smooth_inequality_4} \\
  \leq  & \sum_{p=1}^P \Big[ \frac{1}{E} \sum_{e=1}^E \Big\{ \frac{\| \tilde{\xi}^{t,e}_p \|^2 }{2(1/\eta^{t} - L)} + \frac{1}{2\eta^{t}} \Big( \|z_p-z^{t,e}_p\|^2- \|z_p-z^{t,e+1}_p\|^2 \Big)\nonumber \\
  &  + \langle  \tilde{\xi}^{t,e}_p, z_p - z^{t,e}_p \rangle \Big\} \Big],  \forall z_1 \in \mathcal{W}_1, \ldots, z_P \in \mathcal{W}_P. \nonumber
\end{align}
Recall that $\mathcal{W} := \cap_{p=1}^P \mathcal{W}_p \neq \emptyset$.
For ease of exposition, we introduce the following notation:
\begin{align}
& 
z := [z_1^{\top}, \ldots, z_P^{\top}]^{\top},  
\ \ \lambda := [\lambda_1^{\top}, \ldots, \lambda_P^{\top}]^{\top}, 
\ \ \tilde{\lambda} := [\tilde{\lambda}_1^{\top}, \ldots, \tilde{\lambda}_P^{\top}]^{\top},  
\label{Notation_smooth}   \\  
& x := \begin{bmatrix}
  w \\ z \\ \lambda
  \end{bmatrix}, \  
  \tilde{x}^t := \begin{bmatrix}
  w^{t+1} \\ z^{t+1} \\ \tilde{\lambda}^t
  \end{bmatrix}, \
  x^* := \begin{bmatrix}
    w^* \\ z^* \\ \lambda
  \end{bmatrix}, \
  A := \begin{bmatrix}
    \mathbb{I}_J \\ \vdots \\ \mathbb{I}_J
  \end{bmatrix}, \
  G := \begin{bmatrix}
      0 & 0 & A^{\top}     \\
      0 & 0 & -\mathbb{I}_{PJ}     \\
      -A & \mathbb{I}_{PJ} & 0     \\
  \end{bmatrix}, \nonumber  \\    
& x^{(T)} := \textstyle \frac{1}{T} \sum_{t=1}^T \tilde{x}^{t}, \ \
w^{(T)} := \textstyle \frac{1}{T} \sum_{t=1}^T w^{t+1}, \ \
z^{(T)} := \textstyle \frac{1}{TE} \sum_{t=1}^T \sum_{e=1}^{E} z^{t,e+1}, \ \
\nonumber \\
&\lambda^{(T)} := \textstyle \frac{1}{T} \sum_{t=1}^T \tilde{\lambda}^{t}, \ \ \textstyle A^{\top} \tilde{\lambda}^t = \sum_{p=1}^P \tilde{\lambda}^t_p , 
\ \  F(z) := \textstyle \sum_{p=1}^P f_p(z_p), \nonumber \\
& \tilde{\xi}^{t,e} := [(\tilde{\xi}^{t,e}_1)^{\top}, \ldots, (\tilde{\xi}^{t,e}_P)^{\top}]^{\top}. \nonumber
\end{align}
Based on the above notation as well as \eqref{Block1:optimality_condition} and \eqref{Block2:smooth_inequality_4}, we derive $\text{LHS}^t(w^*,z^*) \leq \text{RHS}^t(z^*)$ at optimal $w^* \in \mathbb{R}^n, z^*_1 \in \mathcal{W}_1, \ldots z^*_P \in \mathcal{W}_P$, where
\begin{subequations}
\begin{align}
  & \text{LHS}^t(w^*,z^*)  := \nonumber \\
  & \frac{1}{E}\sum_{e=1}^E F(z^{t,e+1}) - F(z^*) - \langle \lambda^{t+1}, z^{t+1}- z^* \rangle + \langle A^{\top} \tilde{\lambda}^t, w^{t+1}-w^* \rangle, \nonumber  \\
  & \text{RHS}^t(z^*)  := \nonumber  \\
  & \frac{1}{E}\sum_{e=1}^E \Big\{ \frac{  \| \tilde{\xi}^{t,e} \|^2}{2(1/\eta^{t}-L)}    + \frac{1}{2\eta^{t}} \big( \|z^* -z^{t,e}\|^2  - \| z^* - z^{t,e+1} \|^2 \big)  + \langle \tilde{\xi}^{t,e}, z^* - z^{t,e} \rangle \Big\}.  \nonumber
\end{align}
\end{subequations}

\subsection{Lower bound} \label{apx-sec:smooth_LB_LHS}
Recall that
\begin{align}
  \lambda^{t+1} = \lambda^t + \rho^t (A w^{t+1} - z^{t+1}), \ \  \tilde{\lambda}^{t} = \lambda^t + \rho^t (A w^{t+1} - z^t).  \label{lambdas}
\end{align}
By utilizing \eqref{lambdas}, we rewrite $\text{LHS}^t(w^*,z^*)$ as follows:
\begin{subequations} 
\begin{align}
& \text{LHS}^t(w^*,z^*)  \label{Smooth_LHS_t_LB_1}  \\
= & \frac{1}{E}\sum_{e=1}^E F(z^{t,e+1}) - F(z^*) \nonumber \\
& + \Biggr\langle
\begin{bmatrix}
w^{t+1} - w^* \\
z^{t+1} - z^* \\
\tilde{\lambda}^t - \lambda
\end{bmatrix}
, \
\begin{bmatrix}
A^{\top} \tilde{\lambda}^t             \\
-\tilde{\lambda}^t \\
- Aw^{t+1} + z^{t+1}
\end{bmatrix}
-
\begin{bmatrix}
0 \\
\rho^t (z^t - z^{t+1})        \\
(\lambda^t- \lambda^{t+1} ) / \rho^t
\end{bmatrix}
\Biggr \rangle.  \nonumber
\end{align}
The third term in \eqref{Smooth_LHS_t_LB_1} can be written as
  \begin{align}
      &\Biggr \langle
      \begin{bmatrix}
          w^{t+1} - w^* \\
          z^{t+1} - z^*  \\
          \tilde{\lambda}^t - \lambda
      \end{bmatrix}
      ,
      \begin{bmatrix}
      A^{\top} \tilde{\lambda}^t \\
      - \tilde{\lambda}^t \\
      -Aw^{t+1} + z^{t+1}
      \end{bmatrix}
      \Biggr \rangle
      = \langle \tilde{x}^{t} - x^*, G\tilde{x}^{t} \rangle  \label{Smooth_LHS_t_LB_2} \\
      =& \underbrace{\langle \tilde{x}^{t} - x^*, G(\tilde{x}^{t} - x^*) \rangle}_{= 0 \text{ as $G$ is skew-symmetric}} + \langle \tilde{x}^{t} - x^* , G x^* \rangle  = \langle \tilde{x}^{t} - x^*, G x^*  \rangle. \nonumber
  \end{align}
Based on \eqref{rule}, the last term in \eqref{Smooth_LHS_t_LB_1} can be written as
    \begin{align}
        & \Biggr\langle
        \begin{bmatrix}
        w^* - w^{t+1} \\
        z^* - z^{t+1}   \\
        \lambda - \tilde{\lambda}^t
        \end{bmatrix}
        , \
        \begin{bmatrix}
        0 \\
        \rho^t (z^t-z^{t+1})  \\
        (\lambda^t- \lambda^{t+1} ) / \rho^t
        \end{bmatrix}
        \Biggr \rangle 
         \nonumber \\
        =  & \textstyle \frac{\rho^t}{2} \big( \| z^* - z^{t+1}\|^2 - \|z^* - z^t\|^2 + \|z^{t+1} - z^t\|^2 \big) \nonumber \\
        &  + \frac{1}{2\rho^t} \big( \| \lambda - \lambda^{t+1} \|^2 - \| \lambda - \lambda^t \|^2 + \underbrace{\|\tilde{\lambda}^{t} - \lambda^t \|^2}_{\geq 0}  - \underbrace{\| \tilde{\lambda}^t - \lambda^{t+1} \|^2}_{=\|\rho^t(z^{t+1}-z^t)\|^2} \big) \nonumber \\
        \geq & \textstyle  \frac{\rho^t}{2} \big( \|z^* - z^{t+1}\|^2  - \|z^* - z^t\|^2 \big) +  \frac{1}{2\rho^t}  \big( \| \lambda - \lambda^{t+1} \|^2  - \| \lambda - \lambda^t \|^2 \big). \nonumber 
    \end{align}  
Therefore, we have
\begin{align}
  & \text{LHS}^t(w^*,z^*)  \label{Smooth_LHS_t_LB_4} \\
  \geq & \frac{1}{E}\sum_{e=1}^E F(z^{t,e+1}) - F(z^*) + \langle \tilde{x}^{t} - x^*, G x^*  \rangle  \nonumber \\
  & + \frac{\rho^t}{2} \big( \|z^* - z^{t+1}\|^2  - \|z^*-z^t\|^2 \big) + \frac{1}{2\rho^t} \big( \| \lambda - \lambda^{t+1} \|^2  - \| \lambda - \lambda^t \|^2 \big).  \nonumber 
\end{align}
\end{subequations}

Let $\text{LHS}(w^*,z^*) :=\frac{1}{T} \sum_{t=1}^T \ \text{LHS}^t(w^*,z^*)$.
We aim to find a lower bound on  $\text{LHS}(w^*,z^*)$.
\begin{align}
& \text{LHS}(w^*,z^*)   \label{LB_LHS_1}\\
\geq & \underbrace{\frac{1}{TE} \sum_{t=1}^T \sum_{e=1}^E F(z^{t,e+1})}_{\geq F(z^{(T)}) \text{ as F is convex}} -  F(z^*) + \underbrace{\langle x^{(T)} - x^* , Gx^* \rangle}_{\text{``\texttt{F}''}} \nonumber \\
& + \frac{1}{T} \Big\{ \underbrace{\sum_{t=1}^T  \frac{ \rho^t}{2}  \big( \|z^* - z^{t+1}\|^2  - \|z^*-z^t\|^2 \big)}_{\text{``\texttt{G}''}} +  \underbrace{\sum_{t=1}^T  \frac{1}{2\rho^t} \big( \| \lambda - \lambda^{t+1} \|^2  - \| \lambda - \lambda^t \|^2 \big)}_{\text{``\texttt{H}''}} \Big\}. \nonumber
\end{align}
The ``\texttt{F}'' term in \eqref{LB_LHS_1} can be written as
\begin{align*}    
& \big\langle x^{(T)} - x^*, Gx \big\rangle   \\
= & \langle A w^{(T)} - z^{(T)} \underbrace{- Aw^* + z^*}_{=0}, \lambda \rangle - \langle \lambda^{(T)} - \lambda, \underbrace{Aw^*-z^*}_{=0} \rangle = \langle \lambda, Aw^{(T)} - z^{(T)} \rangle.
\end{align*}
The ``\texttt{G}'' term in \eqref{LB_LHS_1} can be written as
\begin{align*}    
& \sum_{t=1}^T \frac{\rho^t}{2}  \big( \|z^* - z^{t+1}\|^2  - \|z^*-z^t\|^2 \big)  \nonumber \\
= & - \frac{\rho^1}{2} \|z^*-z^1\|^2 + \sum_{t=2}^T \underbrace{\Big( \frac{\rho^{t-1} - \rho^t}{2} \Big)}_{\leq 0}\|z^* - z^t\|^2 + \underbrace{\frac{\rho^T}{2}\|z^* - z^{T+1}\|^2}_{\geq 0} \nonumber  \\
\geq & - \frac{\rho^1}{2} U_2^2 + \sum_{t=2}^T \Big( \frac{\rho^{t-1} - \rho^t}{2} \Big) U_2^2 = \frac{ - U_2^2 \rho^{T}}{2} \geq \frac{ - U_2^2 \rho^{\text{max}}}{2}. 
\end{align*}
The ``\texttt{H}'' term in \eqref{LB_LHS_1} can be written as
\begin{align*}              
& \sum_{t=1}^T \frac{1}{2\rho^t} \big( \| \lambda - \lambda^{t+1} \|^2  - \| \lambda - \lambda^t \|^2 \big) \nonumber \\
= & - \frac{1}{2\rho^1} \|\lambda-\lambda^1\|^2 + \underbrace{\sum_{t=2}^T \Big( \frac{1}{2 \rho^{t-1}} - \frac{1}{2 \rho^{t}} \Big)\|\lambda - \lambda^t\|^2 + \frac{1}{2\rho^T}\|\lambda - \lambda^{T+1}\|^2}_{\geq 0}  \nonumber \\
\geq & - \frac{1}{2\rho^1} \|\lambda-\lambda^1\|^2. \nonumber
\end{align*}
Therefore, we derive
\begin{align}
  & \text{LHS}(w^*,z^*) \label{Smooth_LHS_LB}\\
  \geq & F(z^{(T)}) - F(z^*) +  \langle \lambda, A w^{(T)} - z^{(T)} \rangle - \frac{1}{T} \Big( \frac{U_2^2  \rho^{\text{max}}}{2} + \frac{1}{2\rho^1} \|\lambda-\lambda^1\|^2 \Big). \nonumber
\end{align}
Since this inequality holds for any $\lambda$, we select $\lambda$ that maximizes the right-hand side of \eqref{Smooth_LHS_LB} subject to a ball centered at zero with the radius $\gamma$:
\begin{subequations}
\label{max_lambdas}  
\begin{align}
\bullet \ & \max_{\lambda: \| \lambda \| \leq \gamma } \ \langle \lambda, Aw^{(T)}-z^{(T)} \rangle = \gamma \| Aw^{(T)}-z^{(T)} \|,  \\  
\bullet \ & \max_{\lambda: \| \lambda \| \leq \gamma } \ \| \lambda - \lambda^1 \|^2 = \| \lambda^1 \|^2 + \max_{\lambda: \| \lambda \| \leq \gamma } \ \{ \| \lambda\|^2 - 2 \langle \lambda, \lambda^1\rangle \} \leq  (\gamma + \|\lambda^1\|)^2.
\end{align}
\end{subequations}
Based on \eqref{Smooth_LHS_LB} and \eqref{max_lambdas}, we derive
\begin{align}
  & \text{LHS}(w^*,z^*) \label{Smooth_LHS_t_LB_3} \\
  \geq & F(z^{(T)}) - F(z^*) + \gamma \| Aw^{(T)}-z^{(T)} \|  -   \frac{U_2^2  \rho^{\text{max}} +  (\gamma + \|\lambda^1\|)^2 / \rho^1 }{2T}. \nonumber
\end{align}

\subsection{Upper bound} \label{apx-sec:RHS_UB}
Let $\text{RHS}(z^*) := \frac{1}{T} \sum_{t=1}^T \ \text{RHS}^t(z^*)$. We aim to find an upper bound on $\text{RHS}(z^*)$.
\begin{align*}
\text{RHS}(z^*) = & \frac{1}{TE} \Big[ \sum_{t=1}^T \sum_{e=1}^E  \Big\{  \frac{ \| \tilde{\xi}^{t,e} \|^2}{2(1/\eta^{t}-L)} + \langle \tilde{\xi}^{t,e}, z^* - z^{t,e} \rangle \Big\} \nonumber \\
&+ \underbrace{ \sum_{t=1}^T \sum_{e=1}^E \frac{1}{2\eta^{t}} ( \|z^* -z^{t,e} \|^2  - \|z^*-z^{t,e+1}\|^2) }_{\text{``\texttt{I}''}} \Big]. 
\end{align*}
The ``\texttt{I}'' term from the above can be written as
\begin{align}
  & \sum_{t=1}^T \sum_{e=1}^E \frac{1}{2\eta^{t}}\big(\|z^*-z^{t,e}\|^2-\|z^*-z^{t,e+1}\|^2 \big) = \sum_{t=1}^T  \frac{1}{2\eta^{t}}\big(\|z^*-z^{t,1}\|^2-\|z^*- \underbrace{z^{t,E+1}}_{= z^{t+1,1}}\|^2 
  \big) \nonumber \\
  & = \frac{1}{2\eta^1} \| z^* - z^{1,1} \|^2 + \sum_{t=2}^T \underbrace{\Big( \frac{1}{2\eta^{t}} - \frac{1}{2\eta^{t-1}} \Big)}_{\geq 0} \|z^*-z^{t,1}\|^2 \underbrace{ - \frac{1}{2\eta^{T}} \| z^*- z^{T,E+1}  \|^2}_{\leq 0}  \leq \frac{U_2^2}{2\eta^{T}}. \nonumber
\end{align}
Therefore, we have
\begin{align}
\text{RHS}(z^*) \leq \frac{U_2^2}{2 T E \eta^{T}} + \frac{1}{TE} \sum_{t=1}^T \sum_{e=1}^E \Big\{ \frac{ \| \tilde{\xi}^{t,e} \|^2 }{2(1/\eta^{t} - L)} +  \langle \tilde{\xi}^{t,e}, z^* - z^{t,e} \rangle \Big\}. \label{Smooth_RHS_t_LB_1}
\end{align}  

\subsection{Taking expectation}
\noindent
By taking expectation on the inequality derived from \eqref{Smooth_LHS_t_LB_3} and \eqref{Smooth_RHS_t_LB_1}, we have
\begin{align*}
  & \mathbb{E} \Big[ F(z^{(T)}) - F(z^*) + \gamma \| Aw^{(T)} - z^{(T)} \| \Big]  \leq \frac{U_2^2  \rho^{\text{max}} +  (\gamma + \|\lambda^1\|)^2 / \rho^1 }{2T} \nonumber \\
  & + \frac{U_2^2}{2 T E \eta^{T}} + \frac{1}{TE} \sum_{t=1}^T \sum_{e=1}^E \Big\{ \frac{ 1 }{2(1/\eta^{t} - L)} \sum_{p=1}^P  \mathbb{E}[\| \tilde{\xi}^{t,e}_p \|^2]  +  \langle \underbrace{\mathbb{E}[\tilde{\xi}^{t,e}]}_{=0}, z^* - z^{t,e} \rangle \Big\}.
\end{align*}
Note that the second moment of each element of $\tilde{\xi}^{t,e}_p$ is the same as its variance due to the zero mean. Thus we have
\begin{align}
  & \mathbb{E}[ (\tilde{\xi}^{t,e}_{pi} )^2] = 
  \begin{cases}
    2 ( \bar{\Delta}^{t,e}_{p,1} / \bar{\epsilon} )^2  & \text{ for the Laplace mechanism }\\
    2 \ln(1.25/\bar{\delta}) ( \bar{\Delta}^{t,e}_{p,2} / \bar{\epsilon} )^2  & \text{ for the Gaussian mechanism }\\
  \end{cases}
  , \ \forall i \in [n].
  \label{U_eps} 
\end{align}
Therefore, 
\begin{align}
\mathbb{E}[ \| \tilde{\xi}^{t,e}_p \|^2] =  \sum_{i=1}^n \mathbb{E}[ (\tilde{\xi}^{t,e}_{pi} )^2] \leq n U_3 / \bar{\epsilon}^2
\end{align}
Also note that
\begin{align*}  
  \bullet \ & \frac{U _2^2}{2 T E \eta^{T}} = \frac{U_2^2 L}{2TE} + \frac{U_2^2 }{2 \bar{\epsilon} E \sqrt{T}},  \\    
  \bullet \ & \sum_{t=1}^T \sum_{e=1}^E \frac{1}{2(1/\eta^{t} - L)} = E \bar{\epsilon} \sum_{t=1}^T \frac{1}{2 \sqrt{t}} \leq E \bar{\epsilon} \sum_{t=1}^T \frac{1}{\sqrt{t} + \sqrt{t-1}} \nonumber \\
  & =  E \bar{\epsilon} \sum_{t=1}^T (\sqrt{t} - \sqrt{t-1}) = E \bar{\epsilon} \sqrt{T}.  
\end{align*}
Therefore, we have
\begin{align*}
  & \mathbb{E} \Big[ F(z^{(T)}) - F(z^*) + \gamma \| Aw^{(T)} - z^{(T)} \| \Big]  \nonumber \\
  & \leq \frac{U_2^2  (\rho^{\text{max}} +L/E) +  (\gamma + \|\lambda^1\|)^2 / \rho^1 }{2T}   + \frac{ nP U_3 + U_2^2/(2E) }{ \bar{\epsilon} \sqrt{T} }.
\end{align*}
This completes the proof.

 
\section{Proof of Theorem \ref{thm:nonsmooth}} \label{apx-thm:nonsmooth}
\noindent
The proof in this section is similar to that in Appendix \ref{apx-thm:smooth} except that
\begin{enumerate}
\item the $L$-smoothness of $f_p$ can no longer be applied to the ``\texttt{D}'' term in \eqref{Block2:smooth_inequality_1} when deriving an upper bound of the term in a nonsmooth setting and
\item the definition of $(\tilde{x}^t, z^{(T)})$ is different from that in \eqref{Notation_smooth} of Appendix \ref{apx-thm:smooth}.
\end{enumerate}

\subsection{Inequality derivation for a fixed iteration $t$ and $e$}
Applying Young's inequality on the ``\texttt{D}'' term in \eqref{Block2:smooth_inequality_1} yields
\begin{align}
& f_p(z^{t,e}_p) - f_p(z_p) -  \langle \lambda^{t+1}_p, z^{t,e+1}_p-z_p \rangle   \label{Block2:nonsmooth_inequality_1} \\
\leq &  \rho^t \langle z^{t+1}_p - z_p^{t,e+1}, z_p^{t,e+1} - z_p \rangle  +\frac{\eta^{t}}{2} \| f'_p(z^{t,e}_p) + \tilde{\xi}^{t,e}_p\|^2   \nonumber \\
& + \frac{1}{2\eta^{t}} \Big\{ \|z_p-z^{t,e}_p\|^2- \|z_p-z^{t,e+1}_p\|^2 \Big\}  + \langle  \tilde{\xi}^{t,e}_p, z_p - z^{t,e}_p \rangle, \ \forall z_p \in \mathcal{W}_p. \nonumber
\end{align}

\subsection{Inequality derivation for a fixed iteration $t$} \label{apx-sec:nonsmooth_ineq}
Summing \eqref{Block2:nonsmooth_inequality_1} over all $e \in [E]$ and dividing the resulting inequalities by $E$, we get 
\begin{align}
  & \frac{1}{E} \sum_{e=1}^E f_p(z^{t,e}_p) - f_p(z_p) -  \langle \lambda^{t+1}_p, \underbrace{\frac{1}{E} \sum_{e=1}^E z^{t,e+1}_p}_{=z^{t+1}_p}-z_p \rangle  \label{Block2:nonsmooth_inequality_2} \\
  \leq &  \rho^t  \underbrace{ \frac{1}{E} \sum_{e=1}^E \langle z_p^{t+1} - z_p^{t,e+1}, z_p^{t,e+1} - z_p \rangle}_{\leq 0 \text{ from \eqref{non_positive_trick}} } + \langle  \tilde{\xi}^{t,e}_p, z_p - z^{t,e}_p \rangle  \nonumber \\
  & + \frac{1}{E} \sum_{e=1}^E \Big\{ \frac{\eta^{t}}{2} \| f'_p(z^{t,e}_p) + \tilde{\xi}^{t,e}_p\|^2  + \frac{1}{2\eta^{t}} \Big( \|z_p-z^{t,e}_p\|^2- \|z_p-z^{t,e+1}_p\|^2 \Big) \Big\}. \nonumber
\end{align}
Summing the inequalities \eqref{Block2:nonsmooth_inequality_2} over $p \in [P]$, we have 
\begin{align}
  &\sum_{p=1}^P \Big[ \frac{1}{E} \sum_{e=1}^E f_p(z^{t,e}_p) - f_p(z_p) -  \langle \lambda^{t+1}_p,  z^{t+1}_p -z_p \rangle \Big]  \label{Block2:nonsmooth_inequality_3}\\
  \leq  & \sum_{p=1}^P \Big[ \frac{1}{E} \sum_{e=1}^E \Big\{ \frac{\eta^{t}}{2} \| f'_p(z^{t,e}_p) + \tilde{\xi}^{t,e}_p\|^2  + \frac{1}{2\eta^{t}} \Big( \|z_p-z^{t,e}_p\|^2- \|z_p-z^{t,e+1}_p\|^2 \Big) \nonumber \\
  & + \langle  \tilde{\xi}^{t,e}_p, z_p - z^{t,e}_p \rangle \Big\} \Big].  \nonumber
\end{align}

\noindent
For ease of exposition, we introduce $z, \lambda, \tilde{\lambda}, x, x^*, A, G, x^{(T)}, w^{(T)}, \lambda^{(T)}, A^{\top} \tilde{\lambda}^t, F(z), \tilde{\xi}^{t,e}$ defined in \eqref{Notation_smooth} with modifications of the following notation:
\begin{subequations}
\label{nonsmooth_notations}
\begin{align}
& \tilde{x}^t := \begin{bmatrix}
  w^{t+1} \\ z^{t} \\ \tilde{\lambda}^t
  \end{bmatrix}, \ \
z^{(T)} := \textstyle \frac{1}{TE} \sum_{t=1}^T \sum_{e=1}^{E} z^{t,e}.
\end{align}
We also define
\begin{align}
f'(z):=[ f'_1(z_1)^{\top}, \ldots, f'_P(z_P)^{\top} ]^{\top}.
\end{align}
\end{subequations}
Based on this notation as well as \eqref{Block1:optimality_condition} and \eqref{Block2:nonsmooth_inequality_3}, we derive $\text{LHS}^t(w^*,z^*) \leq \text{RHS}^t(z^*)$ at optimal $w^* \in \mathbb{R}^n, z^*_1 \in \mathcal{W}_1, \ldots z^*_P \in \mathcal{W}_P$, where
\begin{align*}
  & \text{LHS}^t(w^*,z^*) :=  \frac{1}{E}\sum_{e=1}^E F(z^{t,e}) - F(z^*) - \langle \lambda^{t+1}, z^{t+1}- z^* \rangle + \langle A^{\top} \tilde{\lambda}^t, w^{t+1}-w^* \rangle, \\
  & \text{RHS}^t(z^*) :=  \\
  & \frac{1}{E}\sum_{e=1}^E \Big\{ \frac{\eta^{t}}{2} \| f'(z^{t,e}) + \tilde{\xi}^{t,e} \|^2   + \frac{1}{2\eta^{t}} \big( \|z^* -z^{t,e}\|^2  - \| z^* - z^{t,e+1} \|^2 \big)  + \langle \tilde{\xi}^{t,e}, z^* - z^{t,e} \rangle \Big\}.
\end{align*}
 
\subsection{Lower bound}
\label{apx:nonsmooth-B}
We aim to find a lower bound on $\text{LHS}^t(w^*,z^*)$.
By following the steps in Appendix \ref{apx-sec:smooth_LB_LHS}, one can derive inequalities similar to \eqref{Smooth_LHS_t_LB_4}, as follows:
\begin{align}
  & \text{LHS}^t(w^*,z^*) \geq \frac{1}{E}\sum_{e=1}^E F(z^{t,e}) - F(z^*) + \langle \tilde{x}^{t} - x^*, G x^*  \rangle \underbrace{- \langle \lambda, z^{t+1} - z^t \rangle}_{\text{``\texttt{J}''}}  \label{Nonsmooth_LHS_t_LB}\\
  &  + \frac{\rho^t}{2} \big( \|z^* - z^{t+1}\|^2  - \|z^*-z^t\|^2 \big) + \frac{1}{2\rho^t} \big( \| \lambda - \lambda^{t+1} \|^2  - \| \lambda - \lambda^t \|^2 \big).  \nonumber  
\end{align}
Note that the ``\texttt{J}'' term in \eqref{Nonsmooth_LHS_t_LB} does not exist in \eqref{Smooth_LHS_t_LB_4} because the definition of $\tilde{x}^t$ in \eqref{nonsmooth_notations} is different from that in \eqref{Notation_smooth}.

We aim to find a lower bound on $\text{LHS}(w^*,z^*) :=\frac{1}{T} \sum_{t=1}^T \ \text{LHS}^t(w^*,z^*)$.
Summing \eqref{Nonsmooth_LHS_t_LB} over $t \in [T]$ and dividing the resulting inequality by $T$, we have
\begin{align}
& \text{LHS}(w^*,z^*) \label{nonsmooth_LB_LHS_1} \\
\geq & \underbrace{\frac{1}{TE}  \sum_{t=1}^T \sum_{e=1}^E F(z^{t,e}) }_{\geq F(z^{(T)}) \text{ as F is convex} }-  F(z^*) + \underbrace{\langle x^{(T)} - x^* , Gx^* \rangle}_{= \text{``\texttt{F}'' term in \eqref{LB_LHS_1}} } - \underbrace{\frac{1}{T} \langle \lambda, z^{T+1} - z^1 \rangle}_{\text{``\texttt{K}''}} \nonumber \\
& + \frac{1}{T} \Big\{ \underbrace{\sum_{t=1}^T  \frac{ \rho^t}{2}  \big( \|z^* - z^{t+1}\|^2  - \|z^*-z^t\|^2 \big)}_{= \text{``\texttt{G}'' term in \eqref{LB_LHS_1}}} +  \underbrace{\sum_{t=1}^T  \frac{1}{2\rho^t} \big( \| \lambda - \lambda^{t+1} \|^2  - \| \lambda - \lambda^t \|^2 \big)}_{= \text{``\texttt{H}'' term in \eqref{LB_LHS_1}}} \Big\}. \nonumber
\end{align}
The ``\texttt{K}'' term in \eqref{nonsmooth_LB_LHS_1} can be written as
\begin{align*}  
- \frac{1}{T} \langle \lambda, z^{T+1} - z^1 \rangle \geq - \frac{1}{T}  \| \lambda \| \| z^{T+1} - z^1 \| \geq - \| \lambda \| U_2.
\end{align*}
Therefore, we have
\begin{align}
  & \text{LHS}(w^*,z^*) \nonumber \\
  \geq & F(z^{(T)}) - F(z^*) +  \langle \lambda, A w^{(T)} - z^{(T)} \rangle - \frac{1}{T} \Big(\| \lambda \| U_2 + \frac{U_2^2  \rho^{\text{max}}}{2} + \frac{1}{2\rho^1} \|\lambda-\lambda^1\|^2 \Big). \nonumber
\end{align}
Based on \eqref{max_lambdas} and $\max_{\lambda: \| \lambda \| \leq \gamma } \ \| \lambda \| U_2 \leq \gamma U_2$, we have
\begin{align}
  & \text{LHS}(w^*,z^*) \label{Nonsmooth_LHS_LB_7} \\
  \geq & F(z^{(T)}) - F(z^*) + \gamma \| Aw^{(T)}-z^{(T)} \|  -   \frac{U_2^2 \rho^{\text{max}} + (\gamma + \|\lambda^1\|)^2 / \rho^1 + 2 \gamma U_2 }{2T}. \nonumber
\end{align}

\subsection{Upper bound}
Let $\text{RHS}(z^*) := \frac{1}{T} \sum_{t=1}^T \ \text{RHS}^t(z^*)$.
We aim to find an upper bound on $\text{RHS}(z^*)$.
By following the steps in Appendix \ref{apx-sec:RHS_UB}, we obtain
\begin{align}
\text{RHS}(z^*) \leq \frac{U_2^2}{2 TE \eta^{T}} + \frac{1}{TE} \sum_{t=1}^T \sum_{e=1}^E \Big\{ \frac{\eta^{t}}{2} \| f'(z^{t,e}) + \tilde{\xi}^{t,e} \|^2 +  \langle \tilde{\xi}^{t,e}, z^* - z^{t,e} \rangle \Big\}. \label{Nonsmooth_RHS_UB_2}
\end{align}  

\subsection{Taking expectation}
\noindent
By taking expectation on the inequality derived from \eqref{Nonsmooth_LHS_LB_7} and \eqref{Nonsmooth_RHS_UB_2}, we have
\begin{align*}
  & \mathbb{E} \Big[ F(z^{(T)}) - F(z^*) + \gamma \| Aw^{(T)} - z^{(T)} \| \Big]  \leq \frac{U_2^2 \rho^{\text{max}} + (\gamma + \|\lambda^1\|)^2 / \rho^1 + 2\gamma U_2}{2T} \nonumber \\
  & + \frac{U_2^2}{2 TE \eta^{T}} + \frac{1}{TE} \sum_{t=1}^T \sum_{e=1}^E \Big\{ \frac{\eta^{t}}{2} \sum_{p=1}^P \underbrace{\mathbb{E}[ \| f'_p(z_p^{t,e}) + \tilde{\xi}^{t,e}_p \|^2 ]}_{\leq U_1^2+ n U_3/\bar{\epsilon}^2} +  \langle \underbrace{\mathbb{E}[\tilde{\xi}^{t,e}]}_{=0}, z^* - z^{t,e} \rangle \Big\}.
\end{align*}
We note that
\begin{align*}
 \bullet \ & \frac{U_2^2}{2 TE \eta^{T}} = \frac{U_2^2 /(2E)}{ \sqrt{T}},  \\
 \bullet \ & \sum_{t=1}^T \sum_{e=1}^E \frac{\eta^{t}}{2}   = E \sum_{t=1}^T \frac{1}{2 \sqrt{t}} \leq E \sum_{t=1}^T \frac{1}{\sqrt{t} + \sqrt{t-1}} = E \sum_{t=1}^T   (\sqrt{t} - \sqrt{t-1}) =  E  \sqrt{T}.
\end{align*}
Therefore, we have
\begin{align*}
  & \mathbb{E} \Big[ F(z^{(T)}) - F(z^*) + \gamma \| Aw^{(T)} - z^{(T)} \| \Big]  \nonumber \\
  \leq & \frac{U_2^2 \rho^{\text{max}} + (\gamma + \|\lambda^1\|)^2 / \rho^1 +  2\gamma U_2}{2T}   + \frac{ nPU_3 /\bar{\epsilon}^2 + PU_1^2 + U_2^2 / (2E) }{ \sqrt{T} }.
\end{align*}
This completes the proof. 


\section{Proof of Theorem \ref{thm:strong}} \label{apx-thm:strong}
\noindent
The proof in this section is similar to that in Appendix \ref{apx-thm:nonsmooth} except that
\begin{enumerate}
  \item the $\alpha$-strong convexity of $f_p$ is utilized to tighten the right-hand side of inequality \eqref{Block2:nonsmooth_inequality_1} and
  \item the definition of $x^{(T)}, w^{(T)}, z^{(T)}, \lambda^{(T)}$ is modified to the following:
    \begin{align}  
    & x^{(T)} := \textstyle \frac{2}{T(T+1)} \sum_{t=1}^T t \tilde{x}^{t}, \ \
    w^{(T)} := \textstyle \frac{2}{T(T+1)} \sum_{t=1}^T t w^{t+1},  \label{Notation_strong}  \\
    & z^{(T)} := \textstyle \frac{2}{T(T+1)} \sum_{t=1}^T t  (\frac{1}{E} \sum_{e=1}^E z^{t,e}), \ \
    \lambda^{(T)} := \textstyle \frac{2}{T(T+1)} \sum_{t=1}^T t \tilde{\lambda}^{t}.  \nonumber
    \end{align}
\end{enumerate}

\subsection{Inequality derivation for a fixed iteration $t$ and $e$}
For a given $p \in [P]$, it follows from the $\alpha$-strong convexity of the function $f_p$ that
\begin{align}
f_p(z^{t,e}_p) - f_p(z_p)  \leq \langle f'_p(z^{t,e}_p), z^{t,e}_p - z_p \rangle - \frac{\alpha}{2}\| z_p - z^{t,e}_p\|^2. \label{Block2:strong_inequality_1}
\end{align} 
By utilizing \eqref{Block2:strong_inequality_1} and \eqref{Block2:nonsmooth_inequality_1}, we obtain
\begin{align}
& f_p(z^{t,e}_p) - f_p(z_p) - \langle  \lambda^{t+1}_p, z^{t,e+1}_p-z_p \rangle  \label{Block2:strong_inequality_2}  \\
\leq & \rho^t \langle z_p^{t+1} - z_p^{t,e+1}, z_p^{t,e+1} - z_p \rangle  - \frac{\alpha}{2}\| z_p - z^{t,e}_p\|^2  + \langle  \tilde{\xi}^{t,e}_p, z_p - z^{t,e}_p \rangle \nonumber \\
& + \frac{\eta^{t}}{2} \| f'_p(z^{t,e}_p) + \tilde{\xi}^{t,e}_p\|^2  + \frac{1}{2\eta^{t}} \Big( \|z_p-z^{t,e}_p\|^2- \|z_p-z^{t,e+1}_p\|^2 \Big ), \ \forall z_p \in \mathcal{W}_p.\nonumber
\end{align}
Note that compared with \eqref{Block2:nonsmooth_inequality_1}, the inequalities \eqref{Block2:strong_inequality_2} have an additional term $-(\alpha/2)\|z_p-z_p^{t,e}\|^2$. 

\subsection{Inequality derivation for a fixed iteration $t$}
Following the steps to derive \eqref{Block2:nonsmooth_inequality_3} in Appendix \ref{apx-sec:nonsmooth_ineq}, we derive the following from \eqref{Block2:strong_inequality_2}:
\begin{align}
  & \sum_{p=1}^P \Big[ \frac{1}{E} \sum_{e=1}^E f_p(z^{t,e}_p) - f_p(z_p) -  \langle \lambda^{t+1}_p,  z^{t+1}_p -z_p \rangle \Big] \label{Block2:strong_inequality_3} \\
  \leq & \sum_{p=1}^P \Big[ \frac{1}{E} \sum_{e=1}^E \Big\{ \frac{\eta^{t}}{2} \| f'_p(z^{t,e}_p) + \tilde{\xi}^{t,e}_p\|^2  \nonumber \\
  &  + \big(\frac{1}{2\eta^{t}} - \frac{\alpha}{2} \big) \|z_p-z^{t,e}_p\|^2 - \frac{1}{2\eta^{t}} \|z_p-z^{t,e+1}_p\|^2  + \langle  \tilde{\xi}^{t,e}_p, z_p - z^{t,e}_p \rangle \Big\} \Big]. \nonumber 
\end{align}
For ease of exposition, we introduce $z, \lambda, \tilde{\lambda}, x, x^*, A, G, A^{\top} \tilde{\lambda}^t, F(z), \tilde{\xi}^{t,e}$ as defined in \eqref{Notation_smooth}, $\tilde{x}^t, f'(z)$ as defined in \eqref{nonsmooth_notations}, and the definition \eqref{Notation_strong}.
Based on this notation as well as \eqref{Block1:optimality_condition} and \eqref{Block2:strong_inequality_3}, we derive $\text{LHS}^t(w^*,z^*) \leq \text{RHS}^t(z^*)$ at optimal $w^* \in \mathbb{R}^n, z^*_1 \in \mathcal{W}_1, \ldots z^*_P \in \mathcal{W}_P$, where
\begin{subequations}
  \begin{align*}
    \text{LHS}^t(w^*,z^*) :=&   \frac{1}{E} \sum_{e=1}^E F(z^{t,e}) - F(z^*) - \langle \lambda^{t+1}, z^{t+1}- z^* \rangle + \langle A^{\top} \tilde{\lambda}^t, w^{t+1}-w^* \rangle,   \\
    \text{RHS}^t(z^*) := &  \frac{1}{E} \sum_{e=1}^E \Big\{  \frac{\eta^{t}}{2} \| f'(z^{t,e}) + \tilde{\xi}^{t,e} \|^2   \\
    &  + \big(\frac{1}{2\eta^{t}} - \frac{\alpha}{2})  \|z^* -z^{t,e}\|^2  - \frac{1}{2\eta^{t}} \| z^* - z^{t,e+1} \|^2  + \langle \tilde{\xi}^{t,e}, z^* - z^{t,e} \rangle \Big\}.  
  \end{align*}
\end{subequations}

\subsection{Lower bound}
Let $\text{LHS}(w^*,z^*) :=\frac{2}{T(T+1)} \sum_{t=1}^T t \ \text{LHS}^t(w^*,z^*)$.
We aim to find a lower bound on $\text{LHS}(w^*,z^*)$:
\begin{align}
\text{LHS}(w^*,z^*)  \geq&\underbrace{\frac{2}{T(T+1)} \sum_{t=1}^T t \big( \frac{1}{E}\sum_{e=1}^E F(z^{t,e}) \big)}_{\geq F(z^{(T)}) \text{ as F is convex}} - F(z^*) +  \underbrace{\langle x^{(T)} - x^*, Gx^* \rangle}_{= \text{``\texttt{F}'' term in \eqref{LB_LHS_1}} }   \label{strong:LHS_LB} \\
& - \frac{2}{T(T+1)} \underbrace{\langle \lambda, \sum_{t=1}^T t (z^{t+1} - z^t ) \rangle}_{\text{``\texttt{L}''}} \nonumber \\
& + \frac{2}{T(T+1)}  \Big\{  \underbrace{\sum_{t=1}^T \frac{t \rho^t}{2}  \big( \|z^* - z^{t+1}\|^2  - \|z^* - z^t\|^2 \big)}_{\text{``\texttt{M}''}} \nonumber \\
& + \underbrace{\sum_{t=1}^T \frac{t}{2\rho^t} \big( \| \lambda - \lambda^{t+1} \|^2  - \| \lambda - \lambda^t \|^2 \big)}_{\text{``\texttt{N}''}} \Big\}. \nonumber
\end{align}
The ``\texttt{L}'' term in \eqref{strong:LHS_LB} can be written as
\begin{align*}
& \langle \lambda, \sum_{t=1}^T t (z^{t+1}-z^t) \rangle \\
= & \langle \lambda,  \sum_{t=1}^T (z^{T+1}-z^t) \rangle \leq \sum_{t=1}^T \| \lambda \| \| z^{T+1} - z^t \| \leq \sum_{t=1}^T \| \lambda \| U_2 =  T U_2 \| \lambda \|.  
\end{align*}
The ``\texttt{M}'' term in \eqref{strong:LHS_LB} can be written as
\begin{align*}
& \sum_{t=1}^T \frac{t \rho^t}{2} \big( \| z^* - z^{t+1} \|^2 - \|z^* - z^t\|^2 \big)   \\
=&  - \frac{\rho^1}{2} \|z^*-z^1\|^2 + \sum_{t=2}^T \Big( \underbrace{\frac{(t-1)\rho^{t-1} - t\rho^t}{2}}_{\leq 0 \text{ as } \rho^t \geq \rho^{t-1}} \Big)\|z^*-z^t \|^2 + \frac{T\rho^T}{2}\|z^*-z^{T+1} \|^2    \\
\geq&  - \frac{\rho^1 U_2^2}{2}  + \sum_{t=2}^T \Big( \frac{(t-1)\rho^{t-1} - t\rho^t}{2} \Big) U_2^2 = \frac{ - U_2^2 T \rho^{T}}{2} \geq \frac{ - U_2^2 T \rho^{\text{max}}}{2}.  
\end{align*}
The ``\texttt{N}'' term in \eqref{strong:LHS_LB} can be written as
\begin{align*}      
& \sum_{t=1}^T \frac{t}{2\rho^t} \big( \| \lambda - \lambda^{t+1}\|^2  - \| \lambda - \lambda^t \|^2 \big) \geq - \frac{1}{2\rho^1} \| \lambda - \lambda^1 \|^2 + \sum_{t=2}^T \Big( \frac{t-1}{2 \rho^{t-1}} - \frac{t}{2 \rho^{t}} \Big)\| \lambda - \lambda^t \|^2. 
\end{align*}
Therefore, we have
\begin{align}
  & \text{LHS}(w^*,z^*)  \label{strong:LHS_LB_1} \\
  \geq &  F(z^{(T)}) - F(z^*) +  \langle \lambda, A w^{(T)} - z^{(T)} \rangle - \frac{2 U_2 \| \lambda \|}{T+1} - \frac{U_2^2 \rho^{\text{max}}}{T+1} \nonumber \\
  &   + \frac{2}{T(T+1)} \Big( - \frac{1}{2\rho^1} \| \lambda - \lambda^1 \|^2 + \sum_{t=2}^T \underbrace{\Big( \frac{t-1}{2 \rho^{t-1}} - \frac{t}{2 \rho^{t}} \Big)}_{ \leq 0 \text{ by Assumption \ref{assump:convergence_stronglyconvex}-(B)}} \| \lambda - \lambda^t \|^2 \Big). \nonumber  
\end{align}
In addition to \eqref{max_lambdas}, by Assumption \ref{assump:convergence_stronglyconvex}-(i), we have
\begin{align*}
\bullet \ & \max_{\lambda: \| \lambda \| \leq \gamma } \ \| \lambda - \lambda^t \|^2 = \| \lambda^t \|^2 + \max_{\lambda: \| \lambda \| \leq \gamma } \ \left\{ \| \lambda\|^2 - 2 \langle \lambda, \lambda^t\rangle \right\} \leq  4 \gamma^2.
\end{align*}
By utilizing this to derive a lower bound of the last term in \eqref{strong:LHS_LB_1}, we have
\begin{align*}
- \frac{1}{2\rho^1} \| \lambda - \lambda^1 \|^2 + \sum_{t=2}^T \Big( \frac{t-1}{2 \rho^{t-1}} - \frac{t}{2 \rho^{t}} \Big) \| \lambda - \lambda^t \|^2 \geq 
- \frac{T}{2\rho^T} 4\gamma^2  \geq - \frac{T}{2\rho^1} 4\gamma^2 .
\end{align*}
Therefore, we have
\begin{align}
  \text{LHS}(w^*,z^*) \geq F(z^{(T)}) - F(z^*) + \gamma \| Aw^{(T)}-z^{(T)} \| - \frac{2 U_2 \gamma + U_2^2 \rho^{\text{max}} + 4 \gamma^2 / \rho^1 }{T+1}. \label{Strong_LHS_t_LB_1}
\end{align}

\subsection{Upper bound}
Let $\text{RHS}(z^*) := \frac{2}{T(T+1)} \sum_{t=1}^T t \  \text{RHS}^t(z^*)$.
We aim to find an upper bound on $\text{RHS}(z^*)$:
\begin{align*}
\text{RHS}(z^*) = & \frac{2}{T(T+1)} \sum_{t=1}^T t \Big[  \frac{1}{E} \sum_{e=1}^E \Big\{  \frac{\eta^{t}}{2} \| f'(z^{t,e}) + \tilde{\xi}^{t,e} \|^2  \\
& + \big(\frac{1}{2\eta^{t}} - \frac{\alpha}{2})  \|z^* -z^{t,e}\|^2  - \frac{1}{2\eta^{t}} \| z^* - z^{t,e+1} \|^2  + \langle \tilde{\xi}^{t,e}, z^* - z^{t,e} \rangle \Big\} \Big]. 
\end{align*}
Note that
\begin{align*}
\bullet \ & \eta^{t} = 2/(\alpha(t+2)), \\
\bullet \ &  \sum_{t=1}^T \sum_{e=1}^E t \frac{\eta^{t}}{2} \| f'(z^{t,e}) + \tilde{\xi}^{t,e} \|^2 = \sum_{t=1}^T \sum_{e=1}^E  \frac{t}{\alpha(t+2)}  \| f'(z^{t,e}) + \tilde{\xi}^{t,e} \|^2, \\
\bullet \ & \sum_{t=1}^T \sum_{e=1}^E t \Big\{ (\frac{1}{2\eta^{t}} - \frac{\alpha}{2} ) \|z^* -z^{t,e} \|^2  - \frac{1}{2\eta^{t}} \|z^*-z^{t,e+1}\|^2 \Big\} \nonumber \\
 = & \frac{\alpha}{4} \sum_{t=1}^T \sum_{e=1}^E  \Big\{ t^2 \|z^* -z^{t,e} \|^2  - (t^2+2t) \|z^*-z^{t,e+1}\|^2 \Big\} \nonumber \\
= &\frac{\alpha}{4} \sum_{t=1}^T \Big\{  t^2   \Big( \|z^* -z^{t,1} \|^2  - \|z^*-z^{t,E+1}\|^2 \Big) -2t  \|z^*-z^{t,E+1}\|^2 \nonumber \\
& - \underbrace{2t\sum_{e=1}^{E-1} \|z^*-z^{t,e+1}\|^2}_{\leq 0}   \Big\} \nonumber \\
\leq &\frac{\alpha}{4} \sum_{t=1}^T \Big\{  t^2  \|z^* -z^{t,1} \|^2  - t(t+2) \|z^*-z^{t,E+1}\|^2 \Big\} \nonumber \\
= & \frac{\alpha}{4} \Big\{  \|z^* -z^{1,1} \|^2  + \sum_{t=2}^T t^2  \|z^* - \underbrace{z^{t,1}}_{ = z^{t-1,E+1}} \|^2  \nonumber \\
& - \underbrace{\sum_{t=1}^{T-1} t(t+2) \|z^*-z^{t,E+1}\|^2}_{=\sum_{t=2}^{T} (t^2-1) \|z^*-z^{t-1,E+1}\|^2} \underbrace{- T(T+2) \|z^*-z^{T,E+1}\|^2}_{\leq 0} \Big\} \nonumber \\
\leq & \frac{\alpha}{4} \Big\{  \|z^* -z^{1,1} \|^2  + \sum_{t=2}^T  \|z^* - z^{t-1,E+1} \|^2  \Big\}   \leq \frac{\alpha}{4} T U_2^2.
\end{align*}  
Therefore, we have
\begin{align}
& \text{RHS}(z^*) \label{Strong_RHS_t_LB_1} \\
\leq &  \frac{\alpha U_2^2/E}{2(T+1)} + \frac{2}{ET(T+1)} \sum_{t=1}^T \sum_{e=1}^E \Big\{ \frac{t}{\alpha(t+2)} \sum_{p=1}^P \| f'_p(z^t_p) + \tilde{\xi}^{t,e}_p \|^2 +   t \langle \tilde{\xi}^{t,e}, z^* - z^{t,e} \rangle \Big\}. \nonumber 
\end{align}  

\subsection{Taking expectation}
By taking expectation on the inequality derived from \eqref{Strong_LHS_t_LB_1} and \eqref{Strong_RHS_t_LB_1}, we have
\begin{align*}
  & \mathbb{E} \Big[ F(z^{(T)}) - F(z^*) + \gamma \| Aw^{(T)} - z^{(T)} \| \Big]  \\
  \leq  & \frac{2 U_2 \gamma + U_2^2 \rho^{\text{max}} + 4 \gamma^2 / \rho^1 + \alpha U_2^2/(2E) }{T+1} \\
  & + \frac{2}{ET(T+1)} \sum_{t=1}^T \sum_{e=1}^E \Big\{ \frac{t}{\alpha(t+2)} \sum_{p=1}^P \underbrace{\mathbb{E}\big[ \| f'_p(z^{t,e}_p) + \tilde{\xi}^{t,e}_p \|^2}_{\leq U_1^2+ nU_3/\bar{\epsilon}^2 } \big] + t  \langle \underbrace{\mathbb{E}[ \tilde{\xi}^{t,e}]}_{=0}, z^* - z^t \rangle \big] \\
\leq &  \frac{2 U_2 \gamma + U_2^2 \rho^{\text{max}} + 4 \gamma^2 / \rho^1 + \alpha U_2^2/(2E) }{T+1} \\
& +  \frac{2}{ET(T+1)} \frac{EP(U_1^2+nU_3/\bar{\epsilon}^2)}{\alpha} \underbrace{\sum_{t=1}^T \frac{t}{t+2}}_{ =  T - \sum_{t=3}^{T+2} \frac{2}{t} \leq T} \nonumber \\
\leq &  \frac{2 U_2 \gamma + U_2^2 \rho^{\text{max}} + 4 \gamma^2 / \rho^1 + \alpha U_2^2/(2E) + 2P(U_1^2+nU_3/\bar{\epsilon}^2 )/\alpha }{T+1}.
\end{align*} 
This completes the proof.

\section{Distributed control of power flow} \label{apx:distpf}
We depict a power network by a graph $(\mathcal{N}, \mathcal{L})$, where $\mathcal{N}$ is a set of buses and
$\mathcal{L}$ is a set of lines. For every line $\ell_{ij}
\in \mathcal{L}$, where $i$ is a \textit{from bus}  and $j$
is a \textit{to bus} of line $\ell$, we are given line
parameters, including bounds $[\underline{\theta}_{ij},
\overline{\theta}_{ij}]$ on voltage angle difference,
thermal limit $\overline{s}_{\ell}$, resistance $r_{\ell}$,
reactance $x_{\ell}$, impedance $z_{\ell} := r_{\ell} +
\textbf{i} x_{\ell}$, line charging susceptance
$b^{\text{\tiny c}}_{\ell}$, tap ratio $\tau_{\ell}$, phase
shift angle $\theta^{\text{\tiny s}}_{\ell}$, and admittance
matrix $Y_{\ell}$: 
\begin{align*}
Y_{\ell}
:= &
\begin{bmatrix}
Y^{{\text{\tiny ff}}}_{\ell} & Y^{{\text{\tiny ft}}}_{\ell} \\
Y^{{\text{\tiny tf}}}_{\ell} & Y^{{\text{\tiny tt}}}_{\ell}
\end{bmatrix}
=
\begin{bmatrix}
  (z^{-1}_{\ell} + \textbf{i} \frac{b^{\text{\tiny c}}_{\ell}}{2}) \frac{1}{\tau_{\ell}^2}
& -z^{-1}_{\ell} \frac{1}{\tau_{\ell} e^{- \textbf{i} \theta^{\text{\tiny s}}_{\ell}}}  \\
-z^{-1}_{\ell} \frac{1}{\tau_{\ell} e^{ \textbf{i} \theta^{\text{\tiny s}}_{\ell}}}
& z^{-1}_{\ell} + \textbf{i} \frac{b^{\text{\tiny c}}_{\ell}}{2}
\end{bmatrix}
,
\end{align*}
\noindent
$G^{\text{\tiny cf}}_{\ell} := \Re (Y^{\text{\tiny
ff}}_{\ell})$, $B^{\text{\tiny cf}}_{\ell} := \Im
(Y^{\text{\tiny ff}}_{\ell})$, $G^{\text{\tiny f}}_{\ell} :=
\Re (Y^{\text{\tiny ft}}_{\ell})$, $B^{\text{\tiny
f}}_{\ell} := \Im (Y^{\text{\tiny ft}}_{\ell})$,
$G^{\text{\tiny ct}}_{\ell} := \Re (Y^{\text{\tiny
tt}}_{\ell})$, $B^{\text{\tiny ct}}_{\ell} := \Im
(Y^{\text{\tiny tt}}_{\ell})$, $G^{\text{\tiny t}}_{\ell} =
\Re (Y^{\text{\tiny tf}}_{\ell})$, and $B^{\text{\tiny
t}}_{\ell} := \Im (Y^{\text{\tiny tf}}_{\ell})$. 
For every bus $i \in \mathcal{N}$, we are given bus parameters,
including bounds $[\underline{v}_i, \overline{v}_i]$ on
voltage magnitude, active (resp., reactive) power demand
$p^{\text{\tiny d}}_i$ (resp., $q^{\text{\tiny d}}_i$),
shunt conductance $g^{\text{\tiny s}}_i$, and shunt
susceptance $b^{\text{\tiny s}}_i$. Furthermore, for every
$i \in \mathcal{N}$, we define subsets
$\mathcal{L}^{\text{\tiny F}}_i := \{ \ell_{ij} : j \in
\mathcal{N}, \ell_{ij} \in \mathcal{L} \} $ and
$\mathcal{L}^{\text{\tiny T}}_i := \{ \ell_{ji} : j \in
\mathcal{N}, \ell_{ji} \in \mathcal{L} \}$ of $\mathcal{L}$
and a set of generators $\mathcal{G}_i$. For every generator
$g \in \mathcal{G}_i$, we are given generator parameters,
including bounds $[\underline{p}^{\text{\tiny G}}_g,
\overline{p}^{\text{\tiny G}}_g]$ (resp.,
$[\underline{q}^{\text{\tiny G}}_g,
\overline{q}^{\text{\tiny G}}_g]$) on the amounts of active
(resp., reactive) power generation and coefficients
($c_{1,g}$, $c_{2,g}$) of the quadratic generation cost
function.

Next we present decision variables. 
For every line $\ell_{ij} \in \mathcal{L}$, we
denote active (resp., reactive) power flow along line $\ell$
by $p^{\text{\tiny F}}_{\ell}$, $p^{\text{\tiny T}}_{\ell}$
(resp., $q^{\text{\tiny F}}_{\ell}$, $q^{\text{\tiny
T}}_{\ell}$).
For every $i \in \mathcal{N}$, we denote the complex voltage by
$V_i  = v^{\text{\tiny R}}_i + \textbf{i} v^{\text{\tiny
I}}_i$, and we introduce the following auxiliary variables:
\begin{align}
& w^{\text{\tiny RR}}_{ij} = v^{\text{\tiny R}}_i v^{\text{\tiny R}}_j, \ \
w^{\text{\tiny II}}_{ij} = v^{\text{\tiny I}}_i v^{\text{\tiny I}}_j, \ \
w^{\text{\tiny RI}}_{ij} = v^{\text{\tiny R}}_i v^{\text{\tiny I}}_j, \ \forall j \in \mathcal{N}. \label{AC_Linking}
\end{align}
For every generator $g \in \mathcal{G}_i$, we
denote the amounts of active (resp., reactive) power
generation by $p^{\text{\tiny G}}_g$ (resp., $q^{\text{\tiny
G}}_g$).

Given the aforementioned parameters and decision variables, one can formulate the following convex program to determine power flow that minimizes the amounts of load shedding:
\begin{subequations}
\label{SOCOPF-rect}
\begin{align}
\min \ &  \sum_{i \in \mathcal{N}}  (s^{\text{\tiny p}}_i)^2 + (s^{\text{\tiny q}}_i)^2  \label{SOCOPF-rect-0} \\ 
\mbox{s.t.} \
&  \sum_{\ell \in \mathcal{L}^{\text{\tiny F}}_i} p^{\text{\tiny F}}_{\ell} + \sum_{\ell \in \mathcal{L}^{\text{\tiny T}}_i}  p^{\text{\tiny T}}_{\ell} - \sum_{g \in \mathcal{G}_i} p^{\text{\tiny G}}_g + p^{\text{\tiny d}}_i + g^{\text{\tiny s}}_i (w^{\text{\tiny RR}}_{ii}+w^{\text{\tiny II}}_{ii}) = s^{\text{\tiny p}}_i , \ \forall i \in \mathcal{N} , \label{SOCOPF-rect-7} \\
&  \sum_{\ell \in \mathcal{L}^{\text{\tiny F}}_i} q^{\text{\tiny F}}_{\ell} + \sum_{\ell \in \mathcal{L}^{\text{\tiny T}}_i}  q^{\text{\tiny T}}_{\ell} - \sum_{g \in \mathcal{G}_i} q^{\text{\tiny G}}_g + q^{\text{\tiny d}}_i - b^{\text{\tiny s}}_i (w^{\text{\tiny RR}}_{ii}+w^{\text{\tiny II}}_{ii}) = s^{\text{\tiny q}}_i , \ \forall i \in \mathcal{N}, \label{SOCOPF-rect-8} \\
&  w^{\text{\tiny RR}}_{ii} + w^{\text{\tiny II}}_{ii}  \in [\underline{v}_i^2, \overline{v}_i^2],  \ \forall i \in \mathcal{N}, \label{SOCOPF-rect-9} \\
& p^{\text{\tiny G}}_g \in [\underline{p}^{\text{\tiny G}}_g, \overline{p}^{\text{\tiny G}}_g], \ \ \ q^{\text{\tiny G}}_g \in [\underline{q}^{\text{\tiny G}}_g, \overline{q}^{\text{\tiny G}}_g], \forall i \in \mathcal{N}, \forall g \in \mathcal{G}_i, \label{SOCOPF-rect-10} \\
& p^{\text{\tiny F}}_{\ell} = G^{\text{\tiny f}}_{\ell} (w^{\text{\tiny RR}}_{ij} + w^{\text{\tiny II}}_{ij}) +  B^{\text{\tiny f}}_{\ell} (w^{\text{\tiny RI}}_{ji} + w^{\text{\tiny RI}}_{ij}) + G^{\text{\tiny cf}}_{\ell} ( w^{\text{\tiny RR}}_{ii}+ w^{\text{\tiny II}}_{ii} ), \ \forall \ell_{ij} \in \mathcal{L}, \label{SOCOPF-rect-1} \\
& q^{\text{\tiny F}}_{\ell} = G^{\text{\tiny f}}_{\ell} (w^{\text{\tiny RI}}_{ji} + w^{\text{\tiny RI}}_{ij}) - B^{\text{\tiny f}}_{\ell} (w^{\text{\tiny RR}}_{ij} + w^{\text{\tiny II}}_{ij}) -B^{\text{\tiny cf}}_{\ell} ( w^{\text{\tiny RR}}_{ii}+ w^{\text{\tiny II}}_{ii} ), \ \forall \ell_{ij} \in \mathcal{L},  \label{SOCOPF-rect-2}\\
& p^{\text{\tiny T}}_{\ell} = G^{\text{\tiny t}}_{\ell} (w^{\text{\tiny RR}}_{ji} + w^{\text{\tiny II}}_{ji}) +  B^{\text{\tiny t}}_{\ell} (w^{\text{\tiny RI}}_{ij} + w^{\text{\tiny RI}}_{ji}) + G^{\text{\tiny ct}}_{\ell} ( w^{\text{\tiny RR}}_{jj}+ w^{\text{\tiny II}}_{jj} ), \ \forall \ell_{ij} \in \mathcal{L}, \label{SOCOPF-rect-3} \\
& q^{\text{\tiny T}}_{\ell} = G^{\text{\tiny t}}_{\ell} (w^{\text{\tiny RI}}_{ij} + w^{\text{\tiny RI}}_{ji}) - B^{\text{\tiny t}}_{\ell} (w^{\text{\tiny RR}}_{ji} + w^{\text{\tiny II}}_{ji}) -B^{\text{\tiny ct}}_{\ell} ( w^{\text{\tiny RR}}_{jj}+ w^{\text{\tiny II}}_{jj} ), \ \forall \ell_{ij} \in \mathcal{L}, \label{SOCOPF-rect-4} \\
& (p^{\text{\tiny F}}_{\ell})^2 + (q^{\text{\tiny F}}_{\ell})^2 \leq (\overline{s}_{\ell})^2, \ \ \ (p^{\text{\tiny T}}_{\ell})^2 + (q^{\text{\tiny T}}_{\ell})^2 \leq (\overline{s}_{\ell})^2, \ \forall \ell_{ij} \in \mathcal{L},  \label{SOCOPF-rect-5} \\
& \ w^{\text{\tiny RI}}_{ji} - w^{\text{\tiny RI}}_{ij} \in [\tan(\underline{\theta}_{ij})(w^{\text{\tiny RR}}_{ij} + w^{\text{\tiny II}}_{ij}), \tan(\overline{\theta}_{ij})(w^{\text{\tiny RR}}_{ij} + w^{\text{\tiny II}}_{ij})], \ \forall \ell_{ij} \in \mathcal{L},  \label{SOCOPF-rect-6} \\
& (w^{\text{\tiny RR}}_{ij}+w^{\text{\tiny II}}_{ij})^2 + (w^{\text{\tiny RI}}_{ji}-w^{\text{\tiny RI}}_{ij})^2 + \Big( \frac{w^{\text{\tiny RR}}_{ii}+w^{\text{\tiny II}}_{ii} - w^{\text{\tiny RR}}_{jj}-w^{\text{\tiny II}}_{jj} }{2} \Big)^2   \label{SOCOPF-rect-11} \\
& \ \leq \Big( \frac{w^{\text{\tiny RR}}_{ii}+w^{\text{\tiny II}}_{ii} + w^{\text{\tiny RR}}_{jj}+w^{\text{\tiny II}}_{jj} }{2} \Big)^2, \forall \ell_{ij} \in \mathcal{L}, \nonumber
\end{align}
\end{subequations} 
where \eqref{SOCOPF-rect-0} is the objective function that minimizes the amounts of load shedding $s^{\text{\tiny p}}_i$ and $s^{\text{\tiny q}}_i$ that are computed by the power balance equation as in \eqref{SOCOPF-rect-7}--\eqref{SOCOPF-rect-8},
\eqref{SOCOPF-rect-9} represent bounds on voltage magnitudes,
\eqref{SOCOPF-rect-10} represent bounds on power generation, 
\eqref{SOCOPF-rect-1}--\eqref{SOCOPF-rect-4} represent power flow, 
\eqref{SOCOPF-rect-5} represent line thermal limit, 
\eqref{SOCOPF-rect-6} represent bounds on voltage angle differences, and 
\eqref{SOCOPF-rect-11} represent SOC constraints that ensure linking between auxiliary variables.

For distributed control of power flow, we consider that the network $(\mathcal{N}, \mathcal{L})$ is decomposed into several zones indexed by $\mathcal{Z}:= \{1, \ldots, Z\}$. Specifically, we split a set $\mathcal{N}$ of buses into subsets $\{ \mathcal{N}_z\}_{z \in \mathcal{Z}}$ such
that $\mathcal{N} = \cup_{z \in \mathcal{Z}} \mathcal{N}_z$ and $\mathcal{N}_z \cap \mathcal{N}_{z'} = \emptyset$ for $z,z' \in \mathcal{Z}: z \neq z'$. For each zone $z \in \mathcal{Z}$ we define a line set $\mathcal{L}_z:= \cup_{i \in \mathcal{N}_z} \big( \mathcal{L}^{\text{\tiny F}}_i \cup \mathcal{L}^{\text{\tiny T}}_i \big)$; an extended node set $\mathcal{V}_z := \cup_{i \in \mathcal{N}_z} \mathcal{A}_i$, where $\mathcal{A}_i$ is a set of adjacent buses of $i$; and
a set of cuts $\mathcal{C}_z = \cup_{z' \in \mathcal{Z}\setminus \{z\} } (\mathcal{L}_z \cap
\mathcal{L}_{z'})$. Note that $\{ \mathcal{N}_z \}_{z \in \mathcal{Z}}$ is a collection of disjoint sets, while $\{\mathcal{L}_z \}_{z \in \mathcal{Z}}$ and $\{ \mathcal{V}_z\}_{z \in \mathcal{Z}}$ are not. 
Using these notations, we rewrite problem \eqref{SOCOPF-rect} as
\begin{subequations}
\label{model:ACOPF_matrix}
\begin{align}
\min \ & \sum_{z \in \mathcal{Z}} f_z(x_z) \\
\mbox{s.t.} \ & (x_z, y_z) \in \mathcal{F}_z, \ \forall z \in \mathcal{Z}, \label{model:ACOPF_matrix-1} \\
& \phi_i = y_{zi}, \ \forall z \in \mathcal{Z}, \forall i \in C(z), \label{model:ACOPF_matrix-2} \\
& \phi_i \in \mathbb{R}, \ \forall i \in \mathcal{C}, \label{model:ACOPF_matrix-3}
\end{align}
\end{subequations}
where
\begin{align*}
  x_z \leftarrow & \big\{p^{\text{\tiny F}}_{z\ell}, q^{\text{\tiny F}}_{z\ell}, p^{\text{\tiny T}}_{z\ell}, q^{\text{\tiny T}}_{z\ell},
  w^{\text{\tiny RR}}_{zij}, w^{\text{\tiny II}}_{zij}, w^{\text{\tiny RI}}_{zij}, w^{\text{\tiny RI}}_{zji} \big\}_{\ell_{ij} \in \mathcal{L}_z \setminus \mathcal{C}_z} \\
  &  \cup \big\{ s^{\text{\tiny p}}_{i}, s^{\text{\tiny q}}_{i} \big\}_{i \in \mathcal{V}_z} \cup \big\{ p^{\text{\tiny G}}_g, q^{\text{\tiny G}}_g \big\}_{i \in \mathcal{N}_z, g \in \mathcal{G}_i}, \\  
y_z \leftarrow &  \big\{p^{\text{\tiny F}}_{z\ell}, q^{\text{\tiny F}}_{z\ell}, p^{\text{\tiny T}}_{z\ell}, q^{\text{\tiny T}}_{z\ell},
w^{\text{\tiny RR}}_{zij}, w^{\text{\tiny II}}_{zij}, w^{\text{\tiny RI}}_{zij}, w^{\text{\tiny RI}}_{zji} \big\}_{\ell_{ij} \in \mathcal{C}_z}, \\
\phi \leftarrow & \cup_{z \in \mathcal{Z}} \big\{ p^{\text{\tiny F}}_{\ell}, q^{\text{\tiny F}}_{\ell}, p^{\text{\tiny T}}_{\ell}, q^{\text{\tiny T}}_{\ell},
w^{\text{\tiny RR}}_{ij}, w^{\text{\tiny II}}_{ij}, w^{\text{\tiny RI}}_{ij}, w^{\text{\tiny RI}}_{ji} \big\}_{\ell_{ij} \in \mathcal{C}_z}, 
\end{align*}
$C(z)$ is an index set that indicates each element of $y_z$,
$\mathcal{C} := \cup_{z \in \mathcal{Z}} C(z)$ is an index set of consensus variable {\color{black} $\phi$},  $\mathcal{F}_z := \{(x_z, y_z): \eqref{SOCOPF-rect-1}-\eqref{SOCOPF-rect-6}, \forall \ell_{ij} \in \mathcal{L}_z; \ \eqref{SOCOPF-rect-7}, \eqref{SOCOPF-rect-8}, \forall i \in \mathcal{N}_z; \ \eqref{SOCOPF-rect-9},\forall i \in \mathcal{V}_z; \ \eqref{SOCOPF-rect-10},\forall i \in \mathcal{N}_z,\forall g \in \mathcal{G}_i; and \ \eqref{SOCOPF-rect-11}, \forall i \in \mathcal{V}_z,\forall j\in \mathcal{V}_z \} $ is a convex feasible region defined for each zone.
We note that \eqref{model:ACOPF_matrix} is the form of \eqref{model:dist_1}.

\section{Hyperparameter Tuning} \label{apx:hyper}
The parameter $\rho^t$ in Assumption \ref{assump:convergence} affects the learning performance because it controls the proximity of the local model parameters from the global model parameters.
For all algorithms, we set $\rho^t \leftarrow \hat{\rho}^t$ given by
\begin{align}
  & \hat{\rho}^t := \min \{ 10^9, \ c_1 (1.2)^{\lfloor t/T_c \rfloor} + c_2/\bar{\epsilon}\}, \ \forall t \in [T], \label{dynamic_rho}
\end{align}
where (i) $c_1=2$, $c_2=5$, and $T_c=10000$ for MNIST and (ii) $c_1=0.005$, $c_2=0.05$, and $T_c=2000$ for FEMNIST. 
Note that the chosen parameter $\hat{\rho}^t$ is nondecreasing and bounded above, thus satisfying Assumption \ref{assump:convergence}-(A).
 
\begin{figure}[!htt]
  \centering
  \begin{subfigure}[b]{0.24\textwidth}
      \centering
      \includegraphics[width=\textwidth]{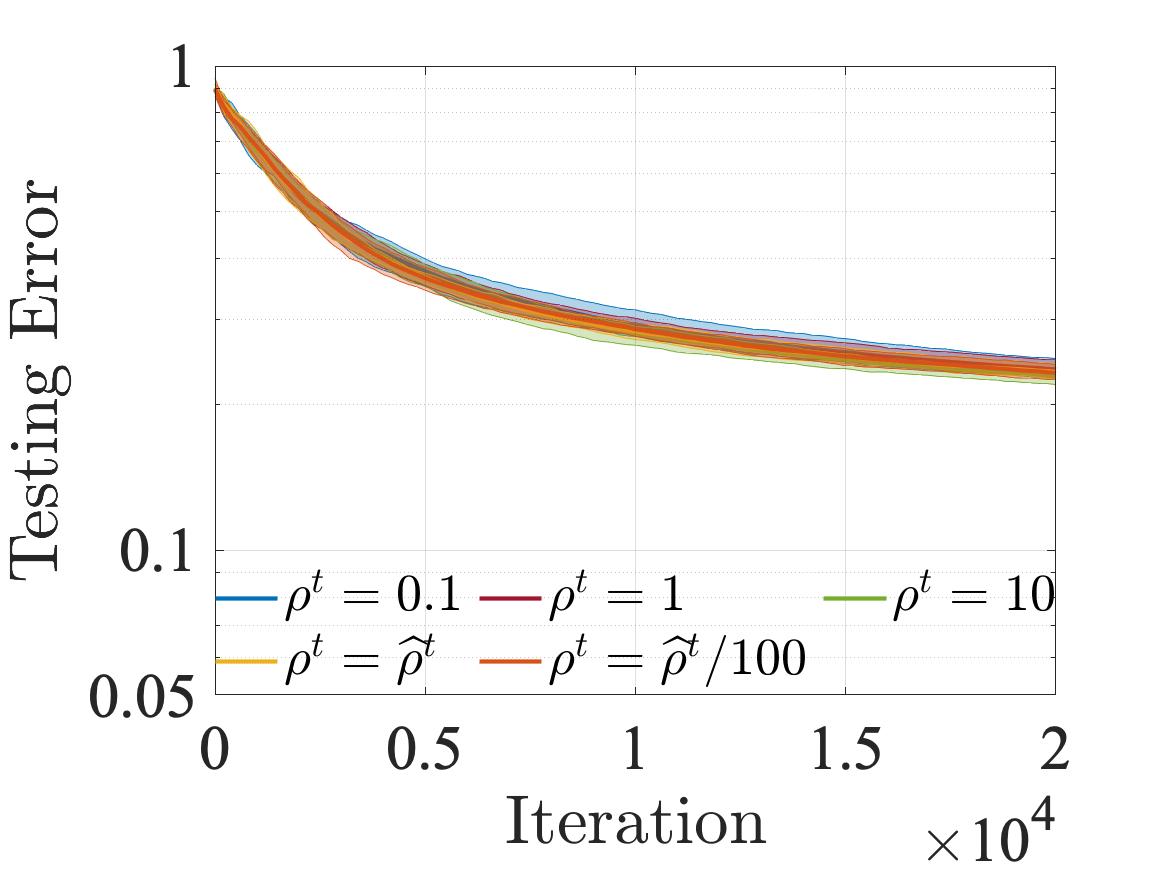}
      \includegraphics[width=\textwidth]{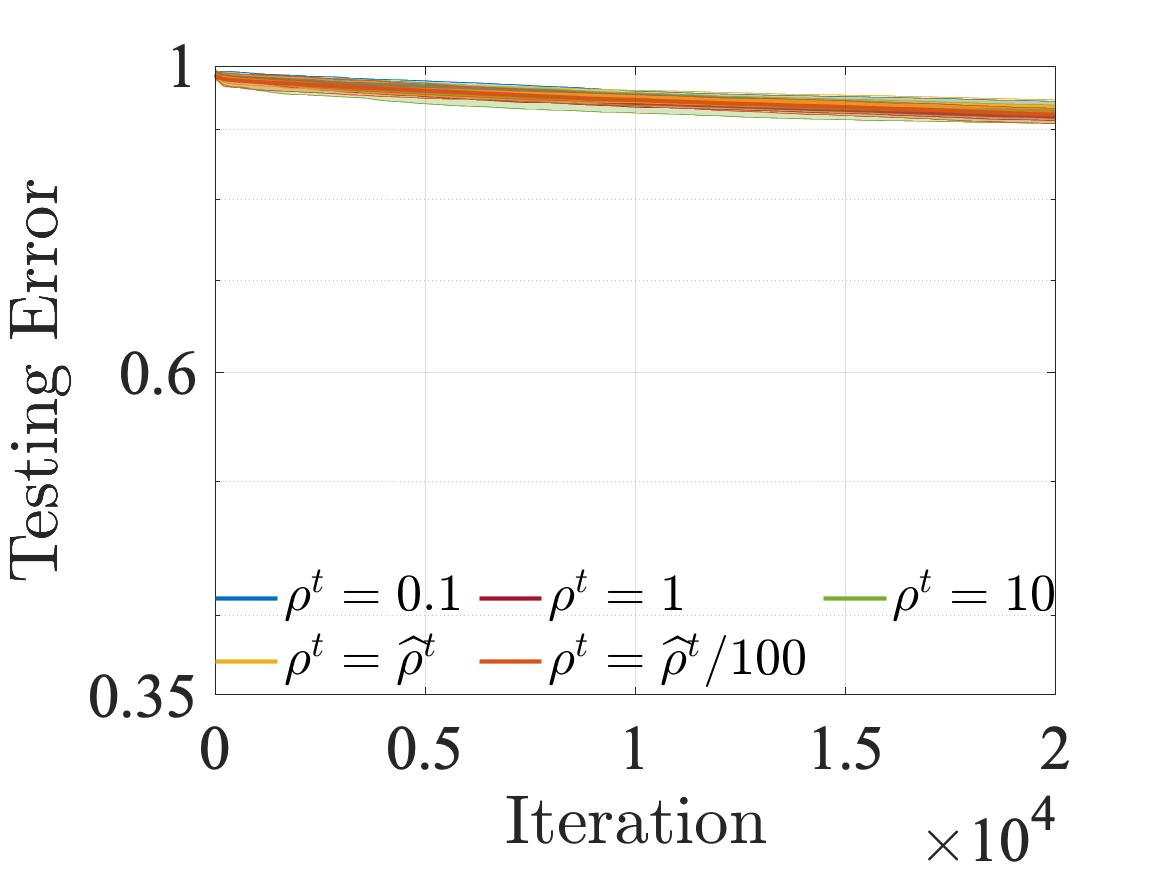}
      \caption{$\bar{\epsilon}=0.05$}
  \end{subfigure}
  \begin{subfigure}[b]{0.24\textwidth}
    \centering
    \includegraphics[width=\textwidth]{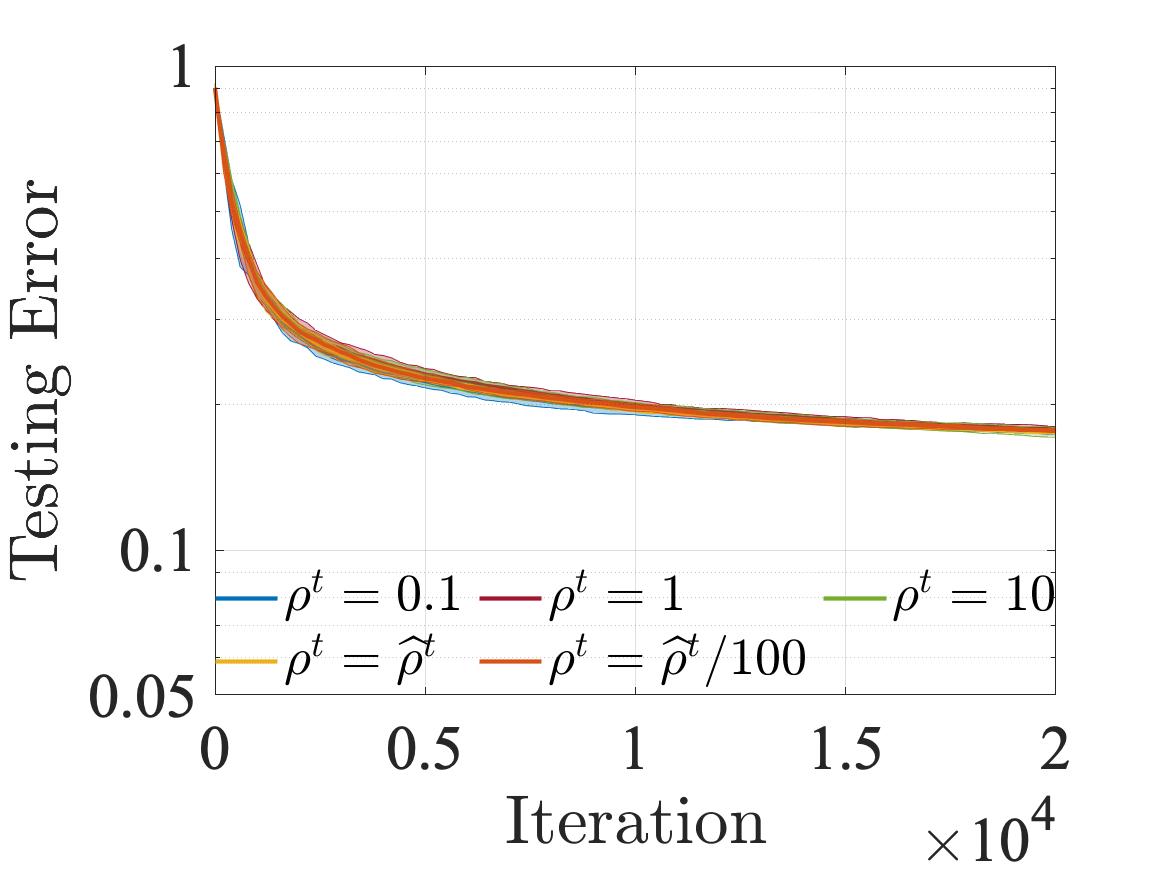}
    \includegraphics[width=\textwidth]{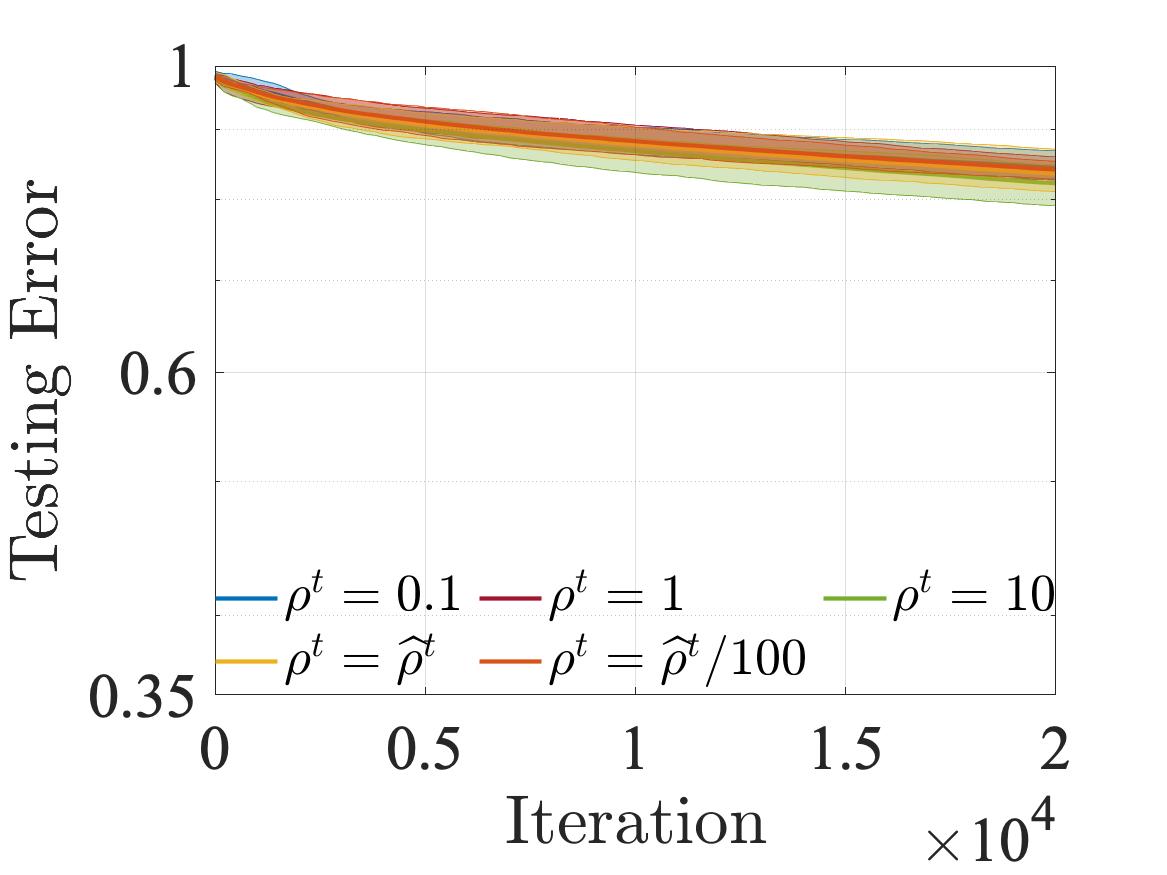}
    \caption{$\bar{\epsilon}=0.1$}
\end{subfigure}
  \begin{subfigure}[b]{0.24\textwidth}
      \centering
      \includegraphics[width=\textwidth]{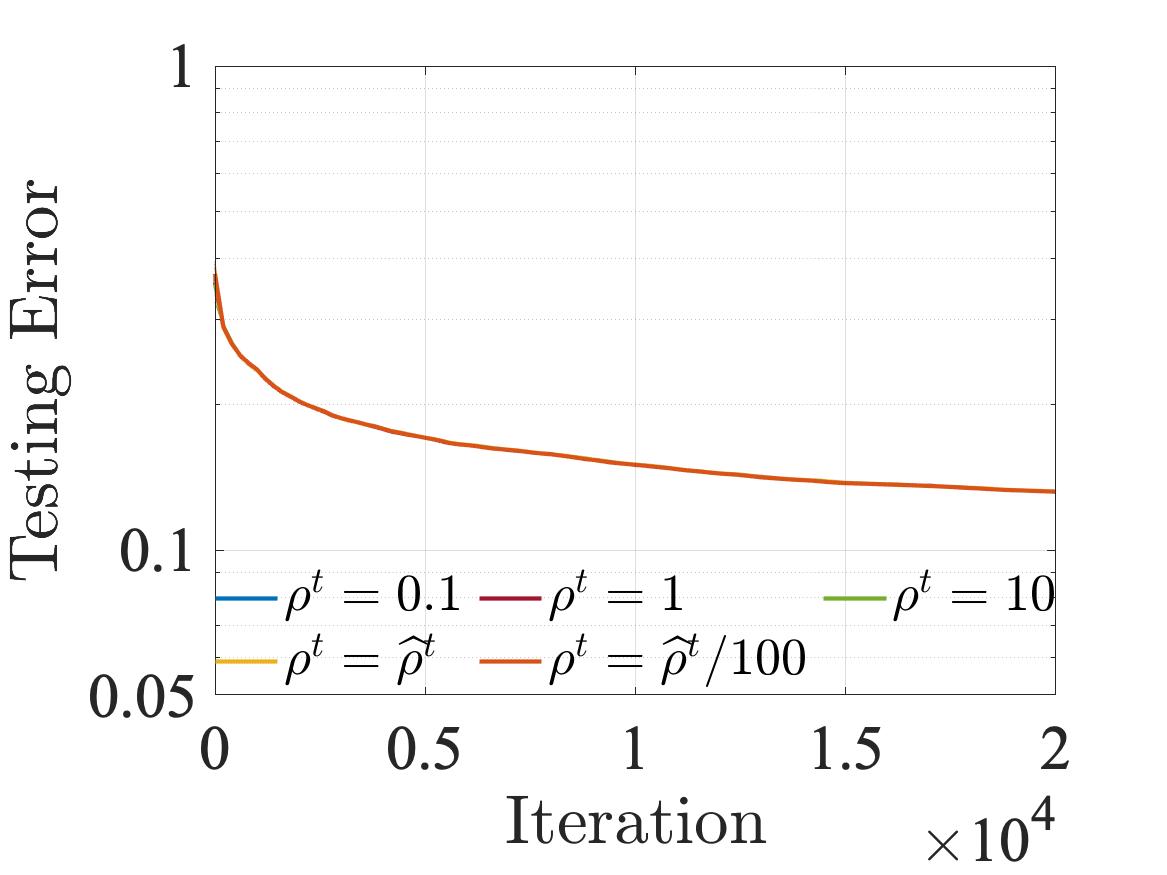}
      \includegraphics[width=\textwidth]{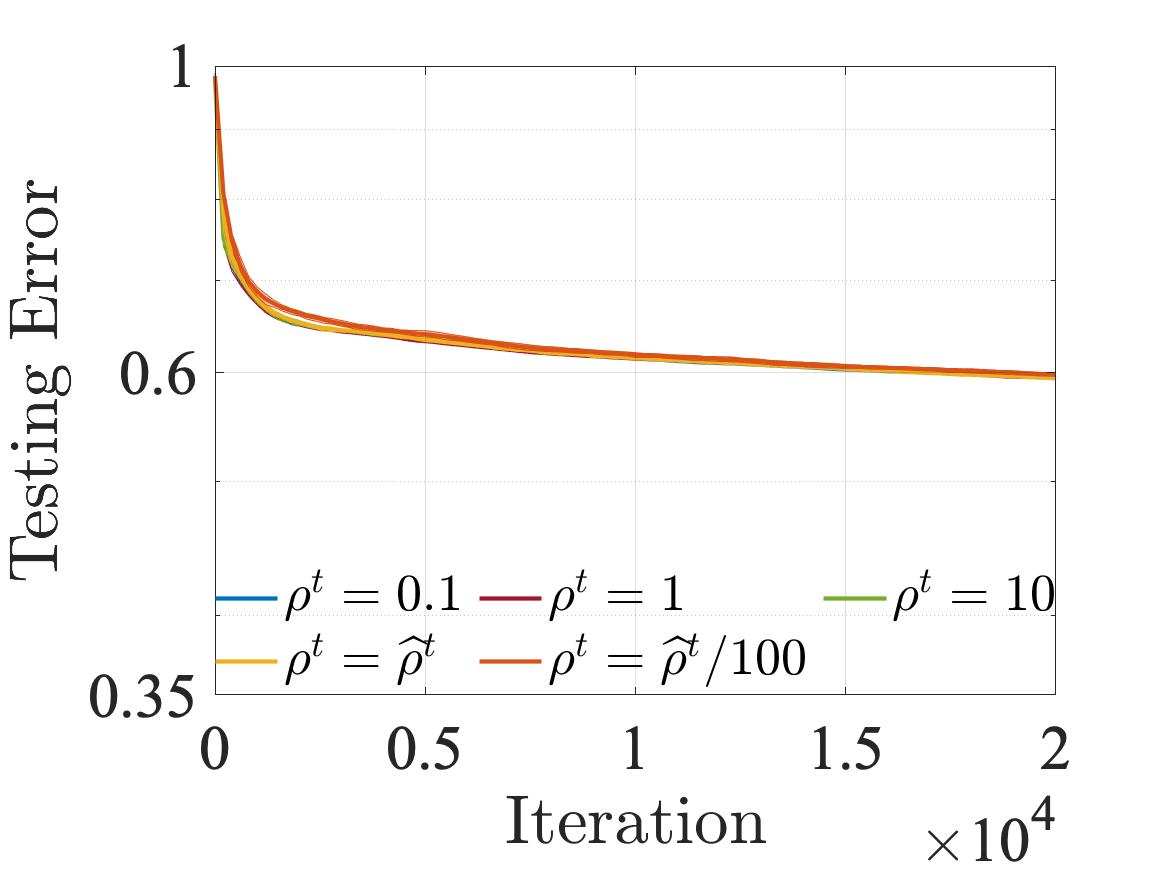}
      \caption{$\bar{\epsilon}=1$}
  \end{subfigure}
  \begin{subfigure}[b]{0.24\textwidth}
    \centering
    \includegraphics[width=\textwidth]{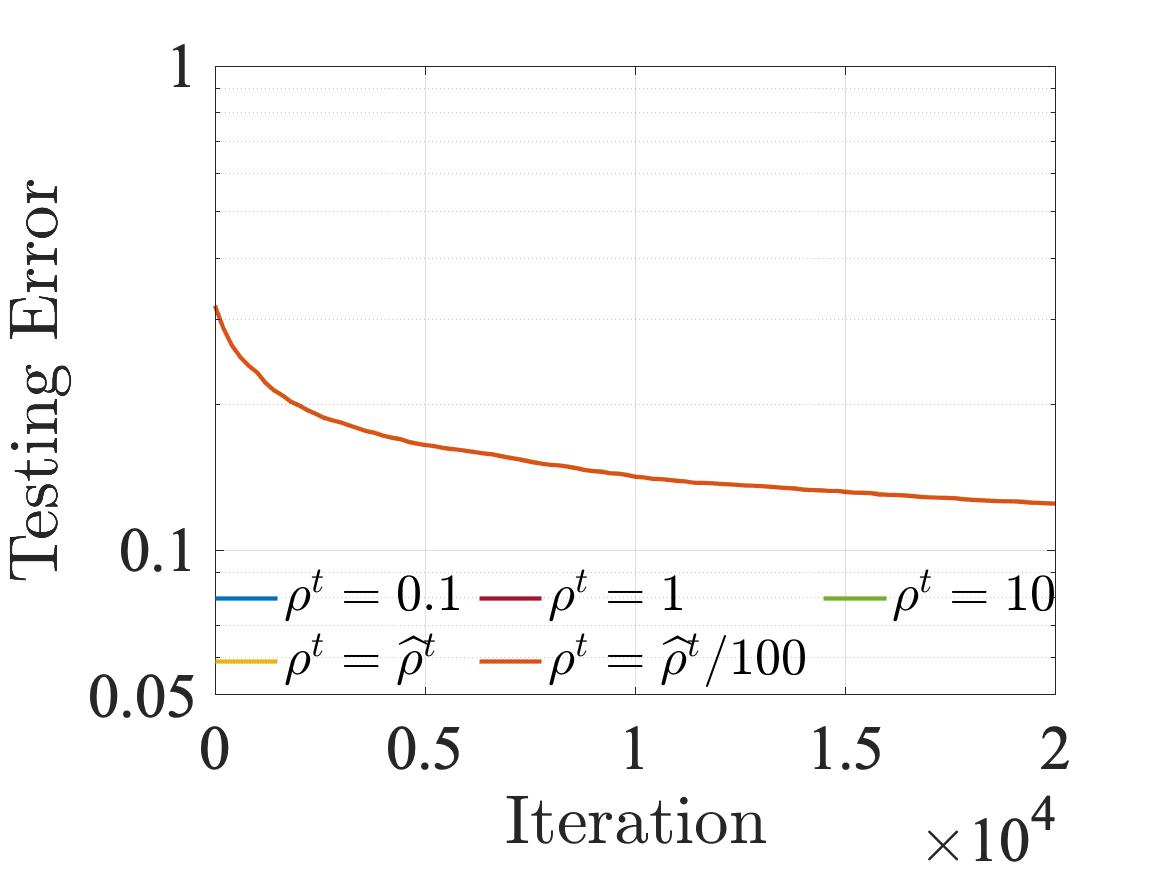}
    \includegraphics[width=\textwidth]{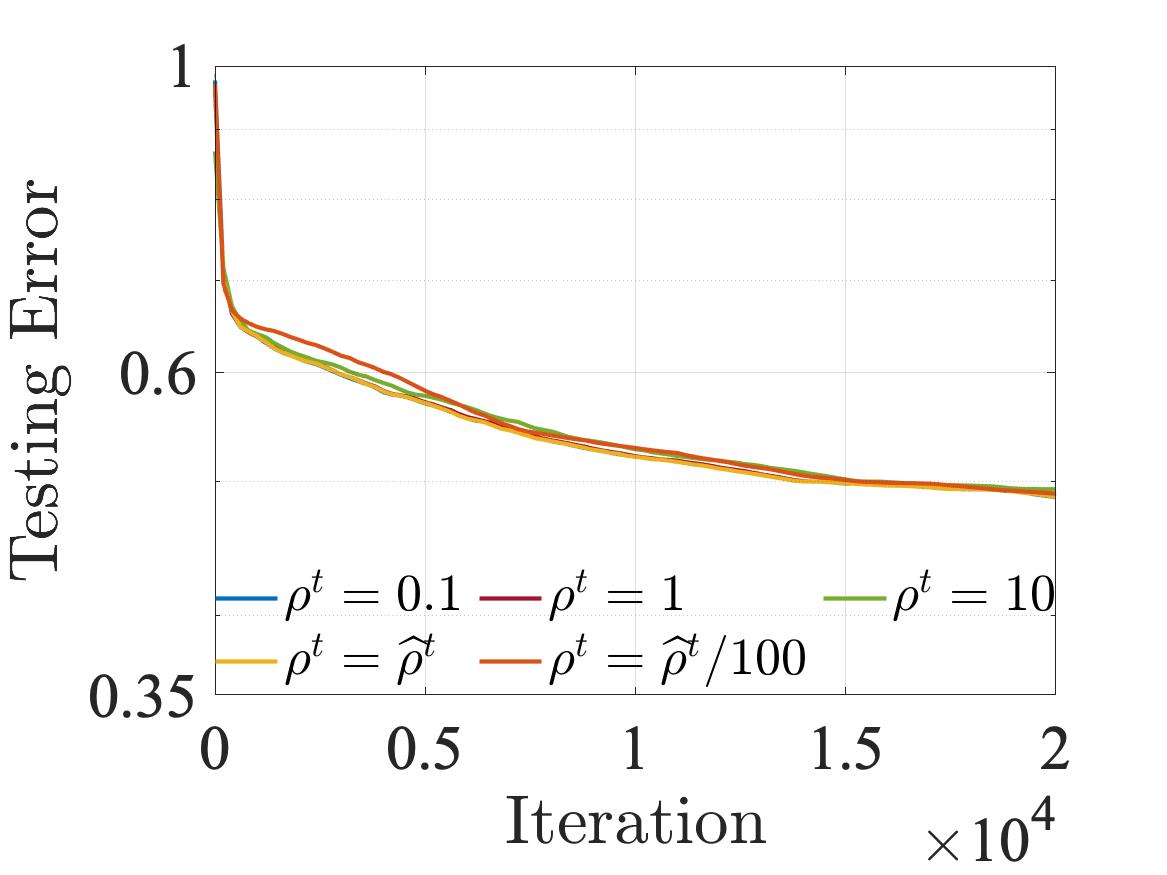}
    \caption{$\bar{\epsilon}=5$}
  \end{subfigure}
     \caption{Testing errors of \texttt{OutP} using MNIST (top) and FEMNIST (bottom).}
     \label{fig:hyper_rho}
\end{figure}

Since these parameter settings may not lead \texttt{OutP} to its best performance, we test various $\rho^t$ for \texttt{OutP} using a set of static parameters, $\rho^t \in \{0.1, 1, 10\}$ for all $t \in [T]$, where $\rho^t=0.1$ is chosen in \cite{huang2019dp}, and dynamic parameters $\rho^t \in \{\hat{\rho}^t, \hat{\rho}^t/100 \}$, where $\hat{\rho}^t$ is from \eqref{dynamic_rho}.
In Figure \ref{fig:hyper_rho} we report the testing errors of \texttt{OutP} using MNIST and FEMNIST under various $\rho^t$ and $\bar{\epsilon}$.
The results imply that the performance of \texttt{OutP} is not greatly affected by the choice of $\rho^t$, but $\bar{\epsilon}$.
Hence, for all algorithms, we use $\hat{\rho}^t$ in \eqref{dynamic_rho}.

\section*{Acknowledgments}
This work was supported by the U.S. Department of Energy, Office of Science, Advanced Scientific Computing Research, under Contract DE-AC02-06CH11357.
We gratefully acknowledge the computing resources provided on Bebop and Swing, high-performance computing clusters operated by the Laboratory Computing Resource Center at Argonne National Laboratory.

\bibliographystyle{siamplain}
\bibliography{references}

\vspace{0.5in}
\noindent\fbox{\parbox{0.97\textwidth}{
The submitted manuscript has been created by UChicago Argonne, LLC, Operator of Argonne National Laboratory (``Argonne''). Argonne, a U.S. Department of Energy Office of Science laboratory, is operated under Contract No. DE-AC02-06CH11357. The U.S. Government retains for itself, and others acting on its behalf, a paid-up nonexclusive, irrevocable worldwide license in said article to reproduce, prepare derivative works, distribute copies to the public, and perform publicly and display publicly, by or on behalf of the Government. The Department of Energy will provide public access to these results of federally sponsored research in accordance with the DOE Public Access Plan (http://energy.gov/downloads/doe-public-access-plan).}
}

\end{document}